\documentclass[10pt]{article}
\usepackage{amssymb,amsmath,latexsym}
\def\Section#1{ \setcounter{equation}{0} \section{#1}}

\newcommand{\preuve}[1][\!\!]{\bigskip\noindent{\bf Proof #1. \ \ }}
\def\ititem#1{\bigskip \par \noindent \it #1 \rm}
\def\fin{\hfill$\Box$\\}
\def\bb{{\cal B}}\def\cc{{\cal C}}

\def\oo{{\cal O}}
\def\ss{{\cal S}}

\def\R{\mathbb R}\def\C{\mathbb C}\def\N{\mathbb N}

\def\D{\partial}
\def\phi{\varphi}
\def\norm#1{\left\Vert#1\right\Vert}

\def\abs#1{\left\vert#1\right\vert}
\def\set#1{\left\{#1\right\}}
\def\seq#1{\left<#1\right>}
\def\sep#1{\left(#1\right)}

\def\defegal{\stackrel{\text{\rm def}}{=}}

\newcommand{\half}{\mbox{$\frac{1}{2}$}}

\newcommand{\tr}{\operatorname{tr}}

\newcommand{\sgn}{\operatorname{sgn}}

\def\ekv#1#2{\begeq\label{#1}#2\endeq}

\def\iint{\int\hskip -2mm\int}
\def\iiint{\int\hskip -2mm\int\hskip -2mm\int}

\def\3{\vert \hskip -1pt\vert\hskip -1pt\vert }

\def\ably{arbitrarily}

\def\asy{asymptotic}
\def\bdd{bounded}
\def\bdy{boundary}

\def\ctf{canonical transformation}
\def\diffeo{diffeomorphism}

\def\ev{eigenvalue}
\def\e{equation}
\def\fu{function}

\def\fop{Fourier integral operator}

\def\hol{holomorphic}
\def\indep{independent}

\def\mfld{manifold}

\def\neigh{neighborhood}
\def\nondeg{non-degenerate}
\def\op{operator}
\def\og{orthogonal}
\def\pb{problem}

\def\pert{perturbation}
\def\plsh{plurisubharmonic}
\def\pol{polynomial}

\def\pop{pseudodifferential operator}

\def\rhs{right hand side}
\def\sa{selfadjoint}

\def\sufly{sufficiently}
\def\tf{transformation}

\def\tf{transform}

\def\trans{^t\hskip -2pt}

\def\ufly{uniformly}
\def\vf{vector field}
\def\wrt{with respect to}

\def\Re{{\mathrm Re\,}}
\def\Im{{\mathrm Im\,}}
\def\dbar{\overline{\partial }}

\def\wrtext#1{\relax\ifmmode{\leavevmode\hbox{#1}}\else{#1}\fi}
\def\abs#1{\left|#1\right|}
\def\begeq{\begin{equation}}
\def\endeq{\end{equation}}

\newcommand{\eps}{\epsilon}
\def\part#1{\frac{\partial}{\partial #1}}
\def\half{\frac{1}{2}}

\def\norm#1{\Vert #1\Vert }

\renewcommand{\Re}{\mbox{\rm Re\,}}
\renewcommand{\Im}{\mbox{\rm Im\,}}

\renewcommand{\exp}{\mbox{\rm exp\,}}

\newtheorem{dref}{Definition}[section]
\newtheorem{lemma}[dref]{Lemma}
\newtheorem{theo}[dref]{Theorem}
\newtheorem{prop}[dref]{Proposition}
\newtheorem{remark}[dref]{Remark}

\renewcommand{\thedref}{\thesection.\arabic{dref}}

\renewcommand{\remark}{\refstepcounter{dref}\bigskip
\noindent\bf Remark \thedref\rm\ \ }

\def\ignore#1{}

\title{Semiclassical analysis for
the  Kramers-Fokker-Planck equation}
\author{Fr{\'e}d{\'e}ric H{\'e}rau\\
Laboratoire de Math{\'e}matiques \\ Universit{\'e} de Reims \\
Moulin de
la Housse B.P. 1039 \\ 51687 Reims cedex 2 \\ France \\
herau@univ-reims.fr \and Johannes Sj{\"o}strand\\Centre de
Math{\'e}matiques\\Ecole Polytechnique\\91128
Palaiseau C{\'e}dex\\France\\johannes@math.polytechnique.fr
\and
Christiaan C. Stolk\\
Centre de Math{\'e}matiques\\
Ecole Polytechnique\\
91128 Palaiseau C{\'e}dex\\
France\\
stolk@math.polytechnique.fr}
\date{}
\begin{document}
\maketitle

\begin{center}
\begin{minipage}{13cm}
\small
\begin{center}  \bf R{\'e}sum{\'e} \rm \end{center}
On {\'e}tudie des estimations semiclassiques sur la r{\'e}solvente
d'op{\'e}rateurs qui ne sont ni elliptiques ni autoadjoints, que l'on
utilise pour {\'e}tudier le probl{\`e}me de Cauchy. En particulier on
obtient une description pr{\'e}cise du spectre pres de l'axe
imaginaire, et des estimations de r{\'e}solvente {\`a} l'int{\'e}rieur du
pseudo-spectre. On applique ensuite les r{\'e}sultats {\`a} l'op{\'e}rateur de
Kramers-Fokker-Planck.

\begin{center}  \bf Abstract \rm \end{center}
  We study some accurate semiclassical resolvent estimates
  for operators that
are neither selfadjoint nor elliptic, and applications to the
Cauchy problem. In particular we get a precise description of the
spectrum near the imaginary axis and precise resolvent estimates
inside the pseudo-spectrum. We apply our results to the
Kramers-Fokker-Planck operator.
\end{minipage}
\end{center}


\newpage
\tableofcontents

\Section{Introduction}\label{section0} \setcounter{equation}{0}

In certain applications one is interested in the long-time
behavior of systems described by a linear partial differential
equation. For example in kinetic equations one studies the decay
to equilibrium of various linear and nonlinear systems. For the
Kramers-Fokker-Planck equation, that will be studied here,
exponential decay was shown in \cite {Tal99} and an explicit rate
was given in \cite{HerNi03}, following earlier results of
Desvillettes and Villani who established explicit decay of any
polynomial order in $t^{-1}$ \cite{DesvillettesVillani2001}. In
\cite{HerNi03,DesvillettesVillani2001}  more general discussions
on decay to equilibrium in kinetic equations can be found.

The study of the long-time behavior naturally leads one to
 study the
spectrum, and, for non-selfadjoint problems that we study
here, to study the growth of the resolvent. For the
Kramers-Fokker-Planck equation this is complicated by the fact
 that the
operator whose time evolution is to be computed is not elliptic, but
only satisfies certain subellipticity conditions. To deal with this
H{\'e}rau and Nier exploit the relation between the
Kramers-Fokker-Planck operator and certain Witten Laplacians. They
obtain estimates for the decay to equilibrium in terms of the first
eigenvalue of this Witten Laplacian. More on the connection between
the Kramers-Fokker-Planck equation and Witten Laplacians can be found
in \cite{HelNi}.

Resolvent estimates have been studied from a different perspective
by a group of authors interested in the notion of pseudospectrum,
i.e.\ the region in the complex spectral plane where the resolvent
may be large. In recent years there has been a great interest in
this area following work of Trefethen, Davies, Zworski and others
\cite{Trefethen1997,Davies1999a,Davies1999b,Zworski2001}. In
\cite{DenckerSjostrandZworski2004} the authors studied the
location of the spectrum inside the spectrum in the semiclassical
limit, and adapted subelliptic estimates to this situation. In
\cite{Hi04}, M. Hitrik has obtained related results for operators
in one dimension.

In the present work we apply such ideas to a class of
pseudodifferential operators that includes the Kramers-Fokker-Planck
operator. We obtain a number of higher eigenvalues, in the
semiclassical limit, for the original operator, i.e.\ not the Witten
Laplacian. We also obtain precise resolvent estimates. We use roughly
the same estimates as \cite{DenckerSjostrandZworski2004} in one region
of phase space, while in other regions we have to make important
changes, and additions. This is then applied for the time evolution.

Evolution problems have also attracted recent interest
\cite{Davies2003a,TangZworski1998,BurqZworski2001}, and here a
difficulty is the generally quite wild growth of the resolvent inside
the pseudospectrum. It is therefore of interest that for a concrete
physically interesting model, we are able to control the resolvent
sufficiently well to get quite precise results about the long-time
evolution.

Spectral properties for some different Fokker-Planck equations
(without the subellipticity property) have been discussed by
Kolokoltsov \cite{Kolokoltsov2000}. In probability theory many other
problems for equations with a small diffusive term have been studied,
see for example the monograph \cite{FreidlinWentzell1984}.

\medskip\smallskip

Our main example, the Kramers-Fokker-Planck operator is given by
\begin{equation} \label{defFP}
  P =v\cdot h\partial _x-V'(x)\cdot h\partial _v+ {\gamma \over
  2}(-(h\partial _v)^2+v^2-hn)
\end{equation}
on $\R ^{2n}$, where $V$ is a $\cc^\infty$ potential, $h$ is
essentially  the temperature, and $x$, $v \in \R^n$. The operator
$P$ is derived from the original equation, introduced in
one-dimensional form by Kramers \cite{Kramers1940}, in Section
$\ref{secFP}$ below (see also \cite{Ri89}, \cite{HerNi03}). The
time evolution problem is given by
\[
  (h \partial_t + P) u(t,x,v) = 0, \qquad u(0,\cdot,\cdot) = u_0 .
\]
As mentioned, we are interested in the low temperature limit
\[
  0<h\ll 1 ,
\]
and the equations are rescaled according to the standard convention in
semiclassical analysis, where each derivative comes with an $h$.
Our main result about the
Kramers-Fokker-Planck equation is the following theorem.

\begin{theo} \label{mainFP} Assume $V$ is a Morse function
and that outside a compact region,
  $|V'(x)| \geq c_0 >0$. Assume also that the derivatives of $V$
of  order 2 or more are bounded. Then there exist constants
$c,C'>0$ such that for every  $C>1$:
\begin{description}
\item{a)} For any fixed \neigh{} $\Omega$ of the \ev{}s of the quadratic
  approximation of ${P _{\vert}}_{h=1}$ at the critical points, 
there exist $h_0$, $C''>0$
  such that for $0<h\le h_0$, $|z|\le C$, $z\not\in \Omega $, 
$$
h\norm{u} \leq C'' \norm{(P -hz)u}, \ \ \ \ \forall u \in \ss.
$$
\item{b)}  There exists $h_1>0$, such that for $0
<h\leq h_1$, $\Re(z) \leq c|z|^{1/3}h^{2/3}$ and $|z|\ge Ch$,
$$
|z|^{1/3}h^{2/3}\norm{u} \leq C'\norm{(P -z)u}, \ \ \ \ \forall u
\in \ss.
$$
\end{description}
\end{theo}

In fact this theorem on the Kramers-Fokker-Planck operator is a
consequence of a more general one. Let us first write the
hypotheses that will be needed for the symbol $p$ of the more
general operator $p^w$ that we shall study.
We assume that $p=p_1+i p_2$
is a smooth function on $\R^{2n}_{x,\xi }$ with $p_1\ge 0$. (The
previous space $\R^{2n}_{x,v}$ now becomes $\R^n_{x}$.)

\bigskip
\bf Assumptions near the critical points: \rm Assume that $p$ has
finitely many critical points $\rho _1,\rho _2,...,\rho _N$ with
$p(\rho _j)=0$. Let $\delta (\rho )\ge 0$ be equivalent to the
distance from $\rho $ to ${\cal C}:=\{ \rho _1,\rho _2,...,\rho
_N\}$, with $\delta ^2\in C^\infty $. We assume in the following
that in a fixed open ball $\bb$ containing $\cc$ we have
\begin{equation} \label{hyppres}
\textbf{(H1)} \ \ \ \ \  p_1+\eps_0 H_{p_2}^2p_1\sim \delta ^2
\end{equation}
for a sufficiently small $\epsilon_0>0$. The assumption that
$p(\rho _j)=0$ is for simplicity only. As we shall later, this
implies that the critical points are \nondeg{}.

\bigskip
\bf Assumptions at infinity: \rm In the following we use the
notions of admissible metrics and weights in the sense of the
Weyl-H{\"o}rmander calculus, that we review in Section
\ref{Secweyl}. We first define an admissible metric on
$\R^{2n}_\rho$ with $\rho = (x,\xi)$:
 $$
 \Gamma_0 = dx^2 + \frac{d\xi^2}{\lambda^2},
 $$
where $\lambda = \lambda(\rho)$.  There is no restriction to
assume  that $1\le \lambda\in \cc^\infty$, and we suppose also that
 \begin{equation} \label{symblambda}
\textbf{(H2)} \ \ \ \ \   \lambda \in S(\lambda,\Gamma_0), \ \ \ \
\D \lambda \in S(1,\Gamma_0).
\end{equation}
 If $m$ is an admissible weight, recall that $S(m,\Gamma_0)$
 is the class of
$\cc^\infty$
 symbols  $p$ satisfying
 $\D^\alpha_x\D^\beta_\xi p (\rho)
 = \oo\sep{ m(\rho) \lambda(\rho)^{-|\beta|}}.$
  We suppose first that $p$ is a symbol of order $2$ but with the
first and second derivatives better than what would be given by
the symbolic calculus:
\begin{equation} \label{symbcal}
\begin{split}
\textbf{(H3)} &\ \ \ \ \ \
  p  \in S(\lambda^{2},\Gamma_0), \ \ \
  \D p  \in S(\lambda^{},\Gamma_0), \ \ \
  \D^2 p_1  \in S(1,\Gamma_0), \ \ \
 \D H_{p_2} p_1   \in S(\lambda,\Gamma_0).
\end{split}
\end{equation}
 We now assume that
outside any fixed neighborhood of $\cc$ we have the following gain
\begin{equation} \label{hyploin}
\textbf{(H4)} \ \ \ \ \  p_1+\eps_0 H_{p_2}^2p_1 \sim \lambda^2.
\end{equation}

\bigskip
Note that these assumptions are satisfied by the symbol of the
Kramers-Fokker-Planck
operator (see Section \ref{secFP}). In order to give a unique
assumption on
the whole space, we extend the function $\delta$ to $\R^{2n}$ to be a
smooth function on $\R^{2n}$, strictly positive away from ${\cal C}$,
and constant outside a fixed neighborhood of that set.
There  is no restriction to assume that
$\lambda = 1 $ inside the same neighborhood. Then
(\ref{hyppres}--\ref{hyploin}) can be summarized in the following
 way
\begin{equation} \label{mainhyp}
p_1+\epsilon_0 H_{p_2}^2p_1 \sim (\lambda\delta)^2.
\end{equation}

\bigskip
We have the following theorem for $P=p^w$:

\begin{theo}  \label{main}
Suppose $p$ satisfies  \bf (H1--H4)\it. Then  there exist
constants $c$, $C'>0$ such that for every $C\ge 1$:
\begin{description}
\item{a)} For any fixed \neigh{} $\Omega $ of the \ev{}s of the
  quadratic approximations of ${P _{\vert}}_{h=1}$ at the critical
  points, 
there exist $h_0,
  C''>0$ such that for $0<h\le h_0$,     $|z| \leq C $, $z\not\in
  \Omega $, 
$$
h\norm{u} \leq C'' \norm{(P -hz)u}, \ \ \ \ \forall u \in \ss.
$$
\item{b)}  There exists $h_1>0$, such that for $0
<h\leq h_1$, $\Re(z) \leq c|z|^{1/3}h^{2/3}$, $|z|\ge Ch$,
$$
|z|^{1/3}h^{2/3}\norm{u} \leq C'\norm{(P -z)u}, \ \ \ \ \forall u
\in \ss.
$$
\end{description}
\end{theo}

\par Here, if $\rho_0$ is a critical point of $p$, we define the
quadratic approximation $P_0$ of $P$, to be the $h=1$ quantization
of $\sum_{|\alpha +\beta |=2}{1\over \alpha !\beta !}\partial
_x^\alpha \partial _\xi ^\beta p(\rho_0)x^\alpha \xi ^\beta $. As
we shall see, $P_ 0$ has discrete spectrum and compact resolvent
in a weighted space and the \ev{}s can be computed explicitly.
(In fact, the spectrum is discrete even without weights and this fact
will be used in Section \ref{SectEv}.) 

\par Staying in the general case, we shall next give results about the
spectrum and the associated heat equation. We then define $P$ to
be the closure of $p^w$ with domain ${\cal S}$. In Section
\ref{Secweyl}, we shall see that the Fefferman-Phong inequality
implies that $\Re (Pu,u)\ge -Ch^2\Vert u\Vert ^2$, $u\in {\cal S}$
and hence also for $u\in {\cal D}(P)$ (this is immediate in the
KFP case). In other words, $P$ is accretive and we shall assume
\ekv{H5}{{\bf (H5)}\quad P\hbox{ is m-accretive,}} i.e. $P$ has no
accretive strict extension. In the KFP-case this has recently been
established in great generality by Helffer--Nier \cite{HelNi} and
their result implies (H5) under our assumptions on $V$. In the
general case, we shall see that the following assumption
\begin{equation*}
\textbf{(H6)} \ \ \  \textrm{ If $u\in L^2$ and $(p^w+1)u\in {\cal
S}$, then $u\in {\cal S}$,}
\end{equation*}
implies for $h$ \sufly{} small, that
${\cal D}(P)=\set{u\in L^2;\, Pu\in L^2 }$, and hence $P$
is m-accretive.

\begin{theo} \label{asympt}
Suppose $P$ satisfies  {\bf(H1--H5)}  and let $C>0$. Then there
exists $h_0>0$ such that for $0<h\leq h_0$, the spectrum of $P$ in
the disc $D(0,Ch)$ is discrete, and the \ev{}s are of the form,
\ekv{EigenValues} { \lambda _{j,k}(h)\sim h(\mu
_{j,k}+h^{1/N_{j,k}}\mu _{j,k,1}+ h^{2/N_{j,k}}\mu _{j,k,2}+...),
} where the ${\mu _{j,k}}$ are the \ev{}s in $D(0,C)$ (repeated
with their multplicity) of the quadratic approximation of ${P_|}_{h=1}$ at
the critical point $\rho_k$ and  $N_{j,k}$ is the dimension of the
corresponding generalized eigenspace.
\end{theo}
Here it is understood that $C$ has been chosen, so that no
quadratic approximation has any \ev{}s on the \bdy{} of the disc
$D(0,C)$. The explicit form of those \ev{}s will be given in
Proposition \ref{quad} and in Section \ref{secFP}. 
Note that they are distributed in an angle
 in $\R^++i\R$ avoiding the imaginary axis (except
 in 0).

As a consequence of the resolvent estimates and the description of
the eigenspaces we give the following theorem on the large time
behavior of the semi-group associated to $P$ :

\begin{theo} \label{evolution}
Suppose $P$ satisfies {\bf(H1)--(H5)}. Consider the set
$\set{\mu_{jk}}$ of eigenvalues of the quadratic approximation of
$P|_{h=1}$ at the critical points (repeated with their
multiplicities) defined in the preceding theorem. Let $b>0$ be
such that the line $\Re z=b$ avoids the set $\set{\mu_{jk}}$ and
define the finite set $J_b = \set{\mu_{j,k} ; \ \Re(\mu_{j,k}) <
b}$. Assume that the $\mu_{j,k}$ in $J_b$ are simple and
distinct.Then we have
\begin{equation}
e^{-tP/h} = \sum_{\mu_{j,k} \in J_b} e^{-t\lambda_{j,k}/h}
\Pi_{j,k} + \oo(1) e^{-tb} \ \ \ \ \text{ in } \ {\cal
L}(L^2,L^2),
\end{equation}
where $\lambda_{j,k}$ is the eigenvalue of $P$ associated to
$\mu_{j,k}$, and $\Pi_{j,k}$  the associated (rank one) spectral
projection. Here the term $\oo(1)$ is with respect to $t\geq 0$
and $h\rightarrow 0$.
\end{theo}

\bigskip
We
construct explicitly a global weight function $G$ with controlled
derivatives, satisfying in particular $G=\oo(h)$,
$G'=\oo(h^{1/2})$ and $G''=\oo(1)$.  The main idea (also used in many
earlier works on resonances and non-\sa{} \op{}s) 
is that we get the
new leading
symbol $p\approx p+{\eps \over i}\{ p,G\}$ with an increased real
part, where $\set{.,.}$ is the Poisson bracket, and $\eps$ is small and
fixed. We will use it both near the critical points of $p$ and at
infinity. Contrary to the earlier works mentioned above (but similarly
to \cite{DenckerSjostrandZworski2004}) we need resolvent and evolution
estimates in the original $L^2$ space and this requires $G/h$ to be
\bdd{} in order to have an equivalent norm on the weighted space. 
Consequently the estimates become more delicate. The
technical realization of this idea can be made either by using the
FBI-Bargmann transform and weighted spaces of \hol{} \fu{}s or using
pseudodifferential calculus (since $G/h$ is \bdd{}). We found it
convenient to use the first method near the critical points and the
second one elsewhere. We choose the semiclassical variant of the
Weyl-H{\"o}rmander calculus with a metric \sufly{} general to cover the
case of the KFP and related \op{}s.

The plan of the article is the following. The next section is devoted
to the construction of $G$.
In Sections 3 to 6 we work near the critical points by using the
Fourier-Bros-Iagolnitzer transform in a modified $L^2$ space
$L^2_{\Phi_\eps}$ associated to $G$. Here $G$ will play the role
of a local escape function. We recall in Section 3 some basic
facts about the FBI transform and construct the spaces
$L^2_{\Phi_\eps}$. In Section 4
 we get local resolvent estimates for a truncated operator
satisfying \bf(H1) \rm. In Section 5 we recall some facts on the
quadratic differential operators from \cite{Sjo74} and give a
localized version of them. Then in Section 6 we compare the
operator $P$ to its quadratic approximations at the critical points
to get precise local resolvent estimates near the critical points.

In Sections 7 to 9 we work away from the critical points of $p$ in
the real phase space using the semiclassical Weyl-H{\"o}rmander
calculus. Here $p$ satisfies hypothesis \bf (H2--H4)\rm. Section 7
is devoted to some basic facts about the semiclassical Weyl
calculus and the construction of a metric adapted to the symbol
$p$. In Sections 8 and 9 we get resolvent estimates  using a
multiplier method, where the symbol of the multiplier is
essentially $1+G/h$.

In Section 10 we combine all the resolvent estimates given in
Sections 3 to 9 and we prove Theorem \ref{main}. Section 11 is
devoted to the proof of Theorem \ref{asympt}, i.e.  the asymptotic
expansion of the eigenvalues of $P$: We solve a Grushin problem
thanks to a slight variation of the resolvent estimates given in
Section 10. In Section 12 we prove Theorem \ref{evolution} about
the large time behavior of the semigroup associated to $P$ under
hypothesis \bf (H5)\rm. Eventually in last section we check that
all the hypotheses \bf(H1--H4) \rm are satisfied for the symbol of
the KFP operator, which proves Theorem \ref{mainFP}.

\Section{Bounded weight function}\label{Section3}
\setcounter{equation}{0}

The aim of this section is to build a weight function $G$ defined
in the whole space, uniformly bounded by a multiple of $h$. Recall
that ${\cal B}$ is the fixed open ball appearing in \bf (H1)\rm.
The result is the following proposition:

\begin{prop} \label{propescape}
Suppose $p$ satisfies {\bf(H1--H4)}. Then  there exists a constant
$C>0$ and  a function $G \in \cc^\infty(\R^{2n}_\rho) $ such that
uniformly in $h$, $\eps>0 $ sufficiently small,  we have
\begin{equation} \label{derivG}
\begin{split}
   & \D^k G = \oo \sep{\delta^{(2-k)_+}} \ \ \
   \text{ for $\delta\lambda \leq
      h^{1/2}$,} \\
   & \D^k  G  = \oo \sep{ h (\delta \lambda h)^{-k/3}}
     \text{ in $\{ \rho \in {\cal B}; \delta \lambda
     \ge h^{1/2}\} $,} \\
   & \D_x^\alpha \D_\xi^\beta  G  = \oo \sep{ h^{1-k/3}
     \lambda^{-( \min(|\alpha|, 1) +|\beta| )/3}} \ \ \ \
      \text{ outside ${\cal B}$, $| \alpha |+|\beta |=k$.}
\end{split}
\end{equation}
Note that this implies $G  = \oo(h)$, $H_G  = \oo(h^{1/2})$
and $\D^2 G = \oo(1)$. Secondly $G$ is  such that
\begin{description}
\item{a)}
In ${\cal B}$, if we let $p$ denote an almost analytic extension
and if we put $\tilde{p}(\rho)
\defegal p(\rho +i\eps H_G(\rho)) =
   \tilde{p}_1(\rho) + i\tilde{p}_2(\rho)$ where
$\rho \in \R^{2n}$, we have
\begin{equation} \label{borneinfp}
\tilde{p}_1 \geq \frac{\eps}{C} \min \sep{ (\delta \lambda)^2,
(\delta  \lambda h)^{2/3}}, \ \ \ \  \tilde{p}_2 = \oo( (\delta
\lambda)^2).
\end{equation}
\item{b)}
Outside ${\cal B}$, we have
\begin{equation} \label{eqp}
p_1 + \eps H_{p_2} G \geq \frac{\eps}{C} (p_1 + (h\lambda)^{2/3}).
\end{equation}
\end{description}
\end{prop}

\subsection{The construction near the critical points}

Let $\rho _j\in {\cal C}$. Fix $T>0$. In a neighborhood of $\rho _j$,
we set
\begin{equation} \label{GT}
 G_T=\int k_T(t)\, p_1\circ \exp
(tH_{p_2})dt,
\end{equation}
 where $k_T(t)=k(t/T)$ and $k\in
\cc({\bf R}\setminus \{ 0\} )$ is the odd function given by:
$k(t)=0$ for $\vert t\vert \ge 1/2$ and $k'(t)=-1$ for $0<\vert
t\vert <1/2$. Notice that $k$ and $k_T$ have a jump of size $1$ at
the origin. $G_T$ is a smooth function satisfying
$$
H_{p_2}G_T=\langle p_1\rangle _T -p_1,\ G_T={\cal O}(\delta ^2),\
\nabla G_T={\cal O}(\delta ),
$$
 where
 $$ \langle
p_1\rangle _T={1\over T}\int_{-T/2}^{T/2}p_1\circ \exp (tH_{p_2})
dt.
$$

 Consider the dilated symbol
$$
\widetilde{p}=\widetilde{p}_\epsilon(\rho )=
p(\rho +i\epsilon H_G(\rho ))=p(\rho
)-i\epsilon H_pG(\rho )+{\cal O}(\epsilon ^2\vert \nabla G\vert
^2),
$$
with real and imaginary parts, given by
\begin{equation}
\begin{split}
\widetilde{p}_1 & =p_1(\rho )+\epsilon H_{p_2}G(\rho )+{\cal
O}(\epsilon ^2\vert \nabla G\vert ^2)  \\
         &=(1-\epsilon
       )p_1(\rho )+\epsilon \langle p_1\rangle _T+{\cal O}_T(\epsilon
               ^2\delta ^2), \\
        \widetilde{p}_2&=p_2(\rho )-\epsilon
H_{p_1}G(\rho )+{\cal O}(\epsilon ^2\vert \nabla G\vert ^2).
\end{split}
\end{equation}
 Using (\ref{mainhyp}) near $\cc$, we see that if we fix $\epsilon
>0$ small enough, depending on $T$, then in an
$(\epsilon ,T)$ - dependent
neighborhood of $\rho _j$, we have
\begin{equation} \label{p1leqh12}
\widetilde{p}_1 \geq \frac{\eps}{C} \delta ^2,\
\widetilde{p}_2={\cal O}(\delta ^2).
\end{equation}
Note in particular that $\widetilde{p}$ takes its values in an
angle around the positive real axis,
 $ \widetilde{p}_1\succ \widetilde{p}_2.$ Note also that another
 choice of weight function near the critical point could have
 been $\tau H_{p_2}p_1$ for $\tau$
 sufficiently small.

\subsection{The construction away from the critical points}

We work in a region
\begin{equation} \label{Ah}
  \{ \rho ;\, \delta\lambda(\rho) \ge h^{1/2}\}.
\end{equation}
Let $\psi\in \cc_0^\infty (]-2,2[$) be a cutoff function
equal to $1$ in $[-1,
1]$. Let $M$ be a large constant to be fixed later.
We choose the following function
\begin{equation} \label{defG}
G = h \frac{H_{p_2}p_1}{(\delta\lambda)^{4/3} h^{1/3}} \psi\sep{
\frac{Mp_1}{(h\delta \lambda)^{2/3}}},
\end{equation}
where we recall that $\delta = \delta(\rho)$ and $\lambda =
\lambda(\rho)$.

We first check  the bounds for the derivatives of $G$. Of course when
$Mp_1 \geq 2
(h \delta \lambda)^{2/3}$, $G=0$ and we have only to
study the derivatives in
the region where $Mp_1 < 2 (h \delta \lambda)^{2/3}$.

Observe that the estimates (\ref{derivG})   for $G$ in $\set{
\delta\lambda \geq h^{1/2 }}$ can be equivalently written using
the following Riemannian metric
$$
\Gamma_h = { dx^2 \over (\delta h )^{2/3}} + {d\xi^2 \over (\delta
\lambda h )^{2/3}} ,
$$
by saying (in the H{\"o}rmander terminology of spaces of symbols, see
Section \ref{Secweyl}) that
\begin{lemma} \label{lemestimG}
 $ G \in S ( h, \Gamma_h)$ and $\nabla G \in S
( h( \delta \lambda h)^{-1/3}, \Gamma_h).$
\end{lemma}

\preuve For the following estimates of the derivatives  we shall
use  this
terminology and stay in the region $\set{ \delta
\lambda \geq h^{1/2 }}
\cap \set{ Mp_1 < 2 (h \delta \lambda)^{2/3}}$. We work step by
step by studying the derivatives of each function entering in the
composition of $G$.

\ititem{Estimates of $p$. } We know that $ p \in S(\lambda^2, dx^2
+ d\xi^2 / \lambda^2)$ from the hypothesis. >From the fact that $p$
is a Morse function we get that $ p \in
S((\delta\lambda)^2, dx^2/\delta^2 + d\xi^2 / (\delta\lambda)^2)$.
For the same reason we have $\nabla p \in S(\delta\lambda,
dx^2/\delta^2 + d\xi^2 / (\delta\lambda)^2)$. Besides we have on
$\set{ \delta\lambda  \geq h^{1/2 }}$
\begin{equation} \label{ordermetric}
\Gamma_h \ge dx^2/\delta^2 + d\xi^2 / (\delta\lambda)^2
\ge C^{-1}\Gamma _0,
\end{equation}
since $\delta^2 \geq (\delta h )^{2/3}$ and
$(\delta\lambda)^2 \geq (\delta
\lambda h )^{2/3}$ in this region. As a consequence we get that
\begin{equation}\label{estimp}
    \nabla p \in S( \delta \lambda, \Gamma_h).
\end{equation}

\ititem{Estimates for $p_1$.} Since $p_1$ is nonnegative with
bounded
second derivatives,  we can apply the well known inequality
for $W^{2, \infty}$ functions
\begin{equation} \label{deriv}
|\nabla f|^2 \leq 2 f \norm{f''}_\infty,
\end{equation}
which yields $\abs{\nabla p_1} \leq C \sqrt{p_1}$.
 Since $p_1 \leq 2 (h \delta
\lambda)^{2/3}$ we get that
$ \nabla p_1 = \oo((h \delta \lambda)^{1/3})$.
Together with the fact that $p_1$ has its second
derivative bounded and
(\ref{ordermetric}) we get that
\begin{equation}\label{estimp1}
     p_1 \in S( (\delta \lambda h)^{2/3}, \Gamma_h),
     \ \ \ \ \text{and} \ \ \ \
     \nabla p_1 \in S( (\delta \lambda h)^{1/3}, \Gamma_h).
\end{equation}
Here we used that $\nabla ^2p_1\in S(1,\Gamma _0)\subset
S(1,\Gamma _h)$.

\ititem{Estimates for powers of $\delta \lambda$.}  Using
(\ref{symblambda}), we first note that
$$
\delta \lambda \in S((\delta\lambda), dx^2/\delta^2 + d\xi^2 /
(\delta\lambda)^2)
$$
Together with the fact that $ \nabla (\delta \lambda)
\in S(1,dx^2/\delta^2 +
d\xi^2 / (\delta\lambda)^2)$, this gives for $\alpha \in \R$,
\begin{equation}\label{estimdeltlamb}
    (\delta \lambda)^\alpha
    \in S( (\delta \lambda)^{\alpha}, \Gamma_h),
    \ \ \ \ \text{ and } \ \ \ \
    \nabla (\delta \lambda)^\alpha
    \in S( (\delta \lambda)^{\alpha-1}, \Gamma_h).
\end{equation}

\ititem{Estimates of $p_1 / (h \delta \lambda)^{2/3}$.} From
(\ref{estimp1}) and (\ref{estimdeltlamb}) with $\alpha = -2/3$ we
get immediately that
$$
p_1 / (h \delta \lambda)^{2/3} \in S( 1 , \Gamma_h).
$$
Besides let us write
$$
\nabla \sep{ p_1 / (h \delta \lambda)^{-2/3} } = (\nabla p_1 )(h
\delta \lambda)^{-2/3} + p_1  \nabla ( h \delta \lambda)^{-2/3}.
$$ From the same estimates for the derivatives  we get
$$
(\nabla p_1 )(h \delta \lambda)^{-2/3}
\in S( (h \delta \lambda)^{-1/3},
\Gamma_h),
$$
and
$$
p_1  \nabla (h \delta \lambda)^{-2/3}
\in S( (h \delta \lambda)^{2/3} \times
h^{-2/3} ( \delta \lambda)^{-5/3}  , \Gamma_h) \subset S( (h \delta
\lambda)^{-1/3}, \Gamma_h),
$$
where in the last inclusion we used the fact that $\delta \lambda
\geq h^{1/2}$. Summing up we have proven that
\begin{equation}\label{estimp1delt}
p_1 / (h \delta \lambda)^{2/3} \in S( 1 , \Gamma_h), \ \ \ \text{
  and } \ \ \ \nabla
\sep{ p_1 / (h \delta \lambda)^{2/3}} \in
 S( (h \delta \lambda)^{-1/3},
\Gamma_h).
\end{equation}

\ititem{Estimates of $\psi(M p_1 / (h \delta \lambda)^{2/3})$.}
An immediate consequence of the first part of  (\ref{estimp1delt})
 is that
$$
\psi(M p_1 / (h \delta \lambda)^{2/3}) \in S( 1 , \Gamma_h),
$$
since $\psi$ is $\cc^\infty$ with compact support. We need to estimate
the derivatives of this expression,
$$
\nabla \psi(M p_1 / (h \delta \lambda)^{2/3}) = M\nabla \sep{ p_1
/ (h \delta \lambda)^{2/3}} \psi'(M p_1 / (h \delta
\lambda)^{2/3}).
$$
For the same reason as before we have
$$
\psi'(M p_1 / (h \delta \lambda)^{2/3}) \in S( 1 , \Gamma_h).
$$
Using the second part of (\ref{estimp1delt}), and summing
 up we
have proven that
\begin{equation}\label{estimpsi}
  \psi(M p_1 / (h \delta \lambda)^{2/3}) \in S( 1 , \Gamma_h), \ \ \
  \text{ and } \ \ \
  \nabla \psi(M p_1 / (h \delta \lambda)^{2/3}) \in  S( (h \delta
    \lambda)^{-1/3}, \Gamma_h).
\end{equation}

\ititem{Estimates for $H_{p_2} p_1$.} We observe that $H_{p_2}
p_1= \sigma (\nabla p_2, \nabla p_1)$ where $\sigma$ is the
canonical symplectic form. Using (\ref{estimp}) for $p_2$ and
(\ref{estimp1}) for $p_1$ we get
$$
 H_{p_2} p_1 \in S(
(\delta \lambda) (h \delta \lambda)^{1/3}, \Gamma_h).
$$
From the hypothesis (\ref{symbcal}) and the fact that $p$ is a
Morse function we can  write
$$
\nabla H_{p_2} p_1 \in S( \delta \lambda, dx^2/\delta^2 + d\xi^2 /
(\delta\lambda)^2 ) \subset S( \delta \lambda, \Gamma_h).
$$
Summing up we have proven that
\begin{equation}\label{estimp2p1}
      H_{p_2} p_1 \in S( h^{1/3}  ( \delta \lambda)^{4/3}, \Gamma_h),
      \ \ \ \text{ and } \ \ \  \nabla H_{p_2} p_1
      \in  S( \delta \lambda,
\Gamma_h).
\end{equation}

\ititem{Estimates for $H_{p_2} p_1 / (h^{1/3} (\delta
\lambda)^{4/3})$.} From the first parts of (\ref{estimp2p1}) and
(\ref{estimdeltlamb}) with $\alpha = -4/3$ we immediately get that
$H_{p_2} p_1 / (h^{1/3} (\delta \lambda)^{4/3}) \in S(1,
\Gamma_h)$. Its derivative is given by
$$
\nabla { H_{p_2} p_1 \over  h^{1/3} (\delta \lambda)^{4/3}} =
{\nabla  H_{p_2} p_1 \over  h^{1/3} (\delta \lambda)^{4/3}} +  {
H_{p_2} p_1 \nabla (h^{-1/3} (\delta \lambda)^{-4/3})}.
$$
Using (\ref{estimp2p1}) and (\ref{estimdeltlamb}) we respectively
 get that
$$
{\nabla  H_{p_2} p_1 \over  h^{1/3} (\delta \lambda)^{4/3}}
\in S((h\delta
\lambda)^{-1/3}, \Gamma_h),
$$
and
$$
{ H_{p_2} p_1 \nabla (h^{-1/3} (\delta \lambda)^{-4/3})}
\in S(h^{1/3} (\delta
\lambda)^{4/3} \times h^{-1/3} \delta \lambda)^{-7/3}, \Gamma_h )
 \subset S(
(\delta \lambda)^{-1}, \Gamma_h).
$$
Using the fact that $\delta \lambda \geq (\delta \lambda h)^{1/3}$
in this formula gives
\begin{equation}\label{estimp2p1delt}
  {H_{p_2} p_1 \over h^{1/3} (\delta \lambda)^{4/3}}
  \in S(1, \Gamma_h),
  \ \ \ \text{ and } \ \ \
  \nabla {H_{p_2} p_1 \over h^{1/3} (\delta \lambda)^{4/3}}
    \in S((h\delta \lambda)^{-1/3}, \Gamma_h).
\end{equation}

\ititem{Estimates for   $G$ and end of the proof of lemma
\ref{lemestimG}.} We can now prove the estimates for $G$. From the
first parts of (\ref{estimpsi}) and (\ref{estimp2p1delt}) and
multiplying by $h$ we get that
$$
G \in S(h, \Gamma_h).
$$
From the second part of the same expressions we also get
immediately that
$$
\nabla G  \in S(h(h\delta \lambda)^{-1/3}, \Gamma_h).
$$
This completes the proof of lemma \ref{lemestimG} and therefore of
the estimates (\ref{derivG}) when $\delta\lambda\geq h^{1/2}$. \fin

\subsection{Proof of (\ref{borneinfp}) in the intermediate
region}

We work here in the region $\set{\rho\in {\cal B}; \ h^{1/2}
 \leq   \delta
\lambda }$, but many of the estimates will be valid also near infinity
and used later, so we indicate when the validity is restricted to a
bounded region.
Consider the function $G$ defined in (\ref{defG}) :
$$
G = h \frac{H_{p_2}p_1}{(\delta\lambda)^{4/3} h^{1/3}} \psi\sep{
\frac{Mp_1}{(h\delta \lambda)^{2/3}}}.
$$
For $\tilde{p}(\rho) \defegal p(\rho+i\eps H_G(\rho)) =
\tilde{p}_1(\rho) + i \tilde{p}_2(\rho)$ in ${\cal B}$, we have
\begin{equation}
\begin{split}
\widetilde{p}_1 & =p_1+\epsilon H_{p_2}G+{\cal
O}(\epsilon ^2\vert \nabla G\vert ^2), \\
        \widetilde{p}_2&=p_2-\epsilon
H_{p_1}G+{\cal O}(\epsilon ^2\vert \nabla G\vert ^2).
\end{split}
\end{equation}
Let us estimate the remainders. >From (\ref{derivG}) we know that
 $ \nabla G =
\oo(h^{2/3} (\delta\lambda)^{-1/3})$. As a consequence
\begin{equation} \label{remainders}
{\cal O}(\epsilon ^2\vert \nabla G\vert ^2)
    =  \eps^2 \oo(h^{4/3} (\delta\lambda)^{-2/3}  )
   \le  \eps^2 \oo( (h\delta\lambda)^{2/3}  ),
\end{equation}
since $h^{4/3} (\delta\lambda)^{-2/3} =
\oo((h\delta\lambda)^{2/3}) $ when $\delta\lambda \geq h^{1/2}$.
Let us now study the first two terms of the expression of
$\tilde{p}_1$ depending on the size of $p_1$.

\ititem{Estimates when $p_1$ is large.}
We work first in the \it elliptic \rm region
$$
\set{ \rho\in \R^{2n}; \ \ Mp_1 \geq  (h\delta\lambda )^{2/3}}.
$$
From (\ref{derivG})  and the fact that $H_{p_2} =
\oo(\delta\lambda)$, we get
\begin{equation} \label{hp2g}
H_{p_2}G = \oo((\delta\lambda h)^{2/3}).
\end{equation}
Restricting the attention to ${\cal B}$ we recall that the remainder
in (\ref{remainders}) is $ \eps^2 \oo((\delta\lambda h)^{2/3} ) $
 and that
$$
\tilde{p}_1 = p_1 + \eps H_{p_2}G
+ \eps^2 \oo((\delta\lambda h)^{2/3} ).
$$
Choosing $\eps$ small enough yields
$$
\tilde{p}_1 \geq \frac{(\delta\lambda  h)^{2/3}}{CM}.
$$
On the other hand we have using the bound on the remainder
and of $H_G$ that
$$
\tilde{p}_2 = \oo((\delta\lambda)^2).
$$

\ititem{Estimates when $p_1$ is small.} In the region
\begin{equation} \label{Hppetit}
\set{\rho\in \R^{2n}; \ \ Mp_1 \leq (h\delta \lambda)^{2/3}},
\end{equation}
we can write
  $G = h \frac{H_{p_2}p_1}{(\delta\lambda)^{4/3}
h^{1/3}}$. We have therefore
\begin{equation} \label{HGp}
 p_1 + \eps H_{p_2}G  = p_1 +  \eps h
\frac{H_{p_2}^2p_1}{(\delta\lambda)^{4/3} h^{1/3}} + \eps h
(H_{p_2}p_1) H_{p_2}\sep{(\delta\lambda)^{-4/3} h^{-1/3}}.
\end{equation}
For the third term of (\ref{HGp}) we use that $|\nabla p_1|\le
C\sqrt{p_1}\le C(h\delta \lambda )^{1/3}/\sqrt{M}$, $|\nabla
p_2|\le C\delta \lambda $ and get using also (\ref{estimdeltlamb})
\begin{equation}
\eps h  (H_{p_2}p_1) H_{p_2}\sep{(\delta\lambda)^{-4/3} h^{-1/3}}
= \frac{\eps}{\sqrt{M}} \oo(h)
=  \frac{\eps}{\sqrt{M}} \oo((\delta \lambda
h)^{2/3}),
\end{equation}
since $\delta\lambda \geq h^{1/2}$. We study next the sum of the
first and the second term. We first observe that
$$
\frac{h}{(\delta\lambda)^{4/3} h^{1/3}}  \leq 1 ,
$$
 and from (\ref{mainhyp}) provided
 $\eps<\eps_0$, we get
\begin{equation} \label{aussiloin}
\begin{split}
p_1 +  \eps h \frac{H_{p_2}^2p_1}{(\delta\lambda)^{4/3} h^{1/3}}
& \geq
\frac{\eps}{\eps_0}
 \frac{h}{(\delta\lambda)^{4/3} h^{1/3}}(p_1 + \eps_0
H_{p_2}^2p_1) \geq \frac{\eps\eps_1}{\eps_0}
(\delta \lambda h)^{2/3}.
\end{split}
\end{equation}
Therefore choosing $M$ sufficiently large (and fixed from now on)
 gives
\begin{equation} \label{sommephg}
p_1 +  \eps  H_{p_2} G \geq \eps  (\delta \lambda h)^{2/3}/C.
\end{equation}
Since the remainder term in (\ref{remainders}) is
$\eps^2\oo((\delta\lambda h)^{2/3})$,
and choosing $\eps$ sufficiently small again, we get on ${\cal B}$:
$$
\tilde{p}_1 \geq \eps (\delta \lambda h)^{2/3}/C \ \ \
 \hbox{ for }\delta \ge
h^{1\over 2}.
$$

\subsection{The global construction}

We shall glue together the two weights constructed in the previous
two subsections. Let us denote by $G_{\rm int}$ the interior
weight $G_T$ defined in (\ref{GT}) and $G_{\rm out}$ the one
defined in (\ref{defG}) where we recall that the constants $T$,
and $M$ appearing in the definitions are fixed. Recall also the
main properties of these weights:

\begin{equation}
\left\{
\begin{array}{ll}
  \D^kG_{\rm int} = \oo(\delta^{(2-k)_+}), &  \\
  p_1 + \eps H_{p_2} G_{\rm int} \geq \frac{\eps}{C} (\delta
  \lambda)^2, &
  \end{array}
\right. \text{in} \ \  {\cal B},
\end{equation}
and
\begin{equation}
\left\{
\begin{array}{ll}
  \D^kG_{\rm out} = \oo(h(\delta\lambda h)^{-k/3}),
   &  \\
  p_1 + \eps H_{p_2} G_{\rm out} \geq \frac{\eps}{C} (\delta
  \lambda h)^{2/3}, &
  \end{array}
\right.
 \text{ for } h^{1/2} \leq  \delta
  \lambda.
\end{equation}
We now build a function $G$ defined everywhere and satisfying
Proposition \ref{propescape}. In the following, we introduce an
additional large real constant $N$ to be fixed later. We first
build modified functions $G_{\rm int}$ and $G_{\rm out}$.

\ititem{Construction of a modified $G_{\rm int}$.} Let us
introduce the following function
 $$
 \widehat{p}_1 \defegal \chi \sep{{\delta \lambda
 \over 4Nh^{1/2}}}p_1,
 $$
 where $\chi \in C_0^\infty (\R ;[0,1])$ is a standard
cut-off near 0, equal to $1$ on $[0,1]$ and to $0$ on $[2, + \infty[$.
Notice that $\widehat{p}_1 = 0$ when $\delta\lambda \geq
8Nh^{1/2}$ and that $\widehat{p}_1 = p_1$ when $\delta\lambda \leq
4Nh^{1/2}$.   Define now  $\widetilde{G}_{\rm int}$ as in
(\ref{GT})  but with $p_1$ there replaced by $\widehat{p}_1$. Then
$\widetilde{G}_{{\rm int}}$ has its support in $\{ \rho ;\, \delta
\lambda\leq 16Nh^{1/2}\}$ and coincides with $G_{\rm int}$ when
$\delta \lambda\leq 2Nh^{1/2}$ (assuming that $T$ has been fixed
 \sufly{} small). As a consequence we get that
\begin{equation} \label{in}
p_1+\epsilon H_{p_2} \widetilde{G}_{\rm int}\geq \left\{
  \begin{array}{ll}
  \frac{\epsilon}{C}(\delta\lambda  )^2,
      & \text{when } \delta \lambda \le 2Nh^{1/2}, \\
   0  &  \text{everywhere}.
\end{array}
\right.
\end{equation}
Note that this implies the following bounds :
\begin{equation} \label{AAA}
p_1+\epsilon H_{p_2} \widetilde{G}_{\rm int}\geq \left\{
  \begin{array}{ll}
  \frac{\epsilon}{C}(\delta\lambda  )^2,
      & \text{when } \delta \lambda \le h^{1/2}, \\
  \frac{\epsilon}{C}(\delta\lambda h )^{2/3},
      & \text{when } h^{1/2} \leq \delta \lambda \le N h^{1/2}, \\
  \frac{\epsilon}{C}N^{4/3}(\delta\lambda h )^{2/3},
      & \text{when } Nh^{1/2}\leq \delta \lambda \le 2N h^{1/2}, \\
   0,  & \text{when }  \delta \lambda \geq 2N h^{1/2},
\end{array}
\right.
\end{equation}
where for the second bound we used the fact that $(\delta
\lambda)^2 \geq (\delta \lambda h)^{2/3}$ when $\delta \lambda
\geq h^{1/2}$, and for the third bound the fact that $(\delta
\lambda)^2 \geq N^{4/3} (\delta \lambda h)^{2/3}$ when $\delta
\lambda \geq N h^{1/2}$.

Let us now study the derivatives of $\widetilde{G}_{\rm int}$.
Since $\widetilde{G}_{\rm int} = {G}_{\rm int}$ when $\delta
\lambda \leq 2Nh^{1/2}$ we get that
\begin{equation} \label{BBB}
\D^k\widetilde{G}_{\rm int}={\cal O}(\delta^{(2-k)_+}) \ \ \
\text{ when } \delta \lambda \leq 2Nh^{1/2}.
\end{equation}
When $2Nh^{1/2} \leq \delta \lambda \leq 16N h^{1/2}$,
$\widetilde{G}_{\rm int}$ inherits the properties of
$\widehat{p}_1$, i.e. $\D^k\widetilde{G}_{\rm int}={\cal
O}((Nh^{1/2})^{(2-k)})$ which yields
\begin{equation} \label{CCC}
 \D^k\widetilde{G}_{\rm int}
 =C_k(N) {\cal O}(h (h\delta\lambda)^{-k/3})
 \ \ \ \text{ when }   h^{1/2} \leq \delta\lambda \leq \oo(1),
\end{equation}
(of course this estimate is true when $ \delta\lambda \geq
16Nh^{1/2}$ since $\widetilde{G}_{\rm int}$ is zero there).

\ititem{Construction of a modified $G_{\rm out}$.} For the same
$\chi $ as before define
$$
\tilde{G}_{\rm out}(\rho) \defegal G_{\rm out} (\rho)
\sep{1-\chi\sep{ {\delta(\rho) \lambda(\rho) \over N h^{1/2}}}}.
$$
We notice that $\tilde{G}_{\rm out} = {G}_{\rm out}$ when
$\lambda\delta \ge 2Nh^{1\over 2}$ and  $\tilde{G}_{\rm out} = 0$
when $\lambda\delta \le Nh^{1\over 2}$. Therefore we have directly
\begin{equation}
p_1 +\epsilon H_{p_2} \tilde{G}_{\rm out}\geq \left\{
  \begin{array}{ll}
     {\epsilon \over C}(\delta \lambda  h)^{2/ 3},\ \ \
        & \hbox{ for }\lambda\delta \ge 2N h^{1\over 2}, \\
  0 ,   & \text{ when } \delta \lambda \leq N h^{1/2}.
  \end{array}
\right.
\end{equation}
In the area
$
Nh^{1\over 2} \leq \lambda\delta \leq 2Nh^{1\over 2},
$
we can write uniformly in $N\geq 1$ that
\begin{equation}\label{2.33}
\begin{split}
p_1 + \eps H_{p_2}\tilde{G}_{\rm out}
& \ge (p_1+\epsilon H_{p_2} G_{\rm out})
(1-\chi) -\eps G_{{\rm out}} H_{p_2}\chi \\
& \geq -\eps  \abs{G_{\rm out}}.
\abs{H_{p_2} \chi} .
\end{split}
\end{equation}
We know that $ G_{\rm out} = \oo(h)$ and that $
\abs{H_{p_2} \chi} \leq \abs{\D p_2} \abs{ \D \chi} =
\oo(\delta\lambda \frac{1}{N h^{1/2}}). $ This gives
$$
\abs{ G_{\rm out}}\abs{H_{p_2} \chi} = \oo\sep{ \frac{\delta
\lambda h^{1/2}}{N}}.
$$
Now using the fact that $\delta \lambda = (\delta
\lambda)^{2/3}(\delta \lambda)^{1/3} \leq (\delta \lambda)^{2/3}
(2N)^{1/3} h^{1/6}$, we deduce that uniformly in $N \geq 1$,
$$
\abs{ G_{\rm out}}\abs{H_{p_2} \chi} = \oo\sep{ \frac{(\delta
\lambda h)^{2/3}}{N^{2/3}}} =  \oo\sep{ (\delta \lambda h)^{2/3}}.
$$
Using this and (\ref{2.33}) we find that
$$
  p_1 + \eps H_{p_2}\tilde{G}_{\rm out}
   \geq - \eps\oo\sep{ (\delta
  \lambda h)^{2/3}}, \ \ \ \text{ when }
  Nh^{1\over 2} \leq \lambda\delta \leq 2Nh^{1\over 2}.
$$
Eventually we have the following bounds in the whole region
$\delta \lambda \leq \oo(1)$ :
\begin{equation} \label{DDD}
 p_1+\epsilon H_{p_2} \widetilde{G}_{\rm
out}\ge \left\{
\begin{array}{ll}
  0,& \text{ when } \delta \lambda \le Nh^{1/2}, \\
  -C\eps (\delta \lambda h)^{2/3},
       & \text{ when } Nh^{1/ 2}
       \leq \delta \lambda \leq 2N h^{1/ 2}, \\
  \frac{\epsilon}{C}(\delta\lambda h)^{2/3}, \ \ \
      & \text{ when } \delta \lambda \ge 2Nh^{1/2}.
\end{array}
\right.
\end{equation}
For the derivatives of $\widetilde{G}_{\rm out}$ we
can write immediately
\begin{equation} \label{EEE}
  \D^k \widetilde{G}_{\rm out} = 0 = \oo(\delta^{(2-k)_+}), \ \ \
  \text{ when } \delta \lambda \leq h^{1/2},
\end{equation}
since $\widetilde{G}_{\rm out} = 0$ there. In the intermediate
region we check that
$$
\D^k \sep{\chi\sep{ \delta \lambda \over N h^{1/2}}} = C_k(N)
\oo(h^{-k/2}) = C'_k(N) \oo((\delta \lambda h)^{-k/3})
$$
since $\delta \lambda \leq 2N h^{1/2}$. Of course the same
estimate is true in the larger region  $\{h^{1/2} \leq  \delta\lambda
\}\cap{\cal B}$, since $\chi$ is compactly supported. Now using the
fact that $\D^kG_{\rm out} = \oo(h(\delta\lambda h)^{-k/3})$, we
get the same estimate for $\tilde{G}_{\rm out}$
\begin{equation} \label{FFF}
  \D^k\tilde{G}_{\rm out} = C_k(N) \oo(h(\delta\lambda h)^{-k/3}), \
  \ \ \text{ when } h^{1/2} \delta \lambda \leq \oo(1).
\end{equation}
The construction of $\tilde{G}_{\rm out}$ is complete.

\ititem{Construction of the weight function $G$.}  We finally pose
\begin{equation} \label{definitiveG}
 G= (\widetilde{G}_{\rm in}+ \widetilde{G}_{\rm out})/2.
 \end{equation}
Using the bounds (\ref{BBB}, \ref{CCC}, \ref{EEE}, \ref{FFF}) for
the derivatives of $\widetilde{G}_{\rm in}$ and
$\widetilde{G}_{\rm out}$ we immediately get that
\begin{equation}
\D^k G = \left\{
  \begin{array}{ll}
  \oo(\delta^{(2-k)_+}), \ \ \
  &  \text{ when } \delta \lambda \leq h^{1/2}, \\
  C'_k(N) \oo(h(\delta \lambda h)^{-k/3}), \ \ \
  & \text{ in ${\cal B}$ when } h^{1/2}
  \leq  \delta\lambda \leq \oo(1),
  \end{array}
  \right.
\end{equation}
i.e. the bounds given in the first two estimates of
(\ref{derivG}). On the other hand, combining (\ref{AAA}) and
(\ref{DDD}) gives
\begin{equation*}
2 p_1+ 2\epsilon H_{p_2} G \ge \left\{
\begin{array}{ll}
  \frac{\epsilon}{C}(\delta\lambda )^{2}, \ \ \
      & \text{ when } \delta \lambda \leq h^{1/2},  \\
  \frac{\epsilon}{C}(\delta\lambda h)^{2/3}, \ \ \
      & \text{ when } h^{1/2} \leq \delta \lambda \leq Nh^{1/2}, \\
  \sep{ \frac{\eps}{C} N^{4/3} -C\eps} (\delta \lambda h)^{2/3},
       & \text{ when } Nh^{1/ 2} \leq
       \delta \lambda \leq 2N h^{1/ 2}, \\
  \frac{\epsilon}{C}(\delta\lambda h)^{2/3}, \ \ \
      & \text{ when } 2Nh^{1/2} \leq \delta \lambda \oo(1) .
\end{array}
\right.
\end{equation*}
Taking $N$ sufficiently large and fixed from now on, and dividing
by $2$  gives with a new constant $C$
\begin{equation} \label{GGG}
 p_1+ \epsilon H_{p_2} G \ge \left\{
\begin{array}{ll}
  \frac{\epsilon}{C}(\delta\lambda )^{2}, \ \ \
      & \text{ when } \delta \lambda \leq h^{1/2},  \\
  \frac{\epsilon}{C}(\delta\lambda h)^{2/3}, \ \ \
      & \text{ when } h^{1/2} \leq \delta \lambda \leq \oo(1).
\end{array}
\right.
\end{equation}
Let us now prove (\ref{borneinfp}). This was already proven in
(\ref{p1leqh12}) in the region $\delta\lambda \leq h^{1/2}$ since
$G= G_T$ there. In the region $h^{1/2} \leq \delta \lambda \leq
\oo(1)$ we follow the same procedure. We write
$$
\widetilde{p}(\rho )=p(\rho +i\epsilon H_G(\rho ))=p(\rho
)-i\epsilon H_pG(\rho )+{\cal O}(\epsilon ^2\vert \nabla G\vert
^2),
$$
with real part given by
\begin{equation*}
\begin{split}
\widetilde{p}_1
        & =p_1(\rho )+\epsilon H_{p_2}G(\rho )
          +\oo(\epsilon ^2\vert \nabla G\vert ^2)
          =p_1(\rho )+\epsilon H_{p_2}G(\rho ) +
         \eps^2\oo( (\delta\lambda h)^{2/3}),
\end{split}
\end{equation*}
since $\nabla G = \oo((\delta\lambda h)^{1/3})$ by Lemma
\ref{lemestimG} and the fact that $\delta \lambda \ge h^{1/2}$.
Using (\ref{GGG}) and taking $\eps $ small enough
yields
$$
\widetilde{p}_1 \geq \frac{\eps}{C}(\delta\lambda h)^{2/3} .
$$
For the imaginary part $\widetilde{p}_2$ we directly write
$$
\widetilde{p}_2 =p_2(\rho )-\epsilon H_{p_1}G(\rho )
        +\oo (\epsilon ^2\vert \nabla G\vert ^2)  = \oo(\delta^2).
$$
This completes the proof of Proposition \ref{propescape} in the
region $\delta \lambda \leq \oo(1)$.

\ititem{End of the proof of  Proposition \ref{propescape}.} We now
work outside ${\cal B}$. We first observe that the estimate
(\ref{sommephg}) remains valid, therefore in the region $
\set{\rho; \ \ Mp_1(\rho) \leq (h\delta
(\rho)\lambda(\rho))^{2/3}}$ we get (\ref{eqp}) from
(\ref{aussiloin}) and (\ref{Hppetit}). In the region $\set{\rho; \
\ Mp_1(\rho) \geq (h\delta (\rho)\lambda(\rho))^{2/3}}$ we use
(\ref{hp2g}) and for $\eps$ small enough  we get
$$
p_1 + \eps H_{p_2} G \geq { \eps \over C} (p_1 + (h\delta
\lambda)^{2/3}).
$$
 The proof of Proposition \ref{propescape} is
complete. \fin

\Section{Review of FBI tools} \setcounter{equation}{0}
\label{secFBI}

The aim of this section is to review the definitions about the FBI
transform and the spaces associated to a function $G$ satisfying
the estimates of Proposition \ref{propescape} in a bounded region
and equal to $0$ elsewhere. Note in particular that it has its
second derivative bounded. The material here is essentially taken
from \cite{Sjo90}. In this section, and in Sections
\ref{locallarge} and \ref{localsmall}, we suppose that the symbol
$p$ satisfies hypothesis \bf (H1) \rm and is bounded with all its
derivatives everywhere.

\subsection{Definitions and main properties}

\par Let $T$ be a FBI-Bargmann \tf{}:
\ekv{7.7} { Tu(x)=Ch^{-{3n\over 4}}\int e^{{i\over h}
\phi (x,y)}u(y)dy, }
where we may choose
$\phi (x,y)={i\over 2}(x-y)^2$ as in the standard Bargmann transform.
Other quadratic $\phi $ with the general properties reviewed in
\cite{Sjo95} are also possible. The associated \ctf{} is given by
\ekv{5.9} { \kappa _T:(y,-\partial _y\phi (x,y))
  \mapsto (x,\partial _x\phi (x,y)). }
We have the associated IR-space (see \cite{Sjo95} for the
 terminology),
\ekv{5.10} { \Lambda _{\Phi _0}=\kappa _{T}(\R^{2n}),\
  \Phi _0(x)=-\Im \phi (x,y_0(x)), }
where $y_0$ is the point where $\R^n\ni y\mapsto -\Im \phi (x,y)$
takes its non-degenerate maximum.

If $P=p^w$,
then by the metaplectic invariance,
\ekv{7.8} { TP=\widehat{P}T, \ \widehat{P}=\widehat{p}^w,}
we have the exact symbol relation: \ekv{7.9} {
\widehat{p}\circ \kappa _T=p. } Shortly, we will recall the definition
of the Weyl quantization on the FBI-transform side.

\par From now on, we work entirely on the FBI-side, and we
 shall write $P$
instead of $\widehat{P}$ and similarly for the symbols. We introduce
the spaces $L^2_{\Phi _0}=L^2(\C ^n;e^{-2\Phi _0/h}L(dx))$,
where $L(dx)$ is the
Lebesgue measure, and $H_{\Phi_0 }$ the subspace of entire functions.
The Weyl-quantization on $H_{\Phi _0}$
takes the form of a contour integral 
\ekv{7.10} { Pu(x)={1\over (2\pi h)^n}\iint_{\theta
={2\over i}{\partial \Phi _0\over
\partial x}({x+y\over 2})}e^{i(x-y)\cdot \theta /h}
p({x+y\over 2},\theta ;h)u(y)dyd\theta .
}  By $p$, we also denote an almost \hol{} extension of $p$ to a
tubular \neigh{} of $\Lambda _{\Phi _0}$. If we introduce a
$\cc^\infty$ function $\psi_0$ equal to $1$ near $0$, we get
 for $u\in
H_{\Phi _0}$:
$$
Pu(x)={1\over (2\pi h)^n}\iint_{\theta ={2\over i}{\partial \Phi
_0\over
\partial x}({x+y\over 2})}e^{i(x-y)\cdot \theta /h}
\psi_0(x-y)p({x+y\over 2},\theta ;h)
u(y)dyd\theta  + R_1u(x),
$$
where $R_1={\cal O}(h^\infty ):L^2_{\Phi _0}\to L^2_{\Phi _0}$.
 We make a  contour deformation:
$$
\Gamma _t \defegal \set{ \theta ={2\over i}{\partial \Phi _0\over
\partial x}({x+y\over 2}) +it(\overline{x-y})},\ 0\le t\le t_0,\
t_0>0.
$$
Stokes' formula gives,
\begin{align*}
Pu(x)={}&{1\over (2\pi h)^n}\iint_{\Gamma _{t_0}}e^{{i\over
h}(x-y)\cdot \theta }\psi_0(x-y)p({x+y\over 2},\theta ;h)
u(y)dyd\theta \\
&+{1\over (2\pi h)^n}\iiint_{\Gamma _{[0,t_0]}}e^{{i\over
h}(x-y)\cdot \theta }u(y)\overline{\partial }_{y,\theta
}(\psi_0(x-y)p({x+y\over 2},\theta ;h))\wedge dy\wedge d\theta
+R_1u(x),
\end{align*}
where $\Gamma _{[0,t_0]}$ is the naturally
defined union of all the $\Gamma _t$ for $t\in [0,t_0]$. The
effective kernel of the first integral, viewed as an \op{} on
$L^2_{\Phi _0}$, is  ${\cal O}(h^{-n})e^{-{t_0\over h}\vert
x-y\vert ^2}$, which implies that this integral does indeed define
a \ufly{} \bdd{} \op{}: $L^2_{\Phi _0}\to L^2_{\Phi _0}$. The
effective kernel of the second integral can be estimated by a
constant times
\begin{align*}
\int_0^{t_0} h^{-n} e^{-{t\over h}\vert x-y\vert ^2}{\rm
dist\,}(({x+y\over 2},\theta ),\Lambda _{\Phi _0})^\infty dt
\hskip 2cm \cr ={\cal O}(1) \int_0^{t_0} h^{-n}e^{-{t\over h}\vert
x-y\vert ^2}(t\vert x-y\vert )^\infty dt={\cal O}(h^\infty ).
\end{align*}
We conclude that \ekv{7.11} { Pu(x)={1\over (2\pi
h)^n}\iint_{\theta ={2\over i}{\partial \Phi _0\over
\partial x}({x+y\over 2})+it_0\overline{(x-y)}}e^{{i\over h}(x-y)\cdot
\theta }\psi_0(x-y)p({x+y\over 2},\theta ;h)u(y)dyd\theta +R_2u,}
for $u\in H_{\Phi _0}$, where $R_2={\cal O}(h^\infty ):L^2_{\Phi
_0}\to L^2_{\Phi _0}$.

\bigskip

The aim of the next subsections is to introduce and study a new
strictly subharmonic function $\Phi_\eps$ related to $G$. As for
$\Phi_0$, the function  $\Phi_\eps$ is   associated to the space
$L^2_{\Phi _\eps}=L^2(\C ^n;e^{-2\Phi _\eps/h}L(dx))$ and its
subspace of entire functions $H_{\Phi_\eps }$. These spaces will
be used later  to get local resolvent estimates.

\subsection{Definition and derivative estimates of  $\Phi_\eps$}
 \label{Section5}

 Recall
that our weight \fu{} $G(\rho )$,  $\rho =(y,\eta )$ defined in
Proposition \ref{propescape}
 satisfies the estimates in the region $\delta\lambda \leq \oo(1)$
\ekv{5.1} { \nabla ^kG={\cal O}(\delta ^{(2-k)_+}),\ \delta (\rho
)\le \sqrt{h}, } \ekv{5.2} { \nabla ^kG={\cal O}(h(h\delta
)^{-{k\over 3}}),\ \delta (\rho )\ge \sqrt{h}. } It follows that
in the same region \ekv{5.3} { \nabla ^kG={\cal O}(hr^{-k}),}
where \ekv{5.4} { r(\rho ):=h^{1\over 3}(h^{1\over 2}+\delta (\rho
))^{1\over 3}. } Notice that \ekv{5.5} { h^{1\over 2}\le r\le
h^{1\over 2}+\delta , } so that $h^{1\over 2}+\delta (\rho )$ is
\ufly{} of constant order of magnitude in $B(\rho _0,{1\over
C_0}r(\rho _0))$ if $C_0>0$ is large enough and \indep{} of $\rho
_0$.

\par In $B(\rho _0,{1\over C_0}r(\rho _0))$ we introduce the scaled
variables $\widetilde{\rho }$, by \ekv{5.6} { \rho =\rho
_0+r_0\widetilde{\rho },\ r_0=r(\delta _0). } Then the scaled
\fu{} $G(\rho _0+r_0\widetilde{\rho })$ satisfies \ekv{5.7} {
\nabla _{\widetilde{\rho }}^k(G(\rho _0+r_0\widetilde{\rho
}))={\cal O}(h),\ |\widetilde{\rho }|<{1\over C_0}. }

\par Let
\ekv{5.11} { \Im (y,\eta )=\epsilon H_G(\Re (y,\eta )). } Then for
$\epsilon >0$ small enough, we have \ekv{5.12} { \kappa _T(\Lambda
_{\epsilon G})=\Lambda _{\Phi _\epsilon }\defegal \set{(x,\xi) \in
\C^{2n}; \ \ \xi ={2\over i}{\partial \Phi _\epsilon \over
\partial x}(x)}, } where $\Phi _\epsilon(x)$ is a critical value
w.r.t. $(y,\eta)$, \ekv{5.13} { \Phi _\epsilon (x)={\rm
v.c.}_{(y,\eta )\in\C ^n\times \R^n}(-\Im \phi (x,y)-(\Im y)\cdot
\eta +\epsilon G(\Re y,\eta )). } We note that, when $\epsilon
=0$, the unique critical point is \nondeg{}.

\par We are in the presence of the following general problem
 (where we
change and simplify the notation), namely to study the critical
value \ekv{5.14} { \Phi _\epsilon (x)={\rm v.c.}_{y}F_\epsilon
(x,y),\ x\in \R ^n,\,\, y\in \R ^n, } where $F_\epsilon
(x,y)$ is a smooth real-valued \fu{} such that \ekv{5.15} {
y\mapsto F_0(x,y)\hbox{ has a unique \nondeg{} critical point
}y_0(x), } \ekv{5.16} {
\partial _\epsilon ^2F_\epsilon (x,y)=0,
} \ekv{5.17} {
\partial _x^\alpha \partial _y^\beta \partial _\epsilon F_\epsilon
(x,y)={\cal O}(hr^{-\vert \beta \vert }), } where $r=h^{1\over
3}(h^{1\over 2}+\delta (y))^{1\over 3}$ and $\delta (y)\geq 0$ is
a Lipschitz \fu{}. From (\ref{5.16})--(\ref{5.17}) we see that
$$
\partial _yF_\epsilon -\partial _yF_0={\cal O}({h\epsilon \over r})\ll
\epsilon ,\ \partial _y^2F_\epsilon -\partial _y^2F_0={\cal
O}({h\over r^2}\epsilon )\ll 1,
$$
for $\epsilon \ll 1$. So, for $0\le \epsilon \le \epsilon _0\ll
1$, we see that $y\mapsto F_\epsilon (x,y)$ has a unique critical
point $y_\epsilon (x)$, depending smoothly on $(x,\epsilon )$.

\par In order to estimate the derivatives of $y_\epsilon (x)$
we work in
an $r_0$-\neigh{} of a variable point $(x_0,y_0)=(x_0,y_0(x_0))$,
$r_0=r(\delta (y_0))$, and put $x=x_0+r_0\widetilde{x}$,
$y_\epsilon
(x)=y_0(x_0+r_0\widetilde{x})+r_0\widetilde{y}_\epsilon
(\widetilde{x})$, with $\widetilde{y}_0(\widetilde{x})=0$ where we
hope that $\widetilde{y}_\epsilon ={\cal O}(\epsilon )$. Then
$\widetilde{y}_\epsilon (\widetilde{x})$ is the critical point of
\ekv{5.18} { \widetilde{y}\mapsto {1\over r_0^2}\big( F_\epsilon
(x_0+r_0\widetilde{x},y_0(x_0+r_0\widetilde{x})+
r_0\widetilde{y})-F_0(x_0+r_0\widetilde{x},y_0(x_0+r_0\widetilde{x}))
\big)
=:G_\epsilon (\widetilde{x},\widetilde{y}), } with \ekv{5.19} {
\partial _{\widetilde{x}}^\alpha \partial _{\widetilde{y}}^\beta
G_\epsilon ={\cal O}(1),\ \partial _{\widetilde{x}}^\alpha
\partial _{\widetilde{y}}^\beta \partial _\epsilon G_\epsilon
(\widetilde{x},\widetilde{y})={\cal O}({h\over r_0^2}), }
\ekv{5.20} {
\partial _{\widetilde{y}}G_0(\widetilde{x},0)=0,\ \vert \det \partial
_{\widetilde{y}}^2G_\epsilon \vert \ge 1/C. }

\par Introducing the rescaled parameter $\widetilde{\epsilon }$ by
$\epsilon ={r_0^2\over h}\widetilde{\epsilon }$, $\partial
_{\widetilde{\epsilon }}={r_0^2\over h}\partial _\epsilon $, we
have uniform bounds on all the derivatives $\partial
_{\widetilde{x}}^\alpha
\partial _{\widetilde{y}}^\beta
\partial _{\widetilde{\epsilon }}^\gamma
G_\epsilon $ while $\partial ^2_{\widetilde{y}}G_\epsilon $ is
\ufly{} \nondeg{}, and the same is therefore true about $\partial
_{\widetilde{x}}^\alpha \partial _{\widetilde{\epsilon }}^\gamma
\widetilde{y}_\epsilon (\widetilde{x})$, so
$$
\partial _{\widetilde{x}}^\alpha \partial _\epsilon ^\gamma
\widetilde{y}_\epsilon (\widetilde{x})={\cal O}(({h\over
r_0^2})^\gamma ),
$$
\ekv{5.21} {
\partial _x^\alpha \partial _\epsilon ^\gamma (y_\epsilon (x)-y_0(x))
={\cal
O}(r({h\over r^2})^\gamma r^{-\vert \alpha \vert }),\
r=r(x)=h^{1\over 3}(h^{1\over 2}+\delta (y_0(x)))^{1\over 3}). }
The critical value $G_\epsilon
(\widetilde{x},\widetilde{y}_\epsilon (x))$ also satisfies
$\partial _{\widetilde{x}}^\alpha \partial _{\widetilde{\epsilon
}}^\gamma (G_\epsilon (\widetilde{x},\widetilde{y}_\epsilon
(x)))={\cal O}(1),$ so \ekv{5.22} {
\partial _x^\alpha \partial _\epsilon ^\gamma
(F_\epsilon (x,y_\epsilon
(x))-F_0(x,y_0(x)))={\cal O}(r^2({h\over r^2})^\gamma r^{-\vert
\alpha \vert }). } We can Taylor expand this \wrt{} $\epsilon $
and get
$$
F_\epsilon (x,y_\epsilon (x))=F_0(x,y_0(x))+F_1(x)\epsilon
+F_2(x)\epsilon ^2+...+F_{N-1}(x)\epsilon ^{N-1}+R_{N}(x,\epsilon
)\epsilon ^N,
$$
where
\begin{align*}
F_1(x)={}&((\partial _\epsilon )_{\epsilon =0}F_\epsilon )
(x,y_0(x)),\\
\partial _x^\alpha F_k(x)={}&{\cal O}
(r^2({h\over r^2})^kr^{-\vert \alpha
\vert }),\,\, k\ge 1,\\
\partial _x^\alpha \partial _\epsilon ^\gamma R_N(x,\epsilon )
={}&{\cal
O}(r^2({h\over r^2})^{N+\gamma }r^{-\vert \alpha \vert }).
\end{align*}
Returning to (\ref{5.13}), we get \ekv{5.23} { \Phi _\epsilon
(x)=\Phi _0(x)+\Phi _1(x)\epsilon +...+\Phi _{N-1}(x)\epsilon
^{N-1}+R_N(x,\epsilon )\epsilon ^N, } where $\Phi _1,...,\Phi
_{N-1},R_N$ satisfy the same estimates and $$\Phi
_1(x)=G(y(x),\eta (x)),\ (y(x),\eta (x))=\kappa _T^{-1}(x,{2\over
i}{\partial \Phi _0\over \partial x}(x)).$$

\subsection{Study of $P$ as an operator on $H_{\Phi _\epsilon }$}

Recall that $H_{\Phi _\epsilon }$  is the subspace of entire
functions of  $L^2_{\Phi _\eps}=L^2(\C ^n;e^{-2\Phi
_\eps/h}L(dx))$. Since $\Phi _\epsilon -\Phi _0={\cal O}(\eps h)$,
we first notice that
$$
e^{-C\epsilon }\le \Vert u\Vert _{\Phi _\epsilon }/\Vert u\Vert
_{\Phi _0} \le e^{C\epsilon },
$$
and hence for instance
$$
R_1={\cal O}(1)e^{C\epsilon }h^\infty :L_{\Phi _\epsilon }^2\to
L^2_{\Phi _\epsilon }.
$$
Similarly, the effective kernel of the integral in (\ref{7.11}) as an
\op{}: $L^2_{\Phi _\epsilon }\to L^2_{\Phi _\epsilon }$ can be
estimated by
$$
{\cal O}(h^{-n})e^{-{t_0\over h}\vert x-y\vert ^2+2C\epsilon },
$$
corresponding to an \op{} of norm ${\cal O}(1)e^{2\epsilon
C}:L^2_{\Phi _\epsilon }\to L^2_{\Phi _\epsilon }$.

\par With the previous $t_0$ fixed, we now make the new contour
 deformation:
$$
\Gamma _t \defegal \set{\theta ={2\over i}{\partial \over \partial
x}((1-t)\Phi _0+t\Phi _\epsilon )({x+y\over
2})+it_0\overline{(x-y)}},\ \ \ 0\le t\le 1.
$$
Along this contour we have, using (\ref{5.22}), (\ref{5.5}):
$$
\overline{\partial }_{y,\theta }\psi_0(x-y)p({x+y\over 2},\theta
;h)={\cal O}(1)(\vert x-y\vert +\epsilon {h\over r({x+y\over
2})})^\infty \le {\cal O}(1)(\vert x-y\vert +\epsilon h^{1/
2})^\infty .
$$
By Stokes' formula, we see that \ekv{7.12} { Pu(x)={1\over (2\pi
h)^n}\iint_{\theta ={2\over i}{\partial \Phi_\eps \over
\partial x}({x+y\over 2})+it_0\overline{(x-y)}}e^{{i\over h}(x-y)\cdot
\theta }\psi_0(x-y)p({x+y\over 2},\theta ;h)u(y)dyd\theta +R_\eps
u, } for $u\in H_{\Phi _\epsilon },$ where \ekv{7.13} {R_\eps
={\cal O}(1)(e^{C\epsilon }h^\infty +h^\infty ):\, L^2_{\Phi
_\epsilon }\to L^2_{\Phi _\epsilon }.}

\subsection{Quantization vs. multiplication}

The aim of this short subsection is to check formula
(\ref{6.2ter}) below i.e. the equivalent of  \cite[formula
1.6]{Sjo90} for the Weyl quantization. Recall that the
$I$-lagrangian manifold $\Lambda_{\eps G}$ is defined by $
\Lambda_{\eps G} = \set{ \rho + i\eps H_G(\rho); \ \ \rho\in
\R^{2n}}$, and that  $G$ has bounded second derivatives. We also
have $\kappa_T(\Lambda_{\eps G}) = \Lambda_{\Phi_\eps} \defegal
\set{\xi = \xi_\eps (x)
\defegal \frac{2}{i} \frac{\D\Phi_\eps}{\D x}(x)}$. Notice that
the second derivatives of $\Phi_\eps$, and the first ones of
$\xi_\eps(x)$ are bounded. Recall  that $p$ is (an almost analytic
extension of) a $\cc^{\infty}$ symbol with all its derivatives
bounded. We get for $u\in H_{\Phi_\eps}$
$$
Pu(x)={1\over (2\pi h)^n}\iint_{\Gamma_\eps }e^{{i\over
h}(x-y)\cdot \theta }\psi_0(x-y)p({x+y\over
2},\theta)u(y)dyd\theta +R_\eps u,
$$
where $\Gamma_\eps  = \set{ \theta ={2\over i}{\partial \Phi_\eps
\over
\partial x}({x+y\over 2})+it_0\overline{(x-y)}}$ is the contour
 of integration
 and  $R_\eps = \oo(h^{\infty}) : L^2_{\Phi_\eps}
\longmapsto L^2_{\Phi_\eps}$.  Sometimes we omit the subscript
$\eps$.

We want to prove that for $h$ sufficiently small

\ekv{6.2ter} {(\chi Pu\vert u)_{H_{\Phi _\epsilon }} =
\int  p_\epsilon \vert u\vert ^2
\chi(x) e^{-2\Phi_\epsilon  (x)/h}L(dx)+{\cal O}(h)\Vert u\Vert ^2, }
where
$\chi$ has bounded derivatives (for example $\chi = 1$) and we define
$p_\epsilon =p(x,{2\over i}{\partial \Phi _\epsilon \over\partial
  x}(x))$ to be the restriction of $p$ to $\Lambda _{\Phi _\epsilon }$.
The proof  is a simple adaptation of the proof given in
\cite{Sjo90}. We first make the Taylor expansion of $p$,
\begin{equation} \label{pxy}
p(\frac{x+y}{2}, \theta) = p(x, \xi(x)) +\sum p^{(j)}(x,
\xi(x))(\theta_j - \xi_j(x)) + \sum p_{(j)}(x,
\xi(x))\sep{\frac{y_j-x_j}{2}} + r(x, y, \theta).
\end{equation}
On $\Gamma_\epsilon (x)$ we have $(\theta_j - \xi_j(x)) =
iC\overline{(x-y)}$, and $r(x, y, \theta) = \oo(|x-y|^2 +h^\infty )$.
 The effective
kernel of the operator $R$  corresponding to $r$ is therefore of
the form
\begin{equation}
\begin{split}
|R(x,y)|
    & = \oo \sep{  h^{-n}  e^{ - C|x-y|^2/2h } |x-y|^2
 \tilde{\psi}(x-y) } \\
    & = \oo \sep{  h^{-n}  e^{ - C|x-y|^2/2h } |x-y|^2/h}
h,
\end{split}
\end{equation}
for $C$ sufficiently large (since the second derivative of $\Phi$
is bounded). As a consequence
\begin{equation}
\begin{split}
\Re(Ru |u)_{\Phi _\epsilon } = \oo(h) \norm{u}_{\Phi_\epsilon }^2.
\end{split}
\end{equation}
For the contribution to (\ref{6.2ter}) from the second term of
(\ref{pxy}) we integrate by part  as in
\cite{Sjo90} and we see that this term is $\oo(h)\Vert u\Vert_{\Phi
  _\epsilon }^2$. For the
third term we simply write
\begin{equation*}\begin{split}
(2\pi h)^{-n} &\iint_{\Gamma}  e^{{i\over h}(x-y)\cdot \theta }
\sep{\sum p_{(j)}(x, \xi(x))\frac{y_j-x_j}{2}} u(y) dyd\theta \\ =
&\sum p_{(j)}(x, \xi(x))\sep{ x_j/2 -x_j/2} u(x) = 0.\end{split}
\end{equation*}
It follows that we have (\ref{6.2ter}).

\par Notice that we can take $\chi $ to be Lipschitz in
(\ref{6.2ter}),
and hence that relation can be iteratated to give:
\begin{equation} \label{mult2}
\begin{split}
( \chi Pu\vert Pu)_{\Phi _\epsilon }
   &  =   \int  |p_\epsilon |^2 \vert u\vert ^2 \chi(x)
    e^{-2\Phi_\epsilon
(x)/h}L(dx)+{\cal O}(h)\Vert u\Vert ^2.
\end{split}
\end{equation}

\Section{Local resolvent estimates for large $z$.}
 \label{locallarge}  \setcounter{equation}{0}

Again in this section we suppose that $p$ satisfies the hypothesis
\bf (H1) \rm  and that it is bounded with all its derivatives
outside a large compact set. The aim of this section is to get
resolvent estimates for functions localized near the critical
points on the FBI side, and for

\ekv{6.1} { h\ll \vert z\vert . }

 We realize $P$ as
an \op{} with leading symbol $p_\epsilon ={p_\vert}_{\Lambda
_{\epsilon G}}$ as $TPT^{-1}:H_{\Phi _\epsilon }\to H_{\Phi
_\epsilon }$, with $\Lambda _{\Phi _\epsilon }=\kappa _T(\Lambda
_{\epsilon G})$, and in the following we identify $P$ with
$TPT^{-1}$. We have seen that $\nabla ^2\Phi _\epsilon $ is
\ufly{} \bdd{} and consequently (see (\ref{6.2ter})) we have with
$\Phi =\Phi _\epsilon $ and scalar products and norms in
$L^2_{\Phi _\epsilon }$:
\ekv{6.2} { (\chi Pu\vert u)=\int p_\epsilon \vert u\vert ^2
\chi(x)  e^{-2\Phi (x)/h}L(dx)+{\cal O}(h)\Vert u\Vert ^2, } where
$p_\epsilon ={p_\vert}_{\Lambda _{\epsilon G}}$ is viewed as a
\fu{} on $\Lambda _{\Phi _\epsilon }$, and $\chi (x)\in C_0^\infty
(\C^n)$. We replace in this section the small parameter $h$
 in the construction of the function $G$ by $Ah$ where
  $A$ is some large constant. As a consequence for
$\eps$ fixed we get from Proposition \ref{propescape} that
$p_\epsilon (\rho )$ satisfies the estimates \ekv{6.3} { \Re
p_\epsilon (\rho )\ge {1\over C_0}\min (\delta (\rho
)^2,(Ah)^{2\over 3}\delta (\rho )^{2\over 3} ),} inside a large
compact set $K$ containing the support of $\chi$. From now on the
inequalities we give are to be understood in $K$. Note that
$C_0>0$ and the uniform estimate on $\nabla ^2\Phi _\epsilon $ do
not depend on $A$.

\par Let $\chi_0 \in C_0^\infty (\R,[0,1])$ be a standard cutoff to a
\neigh{} of $0\in \R$ and consider \ekv{6.4} {
\widetilde{p}_\epsilon (\rho )=p_\epsilon (\rho )+{1\over C_0}\min
(\vert z\vert ,(Ah)^{2\over 3}\vert z\vert ^{1\over 3})\chi_0
({\delta (\rho )^2\over \vert z\vert }).} Then there exists a
$C_1>0$ such that \ekv{6.5} { \Re \widetilde{p}_\epsilon (\rho
)\ge {1\over C_1}(\min (\delta (\rho )^2,(Ah)^{2\over 3}\delta
(\rho )^{2\over 3}))+\min (\vert z\vert,(Ah)^{2\over 3}\vert
z\vert ^{1\over 3})).}

Let us mention for further use that we can choose  the support of
$\chi_0 $ to be contained in a \sufly{} small \neigh{} of $0$, so
that

\ekv{6.6} { \vert p_\epsilon (\rho )-z\vert \ge \vert z\vert
/C_2,\hbox{ when }\chi_0 ({\delta (\rho )^2\over \vert z\vert
})\ne 0. }

 Write
 $$
 \Lambda^2 \defegal \min (\delta (\rho )^2,(Ah)^{2\over 3}\delta
(\rho )^{2\over 3}) ,
 \ \ \  \text { and } \ \ \
 Z
\defegal \min (\vert z\vert ,(Ah)^{2\over 3}\vert z\vert ^{1\over
3}).
$$
and denote
$$
 \chi _{\vert z\vert }(\rho ) \defegal
 \chi_0 ({\delta (\rho )^2 \over
\vert z\vert }) ,
$$
then (\ref{6.4}--\ref{6.5}) can be written as
$$
p_\eps + \frac{Z}{C_0} \chi_{\abs{z}}
 \geq \frac{1}{C_1}(\Lambda^2 + Z) .
$$
Considering $ \chi _{\vert z\vert }$ as a \fu{} of $x$ on the
FBI-\tf{} side, we get from (\ref{6.2})
\begin{equation}\label{6.7}\Re (\chi(P+{Z\over C_0} \chi
_{\vert z\vert }-z)u\vert u) + \oo(h) \norm{u}^2
\ge {1\over C_3}\sep{ \int
\chi \Lambda^2 \vert u\vert ^2e^{-2\Phi /h}L(dx) +Z (\chi u,u)},
\end{equation}
provided that $\chi $ is nonnegative and
(in addition to (\ref{6.1})): \ekv{6.8} { \Re
z\le Z/C_3. } Here $C_3>0$ is some \sufly{} large constant
which is \indep{} of
$A$, and $\chi _{\vert z\vert }$ in (\ref{6.7}) denotes the natural
multiplication \op{} on the FBI-side.

\par We shall combine (\ref{6.7}) with an estimate for
 $(\chi _{\vert
z\vert }u\vert u)$, that we shall obtain using the
ellipticity property
(\ref{6.6}). This will be obtained using an estimate
analogous to (\ref{6.2})
(that can also be found in \cite{Sjo90}) but since the support
of $\chi _{\vert
z\vert }$ may be very small we shall use a rescaling  which
also dilates the
Planck constant.

\begin{prop} \label{Prop6.1} Under the assumptions
  (\ref{6.1},\ref{6.8})
we have \ekv{6.21} { \Vert \chi _{\vert z\vert }u\Vert \le
C({1\over \vert z\vert }\Vert (P-z)u\Vert +\sqrt{{h\over
\min(1,\vert z\vert) }}\Vert u\Vert ), } for all $u\in H_{\Phi
_\epsilon }$.
\end{prop}

\preuve First assume $|z|\leq1$.
 Make the change of variables on the FBI-\tf{} side \ekv{6.9} {
x=\vert z\vert ^{1\over 2}\widetilde{x},\ hD_x=\vert z\vert
^{1\over 2}\widetilde{h}D_{\widetilde{x}},\ \widetilde{h}={h\over
\vert z\vert} .} Then, \ekv{6.10} { P(x,hD_x;h)-z=\vert z\vert
(\widetilde{P}(\widetilde{x},\widetilde{h}D_{\widetilde{x}};
\widetilde{h})-
\widetilde{z}), }
\ekv{6.11} { \widetilde{z}={z\over \vert z\vert
},\ \widetilde{P}(\widetilde{x},\widetilde{\xi
};\widetilde{h})={1\over \vert z\vert }P(x,\xi ;h)  ,\ (x,\xi
)=\vert z\vert ^{1\over 2}(\widetilde{x},\widetilde{\xi }). } If
$P(x,\xi ;h)=p(x,\xi )+hp_1(x,\xi )+h^2p_2(x,\xi )+...$, (where we
now consider the symbols in the complex domain), we see that
$$
\widetilde{P}(\widetilde{x},\widetilde{\xi };\widetilde{h})\sim
\sum_0^\infty \widetilde{p}_j(\widetilde{x},\widetilde{\xi
})\widetilde{h}^j,
$$
where $\widetilde{p}=\widetilde{p}_0={1\over \vert z\vert }p(\vert
z\vert ^{1\over 2}(\widetilde{x},\widetilde{\xi }))$,
$\widetilde{p}_j(\widetilde{x},\widetilde{\xi })=\vert z\vert
^{j-1}p_j(\vert z\vert ^{1\over 2}(\widetilde{x},\widetilde{\xi
})) $ are nice \bdd{} symbols, since $p(x,\xi )={\cal O}((x,\xi
)^2)$. Then using  (\ref{6.10}) \ekv{6.22} { {1\over \vert z\vert
}(P(x,hD_x;h)-z)=(\widetilde{P}(\widetilde{x},
\widetilde{h}D_{\widetilde{x}};\widetilde{h})-\widetilde{z}),\
\widetilde{h}= {h\over \vert z\vert }. } $L^2_\Phi $ \tf{}s into
$$
L^2_{\widetilde{\Phi }}=\{ \widetilde{u};\, \int \vert \widetilde{u}
(\widetilde{x})\vert ^2 e^{-2\widetilde{\Phi
}(\widetilde{x})/\widetilde{h}} L(d\widetilde{x})<\infty \} ,
$$
with the naturally associated norm and with $\widetilde{\Phi
}(\widetilde{x}) /\widetilde{h}=\Phi (x)/h$, so that
$$
\widetilde{\Phi }(\widetilde{x})=\Phi (\vert z\vert ^{1\over
2}\widetilde{x} )/\vert z\vert ,
$$
has a \ufly{} \bdd{} Hessian. Further, $\chi _{\vert z\vert }(x)=\chi
_1(\widetilde{x})$.

\par We have (omitting the Jacobians)
\begin{align*}\Vert {1\over \vert z\vert }(P-z)u\Vert _\Phi ^2={}
&\Vert
(\widetilde{P}-\widetilde{z})\widetilde{u}\Vert ^2_{\widetilde{\Phi }}
 \ge
\Vert \chi _1(\widetilde{P}-\widetilde{z})\widetilde{u}\Vert
_{\widetilde{\Phi }}^2\\
={}&\int \vert \chi _1(\widetilde{x})\vert ^2\vert
 \widetilde{ p_\epsilon \ }
-\widetilde{z}\vert ^2\vert \widetilde{u}\vert ^2
e^{-2\widetilde{\Phi
}/\widetilde{h}}L(d\widetilde{x})-{\cal O}(\widetilde{h})\Vert
\widetilde{u}\Vert _{\widetilde{\Phi }}^2\\
={}&\int\vert \chi _{\vert z\vert }(x)\vert ^2{1\over
\vert z\vert^2 }\vert
p_\epsilon -z\vert ^2 \vert u\vert ^2e^{-2\Phi /h}L(dx)-{\cal O}
({h\over
\vert z\vert })\Vert u\Vert _\Phi ^2\\
\ge{}& {1\over C}\Vert \chi _{\vert z\vert }u\Vert _\Phi ^2-{\cal
O}({h\over \vert z\vert })\Vert u\Vert _\Phi ^2.
\end{align*}
Here we used (\ref{mult2}) to obtain the second equality and
(\ref{6.6}) to get the last estimate.

In the case $|z| \geq 1$, we get more directly
\begin{equation}
\begin{split}
\Vert {1\over \vert z\vert }(P-z)u\Vert _\Phi ^2 & = \int\vert
\chi _{\vert z\vert }(x)\vert ^2{1\over \vert z\vert^2 }\vert
p_\epsilon -z\vert ^2 \vert u\vert ^2e^{-2\Phi /h}L(dx)+{\cal O}
({h})\Vert u\Vert _\Phi ^2\\
& \geq {1\over C}\Vert \chi _{\vert z\vert }u\Vert _\Phi ^2-{\cal
O}({h })\Vert u\Vert _\Phi ^2.
\end{split}
\end{equation}
This completes the proof of Proposition \ref{Prop6.1}. \fin

\par We can therefore write
\begin{align*}
&\Re (\chi(P-z)u\vert u)+{1\over C_0}Z(\chi \chi _{\vert z\vert
}u\vert u) \\ & \hskip 1cm \le  \Vert (P-z)u\Vert \Vert u\Vert
+{C\over C_0}Z \sep{ {1\over \vert z\vert }\Vert (P-z)u\Vert
+\sqrt{h\over \min(1,\vert z\vert) }\Vert u\Vert }\Vert u\Vert +
\oo(h)\norm{u}^2.
\end{align*}
Combining this with (\ref{6.7}), we get
$$
{Z\over C_2}(\chi u | u)  \le (1+{C\over C_0})\Vert (P-z)u\Vert
\norm{u}  +{C\over C_0}\sqrt{h\over\min(1,\vert z\vert) }Z\Vert
u\Vert + \oo(h)\norm{u}^2.
$$
and writing $\chi = 1 + (\chi-1)$ yields
$$
{Z\over C_2}\norm{u}^2  \le (1+{C\over C_0})\Vert (P-z)u\Vert
\norm{u} +{C\over C_0}\sqrt{h\over \min(1,\vert z\vert) }Z\Vert
u\Vert + \oo(h)\norm{u}^2 + C Z \norm{(1-\chi)u}\norm{u}
$$
 Assuming
$h/\min(1,\vert z\vert)$ \sufly{} small \indep{}ly of $A$, we get
the main result of this section: \ekv{6.23} { Z \Vert u\Vert \le
{\cal O}(1)\sep{ \Vert (P-z)u\Vert + Z \norm{(1-\chi)u} } , }
where we recall the assumptions (\ref{6.1}) and (\ref{6.8}) on
$z$.

\Section{The quadratic case}\label{SectQuad} The main
 purpose of this first
section is to get resolvent estimates  for  operators with
quadratic symbol. The main reference  for this is \cite{Sjo74},
and  all the computations are explicit. In the special case of
the quadratic Kramers-Fokker-Planck operator, the form of the
spectrum is well known (see for example \cite{Ri89}) and we
compute it explicitly in section \ref{secFP}.

\subsection{Sectorial property in a linear weighted space
 and applications}

Let $P_0$ be a quadratic operator in the sense that
the symbol $p = p_1+
ip_2$ is a complex-valued quadratic form and assume that the symbol
satisfies $p_1\ge 0$ and a subelliptic estimate
\begin{equation} \label{mainhyp4}
p_1 + \eps_0 H_{p_2}^2 p_1 \geq \frac{\eps_0}{C} d_0^2 ,
\end{equation}
where $d_0(\rho) = |\rho|^2$. Note that this implies that $p$ has
$0$ as unique critical point.

Now we use the  weight
$$
G^0= G_T,
$$
introduced in (\ref{GT}) near the critical point, and use the
definition there in the whole space.  Since $p$ is quadratic, so
is  $G^0$, and we have for $0<\eps\leq \eps_0$
\begin{equation} \label{mainhyp44}
p_1 + \eps H_{p_2}G^0 \geq \frac{\eps}{C} d_0^2 ,
\end{equation}
 As in section \ref{secFBI} we use the
global FBI transform with quadratic phase $\phi$
$$
 Tu(x)=Ch^{-{3n \over
4}} \int e^{{i\over h} \phi (x,y)}u(y)dy.
$$
 The  canonical transformation
associated with the FBI transform $T$ is given by $  \kappa
_T:(y,-\partial _y\phi (x,y))\mapsto (x,\partial _x\phi (x,y)) $
and  we define  $ \Lambda _{\Phi _0}=\kappa _{T}(\R^{2n})$ and $
\Phi _0(x)=-\Im \phi (x,y_0(x)), $ where $y_0(x)$ is the point
where $\R^n\ni y\mapsto -\Im \phi (x,y)$ takes its non-degenerate
maximum.

We define
\ekv{gzero}{
\C^{2n} \supset \Lambda_{\eps G^0} \defegal  \set{ (y, \eta); \
\Im (y,\eta )=\epsilon H_{G^0}(\Re (y,\eta ))}
}
and for $\eps$ small enough we check that
$$
\kappa _T(\Lambda _{\epsilon G^0})=\Lambda _{\Phi^0_\epsilon }
\defegal \set{ (x, \xi) \,\, \xi ={2\over i}
{\partial \Phi _\epsilon \over \partial x}(x)},
$$
where $ \Phi^0_\epsilon$ is defined using the following procedure:
the function
$$
F^0_\eps (x,y,\eta) = -\Im \phi (x,y)-(\Im y)\cdot \eta +\epsilon
G^0(\Re y,\eta )
$$
 is quadratic and when $\epsilon =0$ it has a  unique
\nondeg{} critical point for $x$ fixed. By homogeneity, this is
also the case for $F^0_\eps$. The unique critical point
$(y_\eps(x), \eta_\eps(x))$ depends linearly on $x$ and smoothly
on $\eps$. We finally  write
$$
\Phi^0_\epsilon (x)={\rm v.c.}_{(y,\eta )\in{\bf
C}^n\times \R^n}(-\Im \phi (x,y)-(\Im y)\cdot \eta +\epsilon
G^0(\Re y,\eta ))
$$

From now on we work entirely on the FBI side, denoting by $u$ a
function on the FBI side (instead of $Tu$), and  by the same
letter
 $P_0$ the (unbounded) operator on $L^2_{\Phi_0}$
$$
P_0 u(x)={1\over (2\pi h)^n}\iint_{\theta ={2\over i}{\partial
\Phi _0\over
\partial x}({x+y\over 2})} e^{{i\over h}(x-y)\cdot \theta }
p({x+y\over 2},\theta )u(y)dyd\theta.
$$
Since the symbol of $P_0$ is quadratic, it is holomorphic and we
also have the following formula for $P_0$ as an unbounded operator
on $L^2_{\Phi_0}$ :
$$
P_0 u(x)={1\over (2\pi h)^n}\iint_{\theta ={2\over i}{\partial
\Phi _0\over
\partial x}({x+y\over 2}) + it_0\overline{(x-y)} }
e^{{i\over h}(x-y)\cdot \theta } p({x+y\over 2},\theta
)u(y)dyd\theta.
$$

We can now make a new contour deformation, and we get an unbounded
operator again denoted $P_0$ on the space $L^2_{\Phi^0_\eps}$
naturally associated to $\Phi^0_\eps$:
$$
P_0 u(x)={1\over (2\pi h)^n}\iint_{\theta ={2\over i}{\partial
\Phi^0_\eps \over
\partial x}({x+y\over 2}) + it_0\overline{(x-y)} }
e^{{i\over h}(x-y)\cdot \theta }
p({x+y\over 2},\theta )u(y)dyd\theta.
$$
Of course coming back to the real side by the FBI transform, $P_0$
can be viewed
 as an unbounded operator on $L^2(\R^n)$ with symbol
 $$
 \tilde{p} = p(\rho + i\eps H_{G^0}) ,
$$
and here the symbol of $P_0$ is quadratic and
satisfies
\begin{equation}
\begin{split}
\widetilde{p}_1 & =p_1(\rho )+\epsilon H_{p_2}G^0(\rho )+{\cal
O}(\epsilon ^2\vert \nabla G^0\vert ^2),  \\
        \widetilde{p}_2&=p_2(\rho )-\epsilon
H_{p_1}G^0(\rho )+{\cal O}(\epsilon ^2\vert \nabla G^0\vert ^2).
\end{split}
\end{equation}
Now each term is quadratic therefore using the homogeneity,
(\ref{mainhyp44}), and choosing  $\epsilon
>0$ small enough yields
\begin{equation*} \label{p1leqh12ZZZ}
\widetilde{p}_1 \geq \frac{\eps}{C} d_0^2,\ \widetilde{p}_2={\cal
O}(d_0^2).
\end{equation*}
 In particular  $\widetilde{p}$ takes its values in an angle
around the positive real axis,
 $ \widetilde{p}_1\ge \eps |\widetilde{p}_2|/C.$
As a consequence we can apply to $P_0$ as an unbounded operator
on $L^2_{\Phi^0_\eps}$ the result of \cite[Theorem 3.5]{Sjo74},
which gives with $d(x)=|x|$~:

\begin{prop} \label{quad}  Consider  $P_0$ as an operator on
$H_{\Phi^0_\eps}$. Then
\begin{description}
\item{a)} the spectrum of $P_0$ is a set $\set{\mu_l}$ given by
$$
 \set{ \frac{h}{i}
\sum_{\Im \lambda_j
>0} \sep{ {1 \over 2} + k_j }\lambda_j; \ \ \ \lambda_j \in Sp( F
), \ \ \ k_j \in \N} ,
$$
where the $\lambda_j$s are the eigenvalues, repeated with their
multiplicities, of the fundamental matrix $F$  of
$\operatorname{Hess}p$.
\item{b)} Let   $z$ vary in a compact set $K\subset \C$
disjoint from the union of
the $\mu_js$, then
\begin{equation} \label{7.19}
\Vert (h+d^2)u\Vert \le C\Vert (P_0-hz)u\Vert ,\ \Vert
(h+d^2)^{1\over 2}u\Vert \le C\Vert (h+d^2)^{-{1\over 2}}
(P_0-hz)u\Vert
\end{equation}
where $d(x)=\vert x\vert $ (essentially equal to $|\rho|^2 $ if we
lift it to $\Lambda _\Phi $), and for $u$ \hol{} with $(h+d^2)u\in
L^2_{\Phi ^0_\eps}$ and $(h+d^2)^{1\over 2}u\in L^2_{\Phi
^0_\eps}$ respectively.
\end{description}
\end{prop}

Recall that the fundamental matrix of the quadratic form $p$ is
the matrix of the (linearized) Hamilton flow and is given  by
$$
F =\left( \begin{array}{cc}
                       p''_{\xi ,x} & p''_{\xi,\xi} \\
                     - p''_{x,x} & - p''_{x,\xi}
\end{array} \right)
$$

\preuve This follows from  \cite[Theorem 3.5]{Sjo74} and some
remarks: First we note that the presence of the small parameter
$h$ is easy to deal with since $P_0$ is linearly   conjugated with
$h(P_0|_{h=1})$ by the symplectic change of  coordinates $(x,\xi)
\rightarrow (h^{1 \over 2}x, h^{-{1 \over 2}} \xi)$. We also
notice that the eigenvalues of the fundamental matrix $F$ of $p$
are the same as the ones of the fundamental matrix $\tilde{F}$ of
$\tilde{p}$, also by a symplectic change of variables.  Point $b)$
of the proposition is a direct consequence of \cite[Theorem
3.5]{Sjo74} and the change of symplectic coordinates $(x,\xi)
\rightarrow (h^{1 \over 2}x, h^{-{1 \over 2}} \xi)$. \fin

 In Section \ref{secFP} we shall explicitly compute the
  eigenvalues in
the case of the Kramers-Fokker-Planck operator. In the next
subsection we shall compare an operator $P$ with its quadratic
approximation near its critical points: In order to get a global a
 priori estimate for $P-hz$, we will need a truncated version of
(\ref{7.19}).

\subsection{Localized resolvent estimates}

\par Let $\chi_0 \in C_0^\infty (\C^n)$, $\chi_0 =1$ near $x=0$.
We fix $\eps>0$ small and write in this subsection  $\Phi ^0$ instead
of $\Phi ^0_\eps$. The simple idea is to apply (\ref{7.19}) with
$u$ replaced by $\chi_0 u$ and then try to estimate the commutator
$[P_0,\chi_0 ]u$. However, $\chi_0 u$ is not \hol{}, so we will
replace $\chi_0 u$ by $\Pi \chi_0 u$, where $\Pi :L^2_{\Phi ^0}\to
H_{\Phi ^0}$ is the \og{} projection.

The main result of this subsection is
\begin{prop} \label{Prop7.1} Let $\chi_0 \in C_0^\infty (\C^n)$
be fixed and equal to $1$ near 0, and fix $k\in \R$. Then for $z$
varying in a compact set which does not contain any \ev{}s of
$P_0|_{h=1}$, we have \ekv{7.37} { \Vert (h+d^2)^{1-k}\chi_0
u\Vert \le C\Vert (h+d^2)^{-k}\chi_0 (P_0-hz)u\Vert +{\cal
O}(h^{1\over 2})\Vert 1_Ku\Vert , } where $K$ is any fixed
\neigh{} of ${\rm supp}(\nabla \chi_0 )$.
\end{prop}

We need a series of technical preparations.

\ititem{Estimates for $[P_0,\chi_0 ]$.}  We have
$$
[P_0,\chi ]=\sum_{\vert \alpha +\beta \vert =1}h\chi _{\alpha
,\beta }(x)x^\alpha (hD_x)^\beta +h^2\chi _{0,0}(x),
$$
where $\chi
_{\alpha ,\beta }\in C_0^\infty (\C^n)$, ${\rm supp\,}(\chi
_{\alpha ,\beta })\subset {\rm supp\,}\nabla \chi_0 $. We can
conclude that

\ekv{7.20} { \Vert [P_0,\chi ]u\Vert \le Ch\Vert 1_Ku\Vert , }
where $C$ depends on $\chi $ and $K$ is an \ably{} small \neigh{}
of ${\rm supp\,}(\nabla \chi )$. Here we also use that $\Vert
1_{{\rm supp\,}\chi }(hD)^\alpha u\Vert \le C\Vert 1_K u\Vert $,
if  $u$ is \hol{} near $K$.

\ititem{Estimates for $[P_0,\Pi ]$.} Recall (from e.g.
\cite{Sjo95}) that $\Pi $ is given by \ekv{7.21} { \Pi
u(x)=Ch^{-n}\int e^{{2\over h}(\Psi ^0(x,\overline{y})-\Phi
^0(y))}u(y)L(dy), } where $\Psi ^0(x,z)$ is the unique (second
order) \hol{} \pol{} on $\C^{2n}$ with $\Psi
^0(x,\overline{x})=\Phi ^0(x)$. Notice that \ekv{7.22} { (\partial
_x\Psi ^0)(x,\overline{x})=\partial _x\Phi ^0(x), } and recall the
well known fact that \ekv{7.23} { 2\Re \Psi
^0(x,\overline{y})-\Phi ^0(x)-\Phi ^0(y)\sim -\vert x-y\vert ^2. }

\par For $\vert \alpha +\beta \vert \le 2$, we get
 by integration by parts,
\[
\begin{split}
  [x^\alpha (hD_x)^\beta ,\Pi ]u(x)
  ={}&Ch^{-n}\int (x^\alpha (hD_x)^\beta
-(-hD_y)^\beta y^\alpha ) e^{{2\over h}(\Psi
^0(x,\overline{y})-\Phi
^0(y))}u(y)L(dy)\\
 ={}&Ch^{-n}\int e^{{2\over h}(\Psi ^0(x,\overline{y})-\Phi
^0(y))}a_{\alpha ,\beta }(x,y;h) u(y)L(dy),
\end{split}
\]
where
$$
a_{\alpha ,\beta }=(x^\alpha (hD_x+{2\over i}\partial _x\Psi
^0(x,\overline{y}))^\beta -(-hD_y+{2\over i}\partial _y\Phi
^0(y))^\beta \circ y^\alpha )(1).$$ Using (\ref{7.22}), we see
that \ekv{7.24} { a_{\alpha ,\beta }=\left\{ \begin{array}{ll} 0,\
& \vert \alpha +\beta \vert =0,\\ b_1(x,y),\ & \vert \alpha +\beta
\vert =1,\\ b_2(x,y)+hb_0,\ & \vert \alpha +\beta \vert
=2,\end{array} \right. } where $b_j$ is a homogeneous \pol{} of
degree $j$, vanishing on the diagonal when $j=1,2$.

\par Relation
(\ref{7.23}) implies that the effective kernel of
$\Pi :L^2_{\Phi ^0}\to
L^2_{\Phi ^0}$ is ${\cal O}(h^{-n}e^{-\vert x-y\vert ^2/(Ch)})$,
so
$$
\Vert [x^\alpha (hD)^\beta ,\Pi ]u\Vert \le {\cal O}(h^{1\over
2})\times \left\{ \begin{array}{ll}  \Vert u\Vert ,\ &\vert \alpha
+\beta \vert =1,\\ \Vert (h+d^2)^{1\over 2}u\Vert ,\ &\vert \alpha
+\beta \vert =2. \end{array} \right.
$$
It follows that \ekv{7.25} { \Vert [P_0,\Pi ]u\Vert \le {\cal
O}(h^{1\over 2})\Vert (h+d^2)^{1\over 2}u\Vert , } since in the
case of $P_0$, we do not have to consider any commutators with
$x^\alpha (hD)^\beta $ with $\vert \alpha +\beta \vert =1$. The
standard inequality
$$
{1+\vert x\vert \over 1+\vert y\vert }\le 1+\vert x-y\vert ,
$$
implies that
$$
{h^{1\over 2}+d(x)\over h^{1\over 2}+d(y)}\le 1+{\vert x-y\vert
\over \sqrt{h}}\le C_\epsilon e^{\epsilon \vert x-y\vert ^2/h},
$$
for every $\epsilon >0$. It is therefore clear that we can
conjugate $[P_0,\Pi ]$ in (\ref{7.25}) by any power of $h^{1\over
2}+d$. Indeed, the proof there shows that the effective kernel of
$[P_0,\Pi ](h^{1\over 2}+d)^{-1}$ is ${\cal O}(1)h^{{1\over 2}-n}
e^{-\vert x-y\vert ^2/(Ch)}$. Hence \ekv{7.26} { \Vert (h^{1\over
2}+d)^{-k}[P_0,\Pi ]u\Vert \le {\cal O}(h^{1\over 2}) \Vert
(h^{1\over 2}+d)^{1-k}u\Vert , } for every $k\in\R$.

\ititem{Estimates for $1-\Pi $.}  We briefly recall
H{\"o}rmander's $L^2$-method for the $h\overline{\partial }$-complex,
following \cite{Sjo95},
$$
C_0^\infty (\C^n)\to C_0^\infty (\C^n;\wedge ^{0,1} \C^n)
\to C_0^\infty (\C^n;\wedge ^{0,2}\C^n)\to ...\to
C_0^\infty (\C^n;\wedge ^{0,n}\C^n).
$$
 We have here the natural
Hilbert space norms induced by the weight $e^{-2\Phi ^0/h}L(dx)$.
Equivalently, we consider the conjugated complex $e^{-\Phi
^0/h}h\dbar{}e^{\Phi ^0/h}=h\dbar{}+\dbar (\Phi ^0)^\wedge $ in
the standard $L^2$-spaces. The adjoint of the last complex is then
given by $h\dbar^*+\partial (\Phi ^0)^\rfloor$. More explicitly,
$$h\dbar +\dbar(\Phi ^0)^\wedge =\sum Z_jd\overline{z}_j^\wedge ,\quad
h\dbar^* +\partial(\Phi ^0)^\rfloor =\sum Z_j^*dz_j^\rfloor ,$$
where $Z_j=h\partial _{\overline{z}_j}+\partial
_{\overline{z}_j}\Phi ^0$. The corresponding Hodge Laplacian is
then
\begin{align*}
&(h\dbar{}+ \dbar (\Phi ^0)^\wedge )  (h\dbar^*+\partial (\Phi
^0)^\rfloor ) + (h\dbar^*+\partial (\Phi ^0)^\rfloor )   (
h\dbar{}+\dbar (\Phi ^0)^\wedge )\\ = {}
 &
\sum_{j,k}(Z_jZ_k^*\otimes d\overline{z}_j^\wedge dz_k^\rfloor
+Z_k^* Z_j \otimes dz_k^\rfloor d\overline{z}_j^\wedge )=
(\sum_jZ_j^*Z_j)\otimes 1+h\sum_{j,k}2\partial
_{\overline{z}_j}\partial _{z_k}\Phi ^0 d\overline{z}_j^\wedge
dz_k^\rfloor ,
\end{align*}
where we used that $[Z_j,Z_k^*]=2h\partial
_{\overline{z}_j}\partial _{z_k}\Phi ^0$ and the standard
identity, $d\overline{z}_j^\wedge dz_k^\rfloor +dz_k^\rfloor
d\overline{z}_j^\wedge =\langle dz_k\vert d\overline{z}_j\rangle
=\delta _{j,k}$. In particular, the Hodge Laplacian
$$
\Delta _1=h\dbar ^* h\dbar+h\dbar\, h\dbar ^*
$$
on $(0,1)$-forms can be identified with \ekv{7.27} {
\widetilde{\Delta }_1=(\sum Z_j^*Z_j)\otimes 1_{{\bf
C}^n}+2h(\partial _{\overline{z}_j}\partial _{z_k}\Phi ^0 ), }
acting on $L^2(\C^n;\C^n)$. The strict \plsh{}ity of
$\Phi ^0$ means that the Hermitian matrix appearing in the last
term in (\ref{7.27}) is $\ge 1/C$, and hence we get the apriori
estimate (now using for a while ordinary $L^2$-norms): \ekv{7.28}
{ {h\over C}\Vert u\Vert ^2+\sum \Vert Z_ju\Vert ^2\le
(\widetilde{\Delta }_1u\vert u), } leading first to \ekv{7.29} {
h\Vert u\Vert \le C\Vert \widetilde{\Delta }_1u\Vert , } and then
to \ekv{7.30} { h^{1\over 2}\Vert Z_ju\Vert \le C\Vert
\widetilde{\Delta }_1u\Vert . } We can also write $\sum
Z_j^*Z_j=\sum Z_jZ_j^*+{\cal O}(h)$, so $\widetilde{\Delta
}_1=\sum Z_jZ_j^*+{\cal O}(h)$, and hence
$$
(\widetilde{\Delta }_1u\vert u)\ge \sum \Vert Z_j^*u\Vert
^2-Ch\Vert u\Vert ^2,
$$
which together with (\ref{7.28}) implies first \ekv{7.30.5} {
h\Vert u\Vert ^2+\sum \Vert Z_ju\Vert ^2+\sum \Vert Z_j^*u\Vert
^2\le C(\widetilde{\Delta }_1u\vert u), } and then \ekv{7.31} {
h^{1\over 2}\Vert Z_j^*u\Vert \le C \Vert \widetilde{\Delta
}_1u\Vert . }

\par We also need to check that these estimates remain valid after
conjugation of $\widetilde{\Delta }_1$ by any power of $h+d^2$ or
equivalently by any power of $\lambda h+d^2$, where $\lambda \gg
1$ is \indep{} of $h$. This will follow from the following
observations:
\smallskip
\par\noindent 1)
$$
(\lambda h+d^2)^{-k}Z_j(\lambda h+d^2)^k=Z_j+k{h\partial
_{z_j}(d^2)\over \lambda h+d^2},
$$
and
$$
\abs{ {h\partial _{z_j}(d^2)\over \lambda h+d^2}} \le {Chd\over
\lambda h+d^2}\le {C\over \lambda ^{1\over 2}}{(\lambda h)^{1\over
2}d\over \lambda h+d^2}h^{1\over 2}\le \alpha (\lambda )h^{1\over
2},
$$
where $\alpha (\lambda )\to 0$ when $\lambda \to \infty $. A
similar remark holds for $(\lambda h+d^2)^{-k}Z_j^*(\lambda
h+d^2)^k$.
\smallskip

\par\noindent 2) We have
\begin{align*}
\widehat{\Delta }_1:={}&(\lambda h+d^2)^{-k}\widetilde{\Delta }_1
(\lambda h+d^2)^k\\
={}&\sum (Z_j^*+o(1)h^{1\over 2})(Z_j+o(1)h^{1\over 2})+2h(\partial
_{\overline{z}_j}\partial _{z_k}\Phi ^0),\end{align*} where
$o(1)$ refers to the limit $\lambda \to \infty $. Thus
\begin{align*}
\Re (\widehat{\Delta }_1u\vert u)\ge{}& {h\over C}\Vert u\Vert
^2+\sum \Vert Z_ju\Vert ^2 -o(1)h^{1\over 2}\Vert u\Vert \Vert
Z_ju\Vert
-o(1)(hu\vert u)\\
\ge{}& {h\over 2C}\Vert u\Vert ^2+{1\over 2}\sum \Vert Z_ju\Vert
^2,
\end{align*}
\ekv{7.33} { h\Vert u\Vert ^2+\sum \Vert Z_ju\Vert ^2\le C\Re
((\widehat{\Delta }_1)u\vert u). }
\smallskip
Then do as before with $\widetilde{\Delta }_1$ replaced by $\Re
\widehat{\Delta }_1$, to get \ekv{7.34} { h\Vert u\Vert ^2+\sum
\Vert Z_ju\Vert ^2+\sum \Vert Z_j^*u\Vert ^2\le C\Re
(\widehat{\Delta }_1u\vert u). }
 Back to the original $\Delta _1$ we thus have (with the norms
now being those of $L^2_{\Phi ^0}$): \ekv{7.34.5} { \Vert \Delta
_1^{-1}\Vert \le {\cal O}({1\over h}),\ \Vert h\dbar^*\Delta
_1^{-1}\Vert \le {\cal O}({1\over \sqrt{h}}), } as well as the
same estimates for
$$
(h+d^2)^k\Delta _1^{-1}(h+d^2)^{-k},\ (h+d^2)^kh\dbar^* \Delta
_1^{-1}(h+d^2)^{-k}.
$$
 Now use the fact, that
\ekv{7.35} { 1-\Pi=h\dbar^* \Delta _1^{-1}h\dbar , } to conclude
that if $\chi_0 \in C_0^\infty (\C^n)$ is fixed and equal to
1 near 0, and $u$ is \hol{} near ${\rm supp\,}\chi_0 $, then
$$
(1-\Pi )(u\chi_0 )=h\dbar^*\Delta _1^{-1}(u(h\dbar \chi_0 ))
$$
satisfies \ekv{7.36} { \Vert (h+d^2)^k(1-\Pi )(u\chi_0 )\Vert \le
C_kh^{1\over 2}\Vert u\dbar \chi_0 \Vert . } Recall that here and
until further notice the norms are those of $L^2_{\Phi
^0}$.\smallskip

\par Let $\chi_0 \in C_0^\infty (\C^n)$ be fixed and $=1$ near 0.
 Recall
that $z$ varies in a compact set which does not contain any \ev{}s
of $(P_0)_{h=1}$.

\preuve[of Proposition \ref{Prop7.1}]  We start from (\ref{7.19}):
\ekv{7.38} { \Vert (h+d^2)^{1-k}u\Vert \le
C\Vert (h+d^2)^{-k}(P_0-hz)u\Vert , } for $u$ \hol{} with
$(h+d^2)^{1-k}u\in L^2_{\Phi ^0}$. Replace $u$ by $\Pi \chi_0 u$:
$$
\Vert (h+d^2)^{1-k}\Pi \chi_0 u\Vert \le C\Vert
(h+d^2)^{-k}(P_0-hz)\Pi \chi_0 u\Vert .
$$
It follows that
\begin{equation}\label{7.39}
\begin{split}
\Vert (h+d^2)^{1-k}\chi_0 u\Vert \le{}& \Vert (h+d^2)^{1-k}\Pi
\chi_0 u\Vert +\Vert (h+d^2)^{1-k}(1-\Pi
)\chi_0 u\Vert \\
\le{} & C\Vert (h+d^2)^{-k}(P_0-hz)\Pi \chi_0 u\Vert +{\cal
O}(h^{1\over 2})
\Vert  u\dbar \chi_0 \Vert .
\end{split}
\end{equation}
where we used (\ref{7.36}) and the fact that $h+d^2\sim 1$ on
${\rm supp\,} \dbar \chi_0 $.

\par Here
\begin{equation}\label{7.40}
\begin{split}
&\Vert (h+d^2)^{-k}(P_0-hz)\Pi \chi_0 u\Vert \\
\le{} & \Vert (h+d^2)^{-k}
\Pi \chi_0 (P_0-hz)u\Vert +\Vert (h+d^2)^{-k}[P_0,\Pi \chi_0 ]u\Vert \\
\le{} & C\Vert (h+d^2)^{-k}\chi_0 (P_0-hz)u\Vert +\Vert
(h+d^2)^{-k}[P_0,\Pi \chi_0 ]u\Vert .
\end{split}
\end{equation}
Now
\begin{equation}\label{7.41}
\begin{split}
[P_0,\Pi \chi_0 ]u={}&[P_0,\Pi ]\chi_0 u+\Pi [P_0,\chi_0 ]u\\
={}&[P_0,\Pi ]\Pi \chi_0 u+[P_0,\Pi ](1-\Pi )\chi_0 u+\Pi
[P_0,\chi_0
]u\\
={}&[P_0,\Pi ](1-\Pi )\chi_0 u+\Pi [P_0,\chi_0 ]u,
\end{split}
\end{equation}
where we used that $[P_0,\Pi ]\Pi =0$, since $P_0$ conserves
\hol{} \fu{}s. Combining (\ref{7.41}), (\ref{7.26}), (\ref{7.36}),
(\ref{7.20}), we see that
\ekv{7.42} { \Vert (h+d^2)^{-k}[P_0,\Pi
\chi_0 ]u\Vert \le {\cal O}(h)\Vert 1_Ku\Vert . }
Combining this with (\ref{7.39}),
(\ref{7.40}), (\ref{7.42}), we get (\ref{7.37}).\fin

\remark \label{remarknormequ} In Proposition
\ref{Prop7.1} we can replace the norm $L^2_{\Phi ^0_\eps}$ by
$L^2_{\Phi_0}$ or any other norm which is equivalent to the $L^2_{\Phi
^0_\epsilon }$ norm for functions with support near $K$.

\Section{Local resolvent estimate for small $z$.}
\label{localsmall}  \setcounter{equation}{0}

Again in this section we suppose that $p$ satisfies the hypothesis
\bf (H1) \rm  and that it is bounded with all its derivatives
outside a large compact set $K$. We also replace for a while the
small parameter $h$ by $Ah$ in the construction of $G$,
where $A$ is some large constant, and
work in $K$.

 Recall that $G=G_{Ah}$ satisfies the estimates:
\ekv{7.1} {\nabla G={\cal O}(\delta ^{2-k)_+}),\ \delta (\rho )\le
\sqrt{Ah},} \ekv{7.2} { \nabla ^kG={\cal O}(Ah(Ah\delta
)^{-k/3}),\ \delta (\rho )\ge \sqrt{Ah}, } implying, \ekv{7.3} {
\nabla ^kG={\cal O}(Ah((Ah)^{1\over 3}(Ah+\delta ^2)^{1\over
6})^{-k})={\cal O}(Ahr^{-k}), } \ekv{7.4} { r(\rho ):=
(Ah)^{1\over 3}(Ah+\delta ^2)^{1\over 6}. } Writing
$${p_\vert}_{\Lambda _{\epsilon G}}=p_\epsilon =p_1+ip_2,$$
we recall that in $K$ \ekv{7.5} { p_1\ge {\epsilon \over C}\min
(\delta (\rho )^2,(\delta Ah)^{2\over 3}), } \ekv{7.6} { p_2={\cal
O}(\delta ^2). }

\par We represent $\Lambda _{\epsilon G}$ on the FBI-\tf{} side by
$$
\xi ={2\over i}{\partial \Phi _\epsilon \over \partial x}(x),\ \Phi
_\epsilon =\Phi _0+\epsilon \widetilde{G}(x;h),
$$
where $\widetilde{G}$ has the same properties as $G$ (cf
(\ref{5.23})).
We also know that
$\widetilde{G}$ and $\Phi _\epsilon $  are \indep{} of $h$
in the region
$\vert x\vert \le \sqrt{Ah}$. From now on $\epsilon >0$ will
 be small and
fixed.

\par Assume for simplicity that ${\cal C}$ consists of just one point,
corresponding to $x=0$. Let
 \ekv{7.14} { p_0(x,\xi)=\sum_{\vert
\alpha +\beta \vert =2}{\partial _x^\alpha
\partial _\xi ^\beta p(0,0)\over \alpha !\beta !}x^\alpha \xi ^\beta
} be the quadratic approximation of $p$, so that \ekv{7.15} {
p-p_0={\cal O}((x,\xi )^3) ={\cal O}((h+(x,\xi )^2)^{3\over 2}). }
We may assume that $\Phi =\Phi _\epsilon (x)$
 is a quadratic \fu{} $\Phi
^0$ in the region $\vert x\vert \le \sqrt{Ah}$
 and for $x$ in that region,
we realize $p_0(x,hD_x)u$ with a contour as in (\ref{7.12}). The
difference between the corresponding effective kernels of $P=p^w$
and $P_0=p_0(x,hD_x)u$ is then
$$
{\cal O}(1)h^{-n}e^{-{t_0\over h}\vert x-y\vert ^2}(h+\vert x\vert ^2
+\vert
y\vert ^2)^{3\over 2}={\cal O}(1)h^{-n}e^{-{t_0\over h}\vert x-y\vert
^2}(h^{3\over 2}+\vert x\vert ^3+\vert x-y\vert ^3).
$$
We conclude that \ekv{7.16} { \Vert Pu-P_0u\Vert _{H_\Phi (\vert
x\vert \le \sqrt{Ah})}\le {\cal O}((Ah)^{3\over 2})\Vert u\Vert
_{H_\Phi }. } Here both $P$ and $P_0$ are realized with a contour
as in (\ref{7.12}). However, $p_0$ is a polynomial and we check
that if we replace $P_0u$ by the corresponding differential
expression
$$
P_0u=\sep{\sum_{\vert \alpha +\beta \vert =2}{\partial _x^\alpha
\partial _\xi ^\beta p(0,0)\over \alpha !\beta !}(x^\alpha
(hD)^\beta )^w}u(x),
$$
then we commit an error $w$, satisfying
\ekv{7.17}
{
\Vert w\Vert _{H_\Phi (\vert x\vert \le \sqrt{Ah})}\le e^{-{1\over
Ch}}\Vert u\Vert _{H_{\Phi }}.
}

Now for $P_0$ we can apply Proposition \ref{Prop7.1} and Remark
\ref{remarknormequ}. We get that for every fixed $k\in \R$ and for
$z$ in a fixed compact set avoiding the eigenvalues of
$P_0|_{h=1}$:
\ekv{rappel} { \Vert (h+d^2)^{1-k}\chi_0 u\Vert \le
C\Vert (h+d^2)^{-k}\chi_0 (P_0-hz)u\Vert +{\cal O}(h^{1\over
2})\Vert 1_Ku\Vert , } where $K$ is any fixed \neigh{} of ${\rm
supp}(\nabla \chi_0 )$.

\par Notice that we can write the last term in (\ref{rappel})
as ${\cal
O}(h^{1\over 2})\Vert (h+d^2)^{1-k}1_Ku\Vert $. We now want to
replace the fixed cutoff $\chi_0 $ in (\ref{rappel}) by $\chi_0
(x/\sqrt{Ah})$ for $A\gg 1$ \indep{} of $h$. Consider the change
of variables, $x=\sqrt{Ah}\widetilde{x}$,
$hD_x=\sqrt{Ah}\widetilde{h}D_{\widetilde{x}}$,
$\widetilde{h}=1/A$. Then
$$
p_0(x,hD_x)={h\over
\widetilde{h}}p_0(\widetilde{x},\widetilde{h}D_{\widetilde{x}} )
=:{h\over \widetilde{h}}\widetilde{P}_0,
$$
and with $d=d(x)$, $\widetilde{d}=d(\widetilde{x})$:
$$
h+d^2={h\over \widetilde{h}}(\widetilde{h}+\widetilde{d}^2),\quad
e^{-2\Phi ^0(x)/h}=e^{-2\Phi ^0(\widetilde{x})/\widetilde{h}}.
$$
Start from (\ref{rappel}) with $x,h$ replaced by
$\widetilde{x},\widetilde{h}$:
$$
\Vert (\widetilde{h}+\widetilde{d}^2)^{1-k}\chi_0
(\widetilde{x})u\Vert \le C\Vert
(\widetilde{h}+\widetilde{d}^2)^{-k}\chi_0
(\widetilde{x})(\widetilde{P}_0-\widetilde{h}z)u\Vert
+C\widetilde{h}^{1\over 2}\Vert
(\widetilde{h}+\widetilde{d}^2)^{1-k}1_Ku\Vert ,
$$
\begin{align*}
&\Vert \big({\widetilde{h}\over h}\big) ^{1-k}(h+d^2)^{1-k}\chi_0
({x\over
\sqrt{Ah}})u\Vert \le\\
&C\Vert \big({\widetilde{h}\over h}\big) ^{1-k}(h+d^2)^{-k}\chi_0
({x\over \sqrt{Ah}})(P_0-hz)u\Vert +C\widetilde{h}^{1\over 2}
\Vert \big({\widetilde{h}\over h}\big) ^{1-k}(h+d^2)^{1-k}
1_K({x\over \sqrt{Ah}})u\Vert ,
\end{align*}
\begin{equation}\label{7.43}\begin{split}
&\Vert (h+d^2)^{1-k}\chi_0 ({x\over
\sqrt{Ah}})u\Vert \le \\ &C\Vert (h+d^2)^{-k}\chi_0 ({x\over
\sqrt{Ah}})(P_0-hz)u\Vert +{C\over \sqrt{A}} \Vert (h+d^2)^{1-k}
1_K({x\over \sqrt{Ah}})u\Vert .
\end{split}
\end{equation}
This estimate will be applied
with $k=1/2$.

\par We now return to the full \op{} $P$ (on the FBI-side)
 and the norms
and scalar products will now be with respect to $e^{-2\Phi /h}$,
$\Phi =\Phi _\epsilon $, $\epsilon >0$ small and fixed. Recall
however that $\Phi =\Phi ^0$ in $\vert x\vert \le \sqrt{Ah}$. Let
$\chi$ be a cutoff function equal to $1$ in a fixed neighboorhood
 of the critical points, but recall however the simplifying
assumption that we only have one critical point corresponding to
$x=0$. Let us denote
\begin{equation} \label{defLambdaa}
\Lambda ^2=h+\min (d^2,(dAh)^{2\over 3}).
\end{equation}
 Using (\ref{6.2ter}) as in
Section \ref{locallarge}, we get for $z={\cal
  O}(1)$,
 \ekv{7.44} {
\norm{\Lambda u}^2\le C(\Re (\chi(x)(P-hz)u\vert u)+C^2(\chi_0
^2({x\over \sqrt{Ah}})\Lambda u\vert\Lambda  u)) +
C'\norm{(1-\chi)\Lambda u}\norm{\Lambda u}. }  Then using
(\ref{7.43}) we get for $\tau>0$:
\begin{equation}\label{7.45}\begin{split}
\norm{\Lambda u}^2 & \le C\Vert \Lambda ^{-1}(P-hz)u\Vert \Vert
\Lambda u\Vert +C\Vert \Lambda \chi_0 ({x\over \sqrt{Ah}})u\Vert ^2
+ C'\norm{(1-\chi)\Lambda u}\norm{\Lambda u} \\
  &\le   {C\over \tau }\Vert \Lambda ^{-1}(P-hz)u\Vert
^2+C\tau \Vert \Lambda u\Vert ^2 +\widetilde{C}\Vert \Lambda
^{-1}\chi_0 ({x\over \sqrt{Ah}})(P_0-hz)u\Vert ^2
 \\
  &  +{\widetilde{C}\over A}
  \Vert \Lambda 1_K({x\over \sqrt{Ah}})u\Vert
^2 + C'\norm{(1-\chi)\Lambda u}\norm{\Lambda u} .
\end{split}\end{equation}
Here we also need (and we can clearly generalize (\ref{7.16}) for that
purpose)
\ekv{7.46}
{
\Vert \chi_0 ({x\over \sqrt{Ah}})\Lambda ^{-1}(P-P_0)u\Vert
\le C(A)h^{1\over
2}\Vert \Lambda u\Vert .
}
Insertion in (\ref{7.45}) gives
\begin{equation}\label{7.47}\begin{split}\Vert \Lambda u\Vert ^2\le
&{C\over \tau }\Vert \Lambda ^{-1}(P-hz)u\Vert ^2+C\tau \Vert
\Lambda u\Vert
^2 \\
&+2\widetilde{C}\Vert \Lambda ^{-1}\chi_0 ({x\over
\sqrt{Ah}})(P-hz)u\Vert ^2+\widetilde{C}(A)h\Vert \Lambda u\Vert
^2+{\widetilde{C}\over A}\Vert \Lambda u\Vert ^2 \\
& + C'\norm{(1-\chi)\Lambda u}\norm{\Lambda u}.\nonumber
\end{split}\end{equation}
Choosing first $\tau $, $1/A$ small enough and then $h$ small
enough, we get \ekv{7.488} { \Vert \Lambda u\Vert \le C\Vert
\Lambda ^{-1}(P-hz)u\Vert + C''\norm{(1-\chi)\Lambda u}. } and
noticing that $h\le\Lambda^2 \leq Ch^{2/3}$, we get the main
result of this section \ekv{7.48} { h\Vert  u\Vert \le C\Vert
(P-hz)u\Vert + C''h^{5/6}\norm{(1-\chi) u}. }

\Section{Review of semiclassical Weyl Calculus} \label{Secweyl}
 In this section
we introduce some tools and make some remarks about the
translation into the semiclassical point of view of some basic
facts on the classical Weyl-H{\"o}rmander Calculus.

\subsection{Weyl-H{\"o}rmander calculus}

First recall the framework of the Weyl-H{\"o}rmander Calculus, which
 can be found in   \cite[Chapter 18]{Hor85}. We put a subscript
$cl$ everywhere here to emphasize the fact the we are in the
original (opposite to  semiclassical) framework of the calculus.
 Recall that the classical Weyl quantization is given for
 an admissible
 symbol $p_{cl}$ (to be defined below) by
 \begin{equation} \label{weyl}
(p_{cl}^{w_{cl}} u)(x) = \frac{1}{(2 \pi)^n}\int \!\!\! \int
         e^{i\seq{x-y,\xi}} p_{cl}(\frac{x+y}{2},\xi)u(y) dy d\xi.
\end{equation}
Consider the symplectic space
 $\R^{2n}$ equipped with the symplectic form
 $\sigma = \sum_{i=1}^n d\xi_i \wedge d x_i$.
 If  $g$ is a  positive definite quadratic form, we define
\begin{equation} \label{defgs}
 g_{cl}^\sigma(T) = \sup_{g_{cl}(Y)=1} \sigma \sep{T,Y}^2,
\end{equation}
which is also a positive definite quadratic form. We say that
$g_{cl}$ is a $cl$-admissible metric if
\begin{equation}  \label{admissible}
\begin{array}{ll}
  \forall X \in \R^{2n}, \ \ g_{cl, X} \leq
g_{cl, X}^\sigma & \textrm{($cl$-Uncertainty Principle)}, \\
\exists C_0>0 \textrm{  such that } g_{cl, X}(X-Y) \leq C_0^{-1}
\Longrightarrow \sep{ g_{cl, X}/g_{cl, Y} }^{\pm 1} \leq C_0
& \textrm{($cl$-slowness)}, \\
  \exists C_1, N_1>0  \textrm{  such that }
  g_{cl, X}/g_{cl, Y}\leq  C_1 \sep{ 1 + g_{cl, X}^\sigma(X-Y)
}^{N_1} & \textrm{($cl$-temperance)},
\end{array}
\end{equation}
for positive constants $C_0, C_1, N_1$. Let us note that if the
metric $g_{cl}$ depends on a parameter (for example $h$), we call
it $cl$-admissible if (\ref{admissible}) occurs uniformly in this
parameter. The same is true for the $cl$-admissible weights we
introduce now. A $cl$-admissible weight is a positive function
$m_{cl}$ on the phase space $\R^{2n}$, for which there exists
 $\tilde{C}_0,\tilde{C}_1,\tilde{N}_1>0$ such that
 \begin{equation} \label{admissibleweight}
\begin{array}{ll}
 g_{cl, X}(X-Y) \leq \tilde{C}_0 \Longrightarrow \sep{
m_{cl}(Y)/m_{cl}(X) }^{\pm 1}
 \leq \tilde{C}_0   & \textrm{($cl$-slowness),}  \\
 m_{cl}(Y)/m_{cl}(X)  \leq \tilde{C}_1 \sep{ 1 +
g_{cl, X}^\sigma(X-Y) }^{\tilde{N_1}} &
\textrm{($cl$-temperance)}.
\end{array}
\end{equation}
We define next the {\it $cl$-uncertainty parameter} $\lambda_{cl}$,
 which is a special admissible weight for $g$,
\begin{equation}  \label{unct}
\lambda_{cl}(X) = \inf_{T \in \R^{2n}/\set{0}} \sep{ g_{cl,
X}^\sigma(T) / g_{cl, X}(T) }^{1/2} \geq 1.
\end{equation}
 Let us  now introduce some  spaces of symbols.
 We say that a function $p_{cl}$ is
a symbol in $S(m_{cl},g_{cl})$ if  $p_{cl} \in
\cc^\infty(\R^{2n})$,
 and if the following semi-norms are finite
 \begin{equation}  \label{eqn:smg}
\sup_{
    X \in \R^{2n},
   g_{cl,X} (T_j)\leq 1  }
 \abs{\seq{ p_{cl}^{(k)}(X),T_1\otimes ...\otimes T_l }}
  m_{cl}^{-1}(X) .
\end{equation}
If $m_{cl}$ is of the form $\lambda_{cl}^\mu$, we say that
$p_{cl}$ is of order $\mu$. For good symbols (in
$S(m_{cl},g_{cl})$ classes for instance),
 we define the composition law $\sharp_{cl} $ such that
 $(p_{cl} \sharp_{cl} q_{cl})^{w_{cl}} =
  p_{cl}^{w_{cl}}\circ q_{cl}^{w_{cl}}$ by
\begin{equation}   \label{dieze}
(p_{cl} \sharp_{cl} q_{cl})(x,\xi) = e^{\frac{i}{2} \sigma
\sep{(D_x,D_\xi), (D_y, D_\eta)}}p_{cl}(x,
\xi)q_{cl}(y,\eta)\vert_{y=x,\eta=\xi} ,
\end{equation}
and for $p_{cl} \in S(m_1, g_{cl})$, $q_{cl} \in S(m_2,g_{cl})$,
if $\set{.,.}$ denotes the Poisson bracket, then   there is
 $r_{cl} \in S(m_1m_2 \lambda_{cl}^{-2},g_{cl})$ such that
\begin{equation} \label{diezedvptcl}
p_{cl} \sharp_{{cl}} q_{cl} = p_{cl}q_{cl} + \frac{1}{2i}
\{p_{cl},q_{cl}\} +r_{cl}.
\end{equation}
Recall eventually the Fefferman-Phong inequality that will be used
in the next sections:
\begin{prop}
Let $p_{cl} \in S(m_{cl}, g_{cl})$.  If $p_{cl} \geq 0$ then there
is a real symbol $r_{cl} \in S(m_{cl}\lambda_{cl}^{-2}, g_{cl})$
such that $p_{cl}^w \geq r^w$. Hence if $m_{cl} =
\lambda_{cl}^2$, then $p_{cl}^w$ is bounded from below.
\end{prop}

\subsection{Semiclassical Weyl-H{\"o}rmander Calculus}

The original calculus already contains a parameter that
plays the role of a Planck's constant, namely the inverse of the
uncertainty parameter. In the semiclassical case this is made more
explicit, but basically this is only a reduction to the original
calculus by a change of variables.

For  an admissible symbol $p$ we first  recall the
definition of semiclassical Weyl quantization
$$
p^wu = \frac{1}{(2\pi h)^n} \iint p\sep{ \frac{x+y}{2}, \xi} e^{
i\seq{x-y, \xi}/h } u(y) dy d\xi, \ \ \ u\in \ss .
$$
A straightforward computation shows that
\begin{equation}
p^w = p_{cl}^{w_{cl}} \ \ \ \  \text{ where $p_{cl}(x, \xi)
 = p(x,h\xi)$}.
\end{equation}
Now observe that  $p$ belongs to a symbol class $S(m, g)$ for
 a Riemanian metric $g$ and  a positive function $m$ if and only if
$p_{cl} \in S(m_{cl}, g_{cl})$ where
$$
m_{cl}(x,\xi) = m(x, h\xi), \ \ \ g_{cl, (x,\xi)}(t, \tau) =
g_{(x, h\xi)}(t, h\tau).
$$
Using definition (\ref{defgs}) and (\ref{unct}) for defining
respectively $g_{cl}$, $g$, and  $\lambda_{cl}$, $\lambda$, we
also get
\begin{equation}
g_{cl, (x,\xi)}^\sigma (t, \tau) = h^{-2}g_{(x, h\xi)}^\sigma (t,
h\tau), \ \ \ \textrm{ and } \ \ \lambda_{cl}(x, \xi) = h^{-1}
\lambda(x, h\xi).
\end{equation}
As a consequence it is natural to introduce the following
definitions in the semiclassical case:
\begin{dref}
We say that $g$ is an admissible (or semiclassically admissible)
metric if
\begin{equation}  \label{admissible2}
\begin{array}{ll}
  \forall X \in \R^{2n}, \ \ g_{ X} \leq
h^{-2} g_{ X}^\sigma \ \ \ (\textrm{i.e. } \lambda \geq h)
& \textrm{(Uncertainty Principle), } \\
\exists C_0>0 \textrm{  such that } g_{X}(X-Y) \leq C_0^{-1}
\Longrightarrow \sep{ g_{X}/g_{Y} }^{\pm 1} \leq C_0
& \textrm{(slowness), } \\
  \exists C_1, N_1>0  \textrm{  such that }
  g_{X}/g_{Y}\leq  C_1 \sep{ 1 + h^{-2}g_{X}^\sigma(X-Y)
}^{N_1} & \textrm{(temperance),  }
\end{array}
\end{equation}
for positive constants $C_0, C_1, N_1$.
\end{dref}
A direct definition holds for semiclassical weights. Using this
(note that all this is simply a change of variables) we can write
\begin{lemma}
The metric $g$ is an admissible metric of uncertainty parameter
$\lambda (\geq  h) $ if and only if  $g_{cl}$ is an admissible
metric of uncertainty parameter $\lambda_{cl} (\geq 1)$, both
uniformly in $0<h\leq 1$.
\end{lemma}
We can therefore translate into the semiclassical point of view
all the classical results. First observe that
 symbols of order $1$ give
bounded operators on $L^2(\R^n)$. Then the product formula is
defined by
$$
p^w \circ q^w = (p\sharp q)^w ,
$$
where
$$
p\sharp q (x, \xi, h) = e^{\frac{ih}{2} \sigma \sep{(D_x,D_\xi),
(D_y, D_\eta)}}p(x, \xi,h)q(y,\eta,h)\vert_{y=x,\eta=\xi}.
$$
The asymptotic expansion is then given  for $p \in S(m_1, g)$, $q
\in S(m_2,g)$  by
\begin{equation} \label{diezedvpt}
p \sharp q = pq + \frac{h}{2i} \{p,q\} + h^2 r,
\end{equation}
where $r \in S(m_1m_2 \lambda^{-2},g)$.  Recall eventually how to
write the semiclassical Fefferman-Phong inequality that will be
used in the text:
\begin{prop}
If $p \geq 0$ then there is a real symbol $r \in S(m\lambda^{-2},
g)$ such that $p^w \geq h^2r^w$. Hence if $m =
h^{-2}\lambda^2$, then $p^w$ is bounded from below uniformly
\wrt{} $h$.
\end{prop}

\remark As an illustration, let us see what happens in the case of
the constant metric $g = dx^2 + d\xi^2$. It is the one generally
used in semiclassical work. We check immediately that it is
admissible in the sense of definition \ref{admissible2}, since
$g_X = g_Y$ for all $X$, $Y$ and that $g^\sigma /g = 1 \geq h$. Of
course the translation procedure  gives  the Fefferman--Phong
inequality:  $ p^w \geq -Ch^2$ if $p$ is real non negative with
all its derivatives  bounded.

 \subsection{The microlocal metric $\Gamma$}
 We study now  a particular metric used in the next sections.
 \begin{lemma} \label{metricgamma} the metric defined
 on $\R^{2n}$ by
 $$
\Gamma= \frac{dx^2}{h^{2/3}} + \frac{d \xi^2}{\mu^2}, \ \ \ \
\textrm{ where }   \mu^2 = p_1 + (h\lambda)^{2/3},
$$
is  (semiclassically) admissible.
\end{lemma}

\preuve
 Recall that we suppose that
$$
  \Gamma_0 = dx^2 + d\xi^2/ \lambda^2,\
  \lambda =\lambda (x,\xi ) \ge 1,
$$
is a $cl$-admissible metric. Let us prove the three points of
(\ref{admissible2}). We first notice that
$$
  \Gamma^\sigma = \mu^2 dx^2 + h^{2/3} d \xi^2 ,
$$
therefore the uncertainty parameter of $\Gamma$ is $\mu h^{1/3}$
and we have for $h$ small
$$
  \mu h^{1/3} \geq \lambda^{1/3} h^{2/3} \geq   h^{2/3} \geq h ,
$$
therefore $\Gamma $ satisfies the uncertainty principle.

\ititem{Slowness of $\ \Gamma$.}
We take $X=(x, \xi)$ and $Y= (y, \eta)$ and we
observe that if $\Gamma_X(X-Y) \leq C_0$ then
\begin{equation} \label{lenteurgm0}
  |x-y|^2 \leq C_0 h^{2/3}, \ \ \
  \textrm{ and } \ \ \
  |\xi-\eta|^2 \leq C_0\sep{ p_1(X) +(h\lambda)^{2/3}(X)}.
\end{equation}
Using a Taylor expansion and the fact that the the second
derivative of
$p_1$ is bounded,  we can write that
\begin{equation*}
\begin{split}
  p_1(Y) & \leq p_1(X) + |\nabla p_1| |X-Y| + C |X-Y|^2 \\
         & \leq p_1(X) + C'\sqrt{p_1} |X-Y| + C |X-Y|^2 \\
         & \leq 2 p_1(X) + C'' |X-Y|^2,
\end{split}
\end{equation*}
where for the second inequality we used inequality (\ref{deriv})
for the non negative function $p_1$. Now use the fact that
$\Gamma_X(X-Y) \leq C_0$. We get
\begin{equation} \label{slow1}
 \begin{split}
 p_1(Y)        & \leq  2p_1(X) + C'' |X-Y|^2 \\
              &  \leq 2p_1(X) + C''C_0\sep{h^{2/3} +  p_1(X) +
              (h\lambda)^{2/3}(X)} \\
              & \leq C(p_1(X) +(h\lambda)^{2/3}(X)) ,
 \end{split}
\end{equation}
since $\lambda\geq 1$. Formula (\ref{lenteurgm0})  implies that
 \begin{equation*} \label{lenteurgm}
 |x-y|^2 \leq C_0 , \ \ \
 \textrm{ and } \ \ \ |\xi-\eta|^2 \leq C_0 \lambda^2(X),
 \end{equation*}
and we get using the slowness of $\Gamma_0$ for $C_0$ sufficiently
small that $(h\lambda)^{2/3}(Y) \leq C'(h\lambda)^{2/3}(X)$.
Using this and (\ref{slow1}) yields
$$
  p_1(Y) + (h\lambda)^{2/3}(Y) \leq C(p_1(X) + (h\lambda)^{2/3}(X)),
$$
that is to say $\mu(X) \leq C\mu(Y)$. This implies immediately that
$\Gamma_Y \leq C \Gamma_X$. Inverting the roles of $X$ and $Y$ proves
the slowness of $\Gamma$.

\ititem{Temperance of $\ \Gamma$.} Again we denote  $X=(x, \xi)$
and $Y= (y, \eta)$. Beginning from the first line of (\ref{slow1})
we write
\begin{equation} \label{temp1}
 \begin{split}
 p_1(Y)        & \leq  2p_1(X) + C'' |X-Y|^2 \\
              &  \leq C \sep{ p_1(X)  +
              (h\lambda)^{2/3}(X)}(1+ h^{-2/3}|X-Y|^2).
 \end{split}
\end{equation}
Notice that
\begin{equation*}
\begin{split}
h^{-2} \Gamma^\sigma_X(X-Y) & = h^{-2} \sep{ \sep{p_1(X)
+(h\lambda)^{2/3}(X) }|x-y|^2 + h^{2/3}|\xi-\eta|^2}\\
&  \geq h^{-4/3}|X-Y|^2 \\
&  \geq h^{-2/3}|X-Y|^2,
\end{split}
\end{equation*}
since $\lambda \geq 1$ and for $h\leq 1$. Hence
\begin{equation} \label{tempp}
p_1(Y)        \leq C \sep{ p_1(X)  +
              (h\lambda)^{2/3}(X)}\sep{1+ h^{-2}
              \Gamma^\sigma_X(X-Y)} .
\end{equation}
Since $\Gamma_0 = dx^2 + d\xi^2/\lambda^2$ is
$cl$-temperate, there exists  $C_0, N \geq 1$ such that
$$
\Gamma_{0, X} \leq C_0 \Gamma_{0, Y} \sep{1 + \Gamma_{0,X}^\sigma
(X-Y)}^N.
$$
Together with the fact that $\Gamma_{0}^\sigma = \lambda^2 dx^2 +
d\xi^2$, this implies that
\begin{equation} \label{temp2}
 \begin{split}
\lambda^2(Y)
     & \leq C_0 \lambda^2(X) \sep{ 1 + \lambda^2(X) |x-y|^2
       + |\xi - \eta|^2}^N \\
       & \leq C_0' \lambda^2(X) \sep{ 1 + \lambda^{2/3}(X) |x-y|^2
       + |\xi - \eta|^2}^{3N} \\
       &  \leq C_0' \lambda^2(X) \sep{ 1 + h^{-2}
        \sep{ ((h\lambda)^{2/3}(X)+p_1(X)) |x-y|^2
       + h^{2/3}|\xi - \eta|^2}}^{3N},
\end{split}
\end{equation}
since $h^{-4/3} \geq 1$ and $p_1 \geq 0$. Now we recognize in the
parentheses a  term of the form $h^{-2} \Gamma^\sigma$.
Multiplying by $h$ and raising to  the power $1/3$ gives
$$
(h\lambda)^{2/3}(Y) \leq  C (h\lambda)^{2/3}(X) \sep{ 1 +
        h^{-2} \Gamma^\sigma_X(X-Y)}^N.
$$
Together with (\ref{temp1}) this gives
$$
\mu^2(Y) \leq C \mu^2(X)\sep{ 1 +
        h^{-2} \Gamma^\sigma_X(X-Y)}^N,
$$
which implies $\Gamma_X \leq \Gamma_Y \sep{ 1 + h^{-2}
\Gamma^\sigma_X(X-Y)}^N$. Consequently $\Gamma$ is
(semiclassically) temperate. Eventually we have proven that
$\Gamma$ is a (semiclassically) admissible metric. \fin

\Section{Resolvent estimates away from the critical points when
$|z| \gg h$} \label{awaylarge}

In this section we  suppose that $p$ satisfies hypotheses \bf
(H2), (H3), (H4) \rm  and we shall work away from a fixed neighborhood
${\cal B}$
of the critical
points and for  $|z| \gg h$. The main result of this section will be
the estimate (\ref{main8}).  At infinity in the phase space, we
shall use  the machinery of the Weyl calculus. Let us consider the
following weight
$$
  \mu^2(x,\xi) = p_1(x,\xi) + (h\lambda(x,\xi))^{2/3}
$$
We notice that $\mu \geq h^{1/3}$ . We  use the metric defined in
lemma \ref{metricgamma}
\begin{equation} \label{defGamma}
  \Gamma =\frac{dx^2}{h^{2/3}} + \frac{d \xi^2}{\mu^2} .
\end{equation}
From the construction of the weight $G$ in  Proposition
\ref{propescape} (cf (\ref{derivG}), (\ref{defG})),  we know that
$$
  g \defegal G/h \in S\sep{ 1,\Gamma} \text{ outside } \bb ,
$$
since $G =0$ when $p_1 \geq  2(h\lambda)^{2/3}/M$. There
 is no restriction to
extend $g$ near the critical points and let it uniformly
 be in the class
$S\sep{ 1,\Gamma}$.

From Proposition \ref{propescape},  we have the following two
estimates for our new $g$:
\begin{equation} \label{derivggg}
 g \in S\sep{ 1,\Gamma}, \  \D g \in S(\mu^{-1}, \Gamma).
\end{equation}
We verify now that some other symbols are good symbols for the
metric $\Gamma$. We first observe the evident fact that
$S(m,\Gamma _0) \subset S(m', \Gamma)$ for all weights $m' \geq
m$, since $\Gamma _0\le \Gamma $. From (\ref{symbcal}) we get
\begin{equation} \label{derivppp}
\D p \in S(\lambda,dx^2 + d\xi^2/\lambda^2) \Rightarrow \D p \in
S(\mu^3h^{-1}, \Gamma),
\end{equation}
since $\mu^3 \geq h\lambda$. Of course in this new class,
$p$, $\D p$ are
no more  symbols of order $2$ and $1$
respectively. Nevertheless the real part $p_1$ has a good behavior:
\begin{equation} \label{derivpp1}
p_1 \in S(\mu^2, \Gamma).
\end{equation}
 Indeed, $0 \leq p_1 \leq \mu^2$ and
since the second derivative of $p_1$ is bounded we use
(\ref{deriv}) to get $|\D p_1| \leq C \sqrt{p_1} \leq C\mu$.
Moreover,  $\D^2 p_1 \in S(1,dx^2 + d\xi^2/\lambda^2)$ gives
$\D^2 p_1 \in S(1, \Gamma)$. This implies
(\ref{derivpp1}).

From the preceding section we know that $\Gamma$ is a
(semiclassically) admissible metric of uncertainty parameter
$h^{1/3} \mu$. We have therefore the following symbolic expansion
for the composition of $q_1 \in S(m_1, \Gamma)$ and $q_2 \in
S(m_2, \Gamma)$:
\begin{equation} \label{qsharpq}
q_1\sharp q_2 (x, \xi, h) = q_1 q_2 (x, \xi, h)
+ {h \over 2i} \set{q_1, q_2}
(x, \xi, h) + h^2R_2(q_1,q_2)(x, \xi, h),
\end{equation}
where
\begin{equation} \label{remain1}
R_2(q_1,q_2) \in S( m_1m_2 (h^{1/3} \mu)^{-2}).
\end{equation}
This means that in the remainder of order two in the asymptotic
expansion of the sharp product, we have a gain of $ (h^{1/3}
\mu)^{-1}$ to the square in addition to the gain of $h^2$ due to
the semiclassical point of view. The  Fefferman-Phong inequality
reads for $\Gamma$:
\begin{lemma} \label{fphong} Let $m$ be an $h$-admissible
weight and $q \in S(m, \Gamma)$. If
$\Re q \geq 0$ then there is a real symbol
$r \in S(m h^2 (h^{1/3} \mu)^{-2})$
such that $\Re(q^w u, u) \geq (r^w u, u)$ for all $u\in \ss$.
 In particular
symbols in $S(h^{-2}(h^{1/3} \mu)^{2}, \Gamma)$ with non-negative
real part   correspond to operators with real part bounded from
below by an $h$-\indep{} constant in the \op{} sense.
\end{lemma}

For the symbols we deal with, we noted in
 (\ref{derivggg}-\ref{derivppp}) that
 $\D p$ and $\D g$ have better symbolic estimates than the one
  given by the symbolic classes of  $p$ and
 $g$. This gives improvements to the symbolic calculus. Let us write
 explicitly the expansion of  $q_1 \sharp q_2$ to the order $d$
\begin{equation}   \label{sjodiezedvpt}
\begin{split}
  (q_1 \sharp q_2)(x, \xi, h) =
  &  \sum_{j=0}^{d-1} \frac{h^j}{j!}
  \sep{
  \frac{i}{2} \sigma \sep{D_{x,\xi }, D_{y,\eta }}}^j  q_1(x,
  \xi,h)q_2(y,\eta,h) \vert_{y=x,\eta=\xi} \\
   &  + h^dR_d(q_1, q_2) (x, \xi, h) ,
\end{split}
\end{equation}
where
\begin{equation}
\begin{split}
 R_d (q_1, q_2) (x, \xi, h) = &
 \int_0^1  \frac{ (1-\theta)^{d-1}}{(d-1)!}
e^{\frac{i\theta
h}{2} \sigma \sep{D_{x,\xi }, D_{y,\eta }}} \\
  &  \sep{{i\over 2}\sigma (D_{x,\xi },D_{y,\eta })}^d
   q_1(x, \xi,h)q_2(y,\eta,h)
\vert_{y=x,\eta=\xi}  d\theta .
\end{split}
\end{equation}
The order (as a symbol in a class $S(m, \Gamma)$),  computed as in
the classical
 case, is exactly  the order of the symbol appearing on
 the second line
 $$
\sep{{i\over 2}\sigma (D_{x,\xi },D_{y,\eta })}^d  q_1(x,
\xi,h)q_2(y,\eta,h) \vert_{y=x,\eta=\xi}.
$$
Now return to the case of $p$ and $g$ with $d=2$. A
straightforward computation using (\ref{derivggg}-\ref{derivppp})
gives that
\begin{equation*}
\begin{split}
\sep{{i\over 2}\sigma (D_{x,\xi },D_{y,\eta })}^2
  & g(x,  \xi,h)p(y,\eta,h) \vert_{y=x,\eta=\xi} \\
  &  \in S(\mu^3h^{-1} \times \mu^{-1} \times
h^{-1/3} \mu^{-1}, \Gamma) \subset S(h^{-4/3} \mu, \Gamma),
\end{split}
\end{equation*}
hence
$$
R_2(g,p) \in S(h^{-4/3} \mu, \Gamma),
$$
 so
\begin{equation} \label{symbpg}
g\sharp p = gp + \frac{h}{2i} \set{g,p} + r \ \ \ \text{ with } \
\ \ r = h^2 R_2(g,p) \in S(h^{2/3} \mu, \Gamma) .
\end{equation}
(Note that this implies $r \in S(h^{1/3}\mu^2, \Gamma)\subset
S(\mu^2, \Gamma)$ since  $h^{1/3} \leq
\mu$).

Let us now fix  $\eps >0$ and take $z\in \C$. We can write for
$u\in \ss$ using (\ref{symbpg}) that
\begin{equation} \label{pdiezeg}
\begin{split}
\Re \sep{ (p^w -z)u, (1 -\eps g )^w u}
& = \Re \sep{ \sep{(1-\epsilon g)\sharp (p -z)}^w u,u} \\
& =    \sep{  ((p_1 -\Re z)(1 -\eps g ) + \eps h \set{p_2, g}/2
  - \eps \Re r)^w u,u},
\end{split}
\end{equation}
where $r \in S(\mu^2, \Gamma)$ was defined in (\ref{symbpg}).
Let us study
the first two terms in the asymptotic development of
 $\Re (p -z)\sharp (1 -\eps g )$. For $\eps$
sufficiently small, we have from (\ref{eqp})
\begin{equation}
 p_1 + \eps h \set{p_2, g}/2  \geq \eps_1
 \sep{(h\lambda)^{2/3} +p_1} =
 \eps_1 \mu^2 ,
\end{equation}
when $|(x, \xi)| \geq \oo(1)$ far from the critical points
 (recall that $G
\defegal hg$ in (\ref{eqp})). This means that
$ p_1 + \eps_0 h \set{p_2, g}/2$ is elliptic in $S(\mu^2, \Gamma)$
far from the
critical points. Choose $\phi \in \cc_0^\infty$ equal to
$1$ in a neighborhood of the critical points, so that
\begin{equation}
 p_1 + \eps h \set{p_2, g}/2  \geq \eps_1 \mu^2 - C\mu ^2
 \phi(x, \xi).
\end{equation}
Recall that $r\in S(h^{1/3}\mu ^2,\Gamma )$.
 Using this and choosing $\eps$
sufficiently small yields
\begin{equation} \label{major0}
\begin{split}
 \Re (p -z)\sharp (1 -\eps g )
& = (p_1 -\Re z)(1 -\eps g ) + \eps h \set{p_2, g}/2 - \eps \Re r \\
& \geq c \mu^2 - 2\max (\Re (z), 0)  - \mu^2 \phi .
\end{split}
\end{equation}
Let us now introduce
$$
Z \defegal h^{2/3} |z|^{1/3}.
$$
We follow the preceding computations, and get with $\eps_2 > 0$ that
\begin{equation} \label{estim2}
\begin{split}
&  c \mu^2 - 2\max (\Re (z), 0)  - C\mu ^2\phi \\
& \hspace*{1cm}
  \geq (c/2)\mu^2 + \frac{c}{2} (\mu^2 - \eps_2Z) +
  \sep{\frac{c\eps_2}{2}Z
-2\max (\Re (z), 0)}  - C\mu ^2\phi.
\end{split}
\end{equation}
We will bound from below each term of the right hand side. We
assume that
\begin{equation} \label{Sigma}
\frac{c\eps_2}{2}Z  \geq 4\Re (z).
\end{equation}
It defines a region $\Sigma$ in the complex plane,
and if $z$ is in this region  the third term of (\ref{estim2}) is
bounded from below by $cZ$. To study the second term we observe
that $\mu^2
\geq \eps_2Z$ since $\lambda^2 \geq |z|\eps_2^3$. Now choose
a cutoff function $\psi_1(t)$ supported in the ball of radius
 $2\eps_2^3$ and equal to one
in the ball of radius $\eps_2^3$. Then
$$
\frac{c}{2} (\mu^2 - \eps_2 Z) \geq  -c'' Z \psi_1^2(\lambda^2/|z|)
$$
Summing up the preceding results we have obtained the following
bound, where $c$, $C$ denote fixed constants
\begin{equation} \label{major}
(p_1 -\Re z)(1 -\eps g ) + \eps h \set{p_2, g}/2
+ \eps \Re r \geq c(\mu^2 + Z)
-CZ\psi_1^2(\lambda^2/|z|) -C\mu^{2}\phi
\end{equation}
Note that $\psi _1^2(\lambda ^2/|z|)\in S(1,\Gamma _0)$.
Now we want to go back to   the operator side. We first notice
that dividing
the two sides of (\ref{major}) by $Z$ yields an inequality in
 $S(h^{-1} \mu^2,
\Gamma)$ uniformly in $z$, which we recall can be arbitrarily
large. Indeed the
terms $p_1$, $h\set{p_2, g}$, $r$ and $\mu^2$ are in
$S( \mu^2, \Gamma)$ and
since $Z \gg h$ (from $|z| \gg h$) we get that these operators
divided by $Z$
are in $S(h^{-1} \mu^2, \Gamma)$. The others (divided by $Z$) are
 bounded by a
constant since by hypothesis $\max\set{\Re(z), 0} \leq CZ$, and a
fortriori are
in $S(h^{-1} \mu^2, \Gamma)$.

Let us apply the inequality of Fefferman-Phong, Lemma \ref{fphong},
in this
class to this operator. We get using (\ref{pdiezeg}-\ref{major})
divided
 by $Z$ and then multiplying by $Z$
\begin{equation}
\begin{split}
((\mu^2 + Z)^wu,u) \leq{}& C\Re((p^w-z)u, (1-\eps g)^w u) + CZ \sep{
  \psi_1^2(\lambda^2/|z|)^w u,u} \\
& + C((\mu ^2\phi)^w u,u) + Z h^2\Re(R^w u, u),
\end{split}
\end{equation}
where $h^2 R$ is of order $h^2 (h^{-1} \mu^2) (\mu h^{1/3})^{-2} =
h^{1/3}$ (recall that $\mu h^{1/3}$ is the uncertainty parameter of
$\Gamma$). Choosing $h$ small enough and using (\ref{chainemanquante})
below, gives
$$Z h^2 R^w \leq {1\over 4}Z\le {1\over 2}(Z+\mu ^2)^w,$$
and we therefore get for $h$ small enough and an other constant
$C$
\begin{equation} \label{bonestim}
((\mu^2 + Z)^wu,u) \leq C\Re((p^w-z)u, (1-\eps g)^w u) + CZ \sep{
\psi_1^2(\lambda^2/|z|)^w u,u } + C((\mu ^2\phi)^w u,u).
\end{equation}

We shall use the following
\begin{lemma} \label{psilz} we have
$ \sep{ \psi_1^2(\lambda^2/|z|)^w u,u } \leq  \frac{ C}{\max
(1,|z|^2)} \norm{(p^w-z) u}^2 + C h \norm{u}^2. $
\end{lemma}

Let us suppose for a while that this lemma is proven. We first
write that for $\eps$ sufficiently small,
$$
\Re((p^w-z)u, (1-\eps g)^w u) \leq C\norm{(p^w-z)u}\norm{u}.
$$
Then we observe that $\mu^2 \geq 0$ and the Fefferman-Phong
inequality in $S(\mu^2, \Gamma)$ yields $(\mu^2)^w \geq
-Ch^{4/3}$. Since $Z \gg h$, we have for $h$ sufficiently small
\ekv{chainemanquante}{
Z \norm{u}^2 \leq 2((Z+ \mu^2)^w u, u).}
Then we use this result and the lemma which yields from
 (\ref{bonestim}) that
$$
Z \norm{u}^2 \leq C\norm{(p^w-z)u}\norm{u} +  Z \frac{ C}{\max(1,
|z|^2)} \norm{(p^w-z) u}^2 + CZh\norm{u}^2
+ C\norm{(\mu ^2\phi )^wu}\norm{u} .
$$
 Choosing $h$ sufficiently small and noticing
 that $Z^2 \leq \max(1,|z|^2)$
yields the main result of this section
\begin{equation}\label{main8}
Z \norm{u} \leq C\norm{(p^w-z)u} + C\norm{(\mu ^2\phi )^wu}
\end{equation}
where we  recall that $|z| \gg h$ and that $\Re(z)  \le C Z
 \defegal C h^{2/3}
|z|^{1/3}$.\\

It remains to prove Lemma \ref{psilz}.

 \preuve[of Lemma \ref{psilz}]
 We first  observe  that for $|z| \ll \oo(1)$, we have
 $\psi_1^2(\lambda^2/|z|) =0$ since $\lambda \geq 1$
 and the support of $\psi_1$ is bounded. Therefore we can
 suppose that $|z|\geq \oo(1)$, since in the other case, the left
 member of the inequality in the lemma is zero.
To prove the result we can go back to the original metric $dx^2 +
d\xi^2/\lambda^2$. We first notice that since $p  =
\oo(\lambda^2)$ we can choose the support of $\psi_1$ (i.e.
$\eps_1$ in (\ref{estim2})) such that
$$
|p-z| \geq |z|/2 \text{ on  the support of } \psi_1.
$$
We notice also that uniformly with respect to $z$, we have
$$
\frac{\sep{p-z}}{|z|} \psi_1(\lambda^2/|z|) \in S(1,dx^2 +
d\xi^2/\lambda^2).
$$
 We therefore
have the following inequality in $S(1,dx^2 + d\xi^2/\lambda^2)$:
\begin{equation}
\begin{split}
\psi_1^2(\lambda^2/|z|) & \leq 4\frac{|{p-z}|^2}{|z|^2}
\psi_1^2(\lambda^2/|z|) \\
& \leq 4\Re \frac{(\overline{p-z})}{|z|} \psi_1(\lambda^2/|z|)
\sharp \frac{\sep{p-z}}{|z|} \psi_1(\lambda^2/|z|) + hR ,
\end{split}
\end{equation}
where by the symbolic calculus, $hR \in S(h\lambda^{-1},dx^2 +
d\xi^2/\lambda^2) \subset S(h, dx^2 + d\xi^2/\lambda^2)$.
Using the G\aa rding inequality for this inequality, we get
\begin{equation} \label{truc1}
\sep{\psi_1^2(\lambda^2/|z|)^w u, u} \leq
\norm{\sep{\frac{\sep{p-z}}{|z|} \psi_1(\lambda^2/|z|)}^w u}^2 +
\oo(h) \norm{u}^2.
\end{equation}
We next use the symbolic calculus  and get from (\ref{sjodiezedvpt})
to first order and using the notation from there
\begin{equation} \label{truc2}
\psi_1(\lambda^2/|z|) \sharp \frac{{p-z}}{|z|} = \psi_1(\lambda^2/|z|)
\sep{\frac{{p-z}}{|z|}} + \frac{h}{|z|}R_1(\psi_1(\lambda^2/|z|), p).
\end{equation}
 Now observe that uniformly in $z \geq \oo(1)$ we have
$$
\D \sep{ \psi_1(\lambda^2/|z|) } \in S(\lambda^{-1}, dx^2
+ d\xi^2/\lambda^2),
\ \ \ \text{ and } \ \ \
\D p \in S(\lambda, dx^2 + d\xi^2/\lambda^2).
$$
The first fact follows from $\D \lambda \in S(1,dx^2 +
d\xi^2/\lambda^2)$  and the second from (\ref{symbcal}).
Consequently we
get a better estimate than the one that would be given by the classical
symbolic calculus in the class associated with the metric $dx^2 +
d\xi^2/\lambda^2$, namely
$$
R_1(\psi_1(\lambda^2/|z|), p) \in S(1, dx^2 + d\xi^2/\lambda^2).
$$
Since $|z| \geq \oo(1)$, we get that
$$
\frac{h}{|z|} R_1(\psi_1(\lambda^2/|z|), p) \in S(h,dx^2
+ d\xi^2/\lambda^2).
$$
Using this together with (\ref{truc1}, \ref{truc2}) yields
\begin{equation} \label{truc3}
\sep{ \psi_1^2(\lambda^2/|z|)^w u, u} \leq \norm{\sep{
\psi_1(\lambda^2/|z|)}^w \sep{ \frac{{p-z}}{|z|}}^w u}^2 +\sep{
\oo(h) + \oo(h^2)} \norm{u}^2.
\end{equation}
Since $\sep{ \psi_1(\lambda^2/|z|)}^w$ is bounded, we get the
lemma. \fin

\Section{Resolvent estimates away from the critical points when
$z$ is small}
 \label{awaysmall}

We work again in this section with  $p$ satisfying  \bf
(H2), (H3), (H4) \rm  and away from the critical points, but for
small $z$. Here the spectral parameter will be denoted $hz$ for $z
= \oo(1)$. We recall some notations of the preceding section,
namely
$$
\mu^2 = p_1 + (h\lambda)^{2/3}, \ \ \ \Gamma =\frac{dx^2}{h^{2/3}}
+ \frac{d \xi^2}{\mu^2}.
$$
As in the preceding section, we fix $\eps
>0$  and work with our operator $p$ satisfying
conditions (\ref{symbcal}-\ref{hyploin}). We can write for $u\in
\ss$
\begin{equation} \label{pdiezeg2}
\begin{split}
\Re \sep{ (p^w -hz)u, (1 -\eps g )^w u}
& = \Re \sep{ \sep{(1 -\eps g )\sharp (p -hz)}^w u,u} \\
& =    \sep{  ((p_1 -\Re hz)(1 -\eps g ) + \eps h \set{p_2, g}/2
  + \eps \Re r)^w u,u} ,
\end{split}
\end{equation}
following the same computations as in (\ref{pdiezeg}--\ref{major}).
 We also
get that
\begin{equation}
\begin{split}
\Re   (1 -\eps g )\sharp (p -hz)
& = (p_1 -\Re hz)(1 -\eps g ) + \eps h \set{p_2, g}/2 + \eps \Re r \\
& \geq c \mu^2 - 2\max (\Re (hz), 0)  - \mu ^2\phi ,
\end{split}
\end{equation}
where   $\phi \in \cc_0^\infty$ is equal to $1$ in a neighborhood
of the
critical points, and where we recall that
$r  \in S(h^{2/3}\mu,\Gamma)$ was
defined in (\ref{symbpg}). Of course outside this fixed neighborhood,
 and for
$h$ small enough, we have, using $\mu\geq h^{1/3}$,
$$
\mu^2 \gg 2\Re(hz),
$$
therefore with a new function $\phi$,
\begin{equation}
\begin{split}
\Re  (p -hz)\sharp (1 -\eps g ) & \geq c \mu^2/2   -\mu ^2 \phi .
\end{split}
\end{equation}
We can now use the Fefferman-Phong inequality (Lemma \ref{fphong}).
Indeed, each term is in $S(\mu^2, \Gamma)$ and we get

\begin{equation}
((\mu^2)^wu,u) \leq C\Re((p^w-hz)u, (1-\eps g)^w u) + C((\mu ^2\phi )^w
u,u) +  \Re(R^w
u, u),
\end{equation}
where $ R$ is of order $h^2 \times \mu^2 \times(\mu h^{1/3})^{-2}
= h^{4/3}$
from lemma \ref{fphong} (recall that $\mu h^{1/3}$ is the uncertainty
parameter
of $\Gamma$). Choosing $h$ small enough and noticing that
$$
(\mu^2)^w \geq ch^{2/3} ,
$$
gives
\begin{equation} \label{bonestim2}
ch^{2/3}\norm{u}^2 \leq \sep{(\mu^2)^wu,u)} \leq C\Re((p^w-hz)u,
(1-\eps g)^w u) +
 C((\mu ^2\phi )^w u,u).
\end{equation}
We next write that  for $\eps$ sufficiently small,
$$
\Re((p^w-hz)u, (1-\eps g)^w u) \leq \norm{(p^w-hz)u}\norm{u}.
$$
From this and (\ref{bonestim2}) we get the main result of this
section
\begin{equation} \label{main9}
ch^{2/3} \norm{u} \leq C\norm{(p^w-hz)u} + C\norm{(\mu ^2\phi )^wu}
\end{equation}

\Section{Proof of Theorem \ref{main}}\label{SectProofMain}

In this section we shall glue together all the results of the
Sections \ref{locallarge}, \ref{localsmall}, \ref{awaylarge} and
\ref{awaysmall}. We give the results here in the original
variables and not on the FBI side.

In the following, we choose  $u \in \ss$ and  we write $U=Tu$
where $T$ is the FBI-Bargmann transform associated with the phase
$i(x-y)^2/2$. We also denote by $P$ the operator $(\chi_0p)^w$ on
the FBI side, where $\chi_0$ is some $\cc^\infty_0$ function equal
to $1$ in a very large compact set (including the critical
points).

 \preuve[of a)] We suppose here that $h |z| \leq \oo(h)$.
Let us first recall the main result (\ref{main9}) of Section
\ref{awaysmall}:
\begin{equation} \label{philarge}
h^{2/3} \norm{u} \leq C \norm{(p^w - hz)u} + C\norm{(\mu^2\phi)^w
u} ,
\end{equation}
where $\phi$ is a cutoff function equal to 1 near the critical
points. We choose once and for all another cutoff function $\psi$
equal to one in a larger neighborhood of the critical points,
so that $\nabla \phi \nabla \psi = 0$. Then
\begin{equation*} \label{philarge2}
h^{2/3} \norm{(1-\psi)^w u} \leq C \norm{(p^w - h z)(1-\psi)^wu} +
C\norm{(\mu^2\phi)^w (1-\psi)^w u}.
\end{equation*}
Notice that $(\mu^2\phi)^w (1-\psi)^w = \oo(h^\infty)$ as a
bounded operator in $L^2$ since the supports are disjoint.
Moreover,
\begin{equation} \label{philarge2bis}
(p^w - h z)(1-\psi)^w = (1-\psi)^w(p^w - h z) + \frac{h}{2i}
\set{\psi, p}^w + \oo(h^2),
\end{equation}
where $q \defegal \frac{1}{2i} \set{\psi, p}$ is a
symbol with $\textrm{supp\,}q\subset \textrm{supp\,}\nabla \psi$,
 so that
the support of $q$ is disjoint from the support of $\phi$. Hence
\begin{equation*} 
h^{2/3} \norm{(1-\psi)^w u} \leq C \norm{(1-\psi)^w(p^w - h z)u} +
Ch\norm{q^w u} + \oo(h^2)\norm{u}.
\end{equation*}
The $L^2$-boundedness of $(1-\psi)^w$ and the fact that $h\le h^{2/3}$
give
\begin{equation} \label{philarge3}
h \norm{(1-\psi)^w u} \leq C \norm{(p^w - h z)u} + Ch\norm{q^w u}
+ \oo(h^2)\norm{u} .
\end{equation}

The main result of Section \ref{localsmall} on the
FBI side states that
\begin{equation} \label{chismallbis}
h \norm{ U}_{\Phi_0} \leq \norm{(P-hz)U}_{\Phi_0}+ h^{5/6}
\norm{(1-\chi) U}_{\Phi_0},
\end{equation}
where $\chi$ is an arbitrary cutoff function equal to $1$ in a
neighborhood of the critical points. We can  choose $\chi$ equal
to 1 in a neighborhood of $\textrm{supp\,}\psi$, where $\psi $ is
viewed as a function on the FBI-side (i.e. $\psi\circ \kappa^{-1}$
where $\kappa$ is the canonical transform associated with the FBI
transform $T$). With these notations we may write that $\phi \prec
\psi \prec \chi \prec \chi_0$ modulo a composition with $\kappa$.
Coming back to the real side for the two first terms of this
inequality, and using the metaplectic invariance gives
\begin{equation} \label{chismall2}
h \norm{u} \leq \norm{((\chi_0p)^w-hz)u} + h^{5/6} \norm{(1-\chi)
U}_{\Phi_0},
\end{equation}
and after replacing $u$ by  $\psi^w u$,
\begin{equation} \label{chismall3}
h \norm{\psi^w u} \leq \norm{((\chi_0p)^w-hz) \psi^w u} + h^{5/6}
\norm{(1-\chi) T \psi^w u}_{\Phi_0}.
\end{equation}
Now we can treat the term $\norm{((\chi_0p)^w-hz) \psi^w u}$ as in
(\ref{philarge2bis}) and get rid of the term $\chi_0$ modulo a
term of order $h^\infty$ and we get with the same $q$
\begin{equation} \label{chismall4}
h \norm{\psi^w u} \leq \norm{(p^w-hz) u} + h\norm{q^w u} +
\oo(h^2)\norm{u} + h^{5/6} \norm{(1-\chi) T \psi^w u}_{\Phi_0}.
\end{equation}

We shall use the following standard lemma for which we briefly review
the proof at the end of this section.

\begin{lemma} \label{chiTpsi}
We have $\norm{(1-\chi) T \psi^w u}_{\Phi_0} = \oo(h^\infty)
\norm{u}$.
\end{lemma}

We can therefore write
\begin{equation} \label{chismall5}
h \norm{\psi^w u} \leq \norm{(p^w-hz) u} + h\norm{q^w u} +
\oo(h^2)\norm{u} .
\end{equation}
Let us now glue together the results
(\ref{philarge3}), (\ref{chismall5}) to get
\begin{equation} \label{finalsmall1}
h \norm{(1-\psi)^w u} + h \norm{\psi^w u} \leq C \norm{(p^w - h
z)u} + Ch\norm{q^w u} + \oo(h^2)\norm{u}.
\end{equation}
For the term $Ch\norm{q^w u}$ we simply
apply (\ref{philarge}) with $u$ replaced by $q^w u$. This gives
\begin{equation} \label{q1}
h^{2/3} \norm{q^w u} \leq C \norm{(p^w - hz)q^w u} + \norm{\phi^w
q^w u}.
\end{equation}
Since $\phi$ and $q$ have disjoint support, we
have $\phi^w q^w = \oo(h^\infty)$ as an operator in $L^2$. Besides
we have
$$
(p^w - hz)q^w = q^w (p^w - hz) + \oo(h),
$$
since $q$ is with compact support. Therefore we get
\begin{equation} \label{q2}
h^{2/3} \norm{q^w u} \leq C \norm{q^w (p^w - hz) u} + C h \norm{
u} \leq \norm{(p^w - hz) u} + C h \norm{ u},
\end{equation}
and eventually
\begin{equation} \label{q3}
h \norm{q^w u}  \leq \norm{(p^w - hz) u} + C h^{4/3} \norm{ u}.
\end{equation}
Together with (\ref{finalsmall1}) this  yields
\begin{equation} \label{finalsmall1bis}
h \norm{(1-\psi)^w u} + h \norm{\psi^w u} \leq C \norm{(p^w - h
z)u} +  \oo(h^{4/3})\norm{u},
\end{equation}
and using the triangle inequality
$
\norm{u}\le \norm{(1-\psi)^w u} +  \norm{\psi^w u} $,
\begin{equation} \label{finalsmall2}
h \norm{ u} \leq C \norm{(p^w - h z)u} + \oo(h^{4/3})\norm{u}.
\end{equation}
Taking $h$ small enough completes the proof of part a) of the
theorem. \fin

\preuve[of b)] In this section we suppose that $|z| \gg h$. We
also denote in the following
$$
Z = |z|^{1/3} h^{2/3} .
$$
 We shall follow the proof of part a). We first recall the main
  result (\ref{main8}) of Section \ref{awaylarge}:
\begin{equation} \label{phismall}
Z \norm{u} \leq C  \norm{(p^w - z)u} +
 C\norm{(\mu ^2\phi )^wu}
,
\end{equation}
where $\phi$ is a cutoff function equal to 1 near the critical
points. As in the preceding section we  choose once and for all
another cutoff function $\psi$ such that $\psi \succ \phi$ and we
write
\begin{equation*} \label{phismall2}
Z \norm{(1-\psi)^w u} \leq C \norm{(p^w - z)(1-\psi)^wu} +
C\norm{(\mu ^2\phi )^w (1-\psi)^w u} .
\end{equation*}
As in
(\ref{philarge2bis}), (\ref{philarge3}) we get
\begin{equation} \label{phismall3}
Z \norm{(1-\psi)^w u} \leq C \norm{(p^w - z)u} + Ch\norm{q^w u} +
\oo(h^2)\norm{u},
\end{equation}
where we recall  $q \defegal \frac{1}{2i} \set{\psi, p}$ is a
symbol with support in $\nabla \psi$.

We now recall the main result of Section \ref{locallarge} on the
FBI side (see equation (\ref{6.23})):
\begin{equation} \label{chilargebis}
Z \norm{ U}_{\Phi_0} \leq C \sep{\norm{(P-z)U}_{\Phi_0}  + Z
\norm{(1-\chi) U}_{\Phi_0}},
\end{equation}
where $\chi$ is an arbitrary cutoff function equal to $1$ in a
neighboorhood of the critical points. We   choose $\chi \succ
\psi$ where  $\psi$ is viewed as a function on the FBI side. With
this notation we  write as in the proof of a)   that $\phi \prec
\psi \prec \chi \prec \chi_0$. Coming back to the real side for
the two first terms of this inequality, and using the metaplectic
invariance gives
\begin{equation} \label{chilarge2}
Z \norm{u} \leq C\sep{\norm{((\chi_0p)^w-z)u} + Z \norm{(1-\chi)
U}_{\Phi_0}}.
\end{equation}
Taking $\psi^w u$ instead of $u$ gives
\begin{equation} \label{chilarge3}
Z \norm{\psi^w u} \leq C \bigg(\norm{((\chi_0p)^w-z) \psi^w u} + Z
\norm{(1-\chi) T \psi^w u}_{\Phi_0}\bigg).
\end{equation}
Now we can treat the term $\norm{(p^w-z) \psi^w u}$ as in
(\ref{philarge2bis}) and we get with the same $q$
\begin{equation} \label{chilarge4}
Z \norm{\psi^w u} \leq C \bigg( \norm{(p^w-z) u} + h\norm{q^w u} +
\oo(h^2)\norm{u} + Z \norm{(1-\chi) T \psi^w u}_{\Phi_0}\bigg).
\end{equation}
Using Lemma \ref{chiTpsi} yields,
\begin{equation} \label{chilarge5}
Z \norm{\psi^w u} \leq C\norm{(p^w-z) u} + C h\norm{q^w u} +
\oo(h^2)\norm{u} + Z \oo(h^{\infty})\norm{u} .
\end{equation}

Let us now combine
(\ref{phismall3}), (\ref{chilarge5}):
\begin{equation} \label{finallarge1}
Z \norm{(1-\psi)^w u} + Z \norm{\psi^w u} \leq C \norm{(p^w - z)u}
+ Ch\norm{q^w u} + \oo(h^2)\norm{u} + Z\oo(h^\infty)\norm{u}.
\end{equation}
Now we can use (\ref{q2}) and we get with new constants
\begin{equation} \label{finallarge1bis}
Z \norm{(1-\psi)^w u} + Z \norm{\psi^w u} \leq C \norm{(p^w - z)u}
+  \oo(h^{4/3})\norm{u} +\oo(h^2)\norm{u} + Z\oo(h^\infty)\norm{u},
\end{equation}
and the triangle inequality gives
\begin{equation*} \label{finallarge2}
Z \norm{ u} \leq C \norm{(p^w - z)u}.
\end{equation*}
The proof of part b) of the theorem is complete. \fin

\preuve[of lemma \ref{chiTpsi}]
We have
\[
  \norm{ (1-\chi) T \psi^w u}_{\Phi_0}^2
  = ( u, \psi^w T^* (1-\chi)^2 T \psi^w u)
  \leq \norm{u} \, \norm{\psi^w T^* (1-\chi)^2 T \psi^w u} ,
\]
where the adjoint $T^*$ is w.r.t.\ the $\Phi_0$ inner product.
We will show that $T^*(1-\chi)^2 T$ is a pseudo\-differential operator
with Weyl symbol that is $\oo(h^\infty)$ where $(1-\chi)^2 \circ
\kappa$ and all its derivatives vanish. Since $\psi$ and
$(1-\chi) \circ \kappa$ have disjoint support, this shows
the result.

To simplify the notation we do the computations for $T^* \chi T$.
Recall that we use the transform (\ref{7.7}), with $\phi(t,y) =
i(t-y)^2/2$. The constant $C$ in (\ref{7.7}) is given by
$2^{-\frac{n}{2}} \pi^{-\frac{3n}{4}}$. The function $\Phi_0$
equals $-(\Im t)^2/2$ and $(\chi \circ \kappa)(x,\xi) =
\chi(x-i\xi)$. If we write $t = x - i \xi$, then the kernel
$K(y,z)$ of $T^* \chi T$ is given by
\begin{equation} \label{eq:K1}
  K(y,z) = 2^{-n} (\pi h)^{-\frac{3n}{2}}
  \int e^{-i(x-y)\cdot\xi/h-(x-y)^2/(2h)
  + i(x-z)\cdot\xi/h -(x-z)^2/(2h)}
    \chi(x-i\xi) \, dxd\xi .
\end{equation}
The phase function can be written as
\[
  i(y-z) \cdot \xi/h - (x-\frac{y+z}{2})^2/h - (y-z)^2/(4h) .
\]
For the last term in this expression we have from a Fourier
transformation
\[
  e^{-(y-z)^2/(4h)}
  = (\pi h)^{-n/2} \int e^{i (y-z)\cdot(\eta-\xi)/h - (\eta-\xi)^2/h}
    \, d\eta .
\]
Entering this in (\ref{eq:K1}) we find that $K(y,z)$ equals
\[
  2^{-n} (\pi h)^{-2n}
  \int e^{i(y-z)\cdot\eta/h - (x-\frac{y+z}{2})^2/h - (\eta-\xi)^2/h}
  \chi(x-i\xi) \, dx d\xi d\eta .
\]
This formally equals a Weyl pseudodifferential operator with symbol
\[
  \tilde{\chi}(y,\eta) = (\pi h)^{-n}
    \int e^{-(y-x)^2/h-(\eta-\xi)^2/h} \chi(x-i\xi) \, dx d\xi .
\]
It is clear that $\tilde{\chi}$ has the correct symbol property, and
that $\tilde{\chi}(y,\eta) = \oo(h^\infty)$ for $(y,\eta)$ such that
$(\chi \circ \kappa)(y,\eta)$ and all its derivatives vanish. This
completes the proof. \fin

\Section{Asymptotic expansion of the eigenvalues.}\label{SectEv}
\setcounter{equation}{0}
\subsection{From apriori estimates to the resolvent.}\label{Subev1}

In the previous sections we obtained apriori estimates for $z$ in
a subset of $\C$, given by
\begin{equation} \label{eq:gen_res_est}
  \| u \| \le C \| (P - z)u \| , \ \
    \forall u \in {\cal S}(\R^n) .
\end{equation}
We will show that such estimates imply the existence of the
resolvent of $P$.

We will first establish this for one particular value of $z$. For
this purpose we will use some functional analysis, and results
given in section 5.2 of \cite{HelNi}. Following \cite{HelNi} we
define $P : {\cal D}(P) \rightarrow L^2(\R^n)$ with domain ${\cal
D}(P) = \cc_0^\infty(\R^{n})$. We let $\overline{P}$ be its
closure (further on we will simply write $P$ instead of
$\overline{P}$ but for the moment we keep the distinction).

We show that ${\cal S}(\R^n) \subset {\cal D}(\overline{P})$. This
follows if for $u \in {\cal S}(\R^n)$ there is a sequence $u_j \in
\cc_0^\infty(\R^n)$ with $u_j \rightarrow u$ in $L^2(\R^n)$ and
$Pu_j \rightarrow Pu$ in $L^2(\R^n)$. Such a sequence is given by
$u_j = \chi(\frac{x}{j}) u(x)$, where $\chi \in
C_0^\infty(\R^n,[0,1])$ is equal to $1$ on a neighborhood of $0$.
We have $Pu_j = \chi(\frac{\cdot}{j}) P u + [
P,\chi(\frac{\cdot}{j})] u \rightarrow Pu$, 
since the symbol of the commutator tends to zero in $S(\lambda ,\Gamma _0)$. 
By the definition of $\overline{P}$ we
have that in fact ${\cal S}(\R^n) \text{ is dense in } {\cal
D}(\overline{P})$.

Next we establish the existence of the resolvent for at least one
value $z_0$ in the left complex half plane, when there is a real
$\lambda_0$ such that $\overline{P}+\lambda_0$ is maximally
accretive. For the Kramers-Fokker-Planck operator this property is
established in proposition 5.5 of \cite{HelNi} with $\lambda_0 =
0$.

\begin{prop}
Assume that $\overline{P}+\lambda_0$ is maximally accretive. Then
there is $\lambda_1 > \lambda_0$ such that $(\overline{P} +
\lambda_1)^{-1}$ exists and is a bounded operator on $L^2(\R^n)$.
\end{prop}

\preuve The accretivity of $\overline{P}+\lambda_0$ means that
$((\overline{P}+\lambda_0)u,u) \ge 0$ for each $u \in {\cal
  D}(\overline{P})$. It follows that for each $\lambda > \lambda_0$
  we
have
\[
  \|(\overline{P} + \lambda) u\| \| u \| \ge
  ((\overline{P}+\lambda_0)u,u) + (\lambda-\lambda_0) \|u\|^2
  \ge  (\lambda-\lambda_0) \|u\|^2 , \ \ u\in {\cal D}(\overline{P}) ,
\]
hence
\begin{equation} \label{eq:accretive_estimate}
  \|u\| \le (\lambda-\lambda_0)^{-1} \|(\overline{P} + \lambda) u\| ,
    \ \ u\in {\cal D}(\overline{P}) .
\end{equation}
Hence $\overline{P}+\lambda$ is injective.

Suppose now that there is a sequence $u_j \in {\cal S}(\R^n)$ such
that $(\overline{P}+\lambda) u_j  \rightarrow v$ in $L^2(\R^n)$
for some $v \in L^2(\R^n)$. Denote $v_j = (\overline{P}+\lambda)
u_j$. Then by the estimate (\ref{eq:accretive_estimate}) it
follows that $\| u_j - u_k \| \rightarrow 0, j,k\rightarrow
\infty$, hence $u_j$ converges to an element $u$ in $L^2(\R^n)$.
Now $(u_j,v_j) \in \operatorname{graph}(\overline{P})$ and $u_j
\rightarrow u, v_j \rightarrow v$ in $L^2(\R^n)$. Therefore, the
range ${\cal R}(\overline{P})$ is closed. Theorem 5.4 of
\cite{HelNi} and the fact that $\overline{P}$ is
  maximally accretive imply that for some $\lambda_1 > \lambda_0$,
  the
range of $P+\lambda_1$ is also dense in $L^2(\R^n)$. It follows
that $\overline{P}+\lambda_1$ is surjective, that the inverse
$(\overline{P} + \lambda_1)^{-1}$ exists and that its norm is
bounded by $\frac{1}{\lambda_1-\lambda_0}$. \fin

\begin{remark}
Alternatively we could use the following additional properties
\begin{align}
  & \| u \| \le C \| (P^* - \overline{z})u \| , \ \
    \forall u \in {\cal S}(\R^n) \label{eq:res_adj_est} , \\
  & u \in L^2(\R^n) , (P-z) u \in {\cal S}(\R^n)
    \Rightarrow u \in {\cal S}(\R^n) \label{eq:res_add_prop} .
\end{align}
(Here we let ${\cal D}(P) = \{ u \in L^2(\R^n) \, ;\, Pu \in
L^2(\R^n)\}$.) The first property is similar to the apriori
estimate (\ref{eq:gen_res_est}). The second property can for
example be derived from hypoelliptic estimates in a chain of
weighted Sobolev spaces as given for the Kramers-Fokker-Planck
case in theorem 3.1d) of \cite{HerNi03} (a result valid under
somewhat different conditions than used here). In short the
argument using (\ref{eq:res_add_prop}) goes as follows. By a
standard argument estimate (\ref{eq:res_adj_est}) and the
Hahn-Banach theorem imply the surjectivity of $P-z$. If $u \in
{\cal D}(P)$ and $(P-z)u=0$, then (\ref{eq:res_add_prop}) implies
that $u \in {\cal S}(\R^n)$ and (\ref{eq:gen_res_est}) that $u=0$.
Hence $P : {\cal D}(P) \rightarrow L^2(\R^n)$ is injective, and
$(P-z)^{-1} : L^2 \rightarrow L^2$ is bounded and $\|(P-z)^{-1} \|
\le C$. One can also show that under these assumptions ${\cal
S}(\R^n)$ is dense in ${\cal D}(P)$ for the graph norm.
\end{remark}

From now on we simply write $P$ instead of $\overline{P}$. To
obtain the resolvent, we consider an abstract situation. Let
$P:L^2({\R}^n)\to L^2({\R}^n)$ be a closed operator and assume (as
we established above) \ekv{ev1.1} { {\cal S}({\R}^n) \hbox{ is
dense in }{\cal D}(P). } Since ${\cal D}(P)$ has the norm $\Vert
u\Vert _{{\cal D}(P)}=\Vert u\Vert +\Vert Pu\Vert $, this means
that for every $u\in {\cal D}(P)$, there is a sequence
$u_j\in{\cal S}$, $j=1,2,...$, such that $u_j\to u$ and $Pu_j\to
Pu$ in $L^2$.

Let $\Omega \subset \C$ be a connected open set. Let $z_0\in
\Omega $ an assume \ekv{ev1.2} {(z_0-P)^{-1}:L^2\to {\cal
D}(P)\hbox{ exists,}} \ekv{ev1.3} { \Vert u\Vert \le C_K\Vert
(P-z)u\Vert ,\ \forall u\in{\cal S},\, z\in K, } for every
$K\subset\subset \Omega $.
\begin{prop}\label{Propev1.1} Under these assumptions,
$(z-P)^{-1}$ exists for every $z\in\Omega $.
\end{prop}

\preuve  Using (\ref{ev1.1}), we see that the apriori estimate in
(\ref{ev1.3}) extends to all $u\in{\cal D}(P)$. In particular,
$z-P:{\cal D}(P)\to L^2$ is injective for $z\in\Omega $, so it
remains to show that $z-P$ is surjective. If $(z_1-P)^{-1}$ exists
for some $z_1\in K\subset\subset \Omega $, then (\ref{ev1.3})
(extended to ${\cal D}(P)$) implies that $\Vert (z_1-P)^{-1}\Vert
\le C_K$. Hence $\Vert (z-z_1)(z_1-P)^{-1}\Vert <1$ for $\vert
z-z_1\vert <1/C_K$, and we conclude that $z-P:{\cal D}(P)\to L^2$
has a right inverse of the form
$(z_1-P)^{-1}(1+(z-z_1)(z_1-P)^{-1})^{-1}$. If in addition, $z\in
\Omega $, this right inverse is equal to the resolvent.

\par If $z\in \Omega $ is any given point, we take a smooth curve
$\gamma $ in $\Omega $ from $z_0$ to $z$, and cover $\gamma $ by
finitely many discs $D(z_j,r)$, $j=0,1,2,...,M$, such that
$r<1/C_\gamma $, $z_{j+1}\in D(z_j,r)$. Hence $(z-P)^{-1}$
exists.\fin

\par The same result is valid for Grushin problems. We keep the
initial hypothesis about $P$. Let $R_-:{\C}^{N_0}\to L^2$,
$R_+:{\cal D}(P)\to {\C}^{N_0}$ be \bdd{} and for simplicity
\indep{} of $z$. Put \ekv{ev1.4} { {\cal P}(z)=\begin{pmatrix}P-z
&R_-\cr R_+ &0\end{pmatrix}:{\cal D}(P)\times {\C}^{N_0}\to
L^2\times {\C}^{N_0}. } Assume still that (\ref{ev1.2}) holds for
some $z_0\in \Omega $. Instead of (\ref{ev1.3}), we assume
\ekv{ev1.5} { \Vert u\Vert +\vert u_-\vert \le C_K(\Vert
(P-z)u+R_-u_-\Vert +\vert R_+u\vert ),\ z\in K,\, u\in{\cal S},\,
u_-\in{\C}^{N_0}, } for every $K\subset\subset \Omega $. (Again,
this extends to the case $u\in{\cal D}(P)$.)
\begin{prop}\label{Propev1.2} \it Under the above assumptions,
${\cal
P}(z)$ has a \bdd{} inverse for every $z\in\Omega $.
\end{prop}

\preuve \rm As before, we notice that (\ref{ev1.5}) implies that
\ekv{ev1.6} { \Vert u\Vert _{{\cal D}(P)}+\vert u_-\vert \le
C_K(\Vert (P-z)u+R_-u_-\Vert +\vert R_+u\vert),\ u\in{\cal
D}(P),\, u_-\in{\C}^{N_0},\, z\in K, } with a new constant $C_K$,
so ${\cal P}(z)$ is injective for all $z\in \Omega $.

\par For $z=z_0$, $(P-z_0):{\cal D}(P)\to L^2$ has a \bdd{}
inverse and is
therefore a Fredholm \op{} of index 0. Hence,
$$
Q:=\begin{pmatrix}P-z_0 &0\cr 0 &O\end{pmatrix}:{\cal D}(P)\times
{\C}^{N_0}\to L^2\times {\C}^{N_0}
$$
is Fredholm of index 0 and ${\cal P}(z_0) $ has the same property,
being a finite  rank \pert{} of $Q$. Being injective by
(\ref{ev1.5}), it is bijective, and as in the preceding proof, we
see that ${\cal P}(z)^{-1}$ exists for all $z\in \Omega $ with
$\vert z-z_0\vert <1/C_K$ if $z_0\in K$. By the same procedure as
above, we get the result. \fin

\subsection{Grushin \pb{} in the quadratic case}\label{Subev2}

\par Let $P_0$ be a quadratic \op{} on $L^2({\R}^n)$,
so that $P_0$ has the Weyl symbol $\sum_{\vert \alpha +\beta \vert
=2}a_{\alpha ,\beta }x^\alpha \xi ^\beta $ that we also denote by
$P_0(x,\xi )$.  (We can also add a constant to our symbol, but we
shall avoid for simplicity to have linear terms in the symbol.) As
in \cite{Sjo74} we assume that $P_0$ is elliptic away from (0,0):
\ekv{ev2.1} { P_0(x,\xi )\ne 0,\ (x,\xi )\in {\R}^{2n}\setminus\{
(0,0)\} . } When $n>1$ this implies that $P_0({\R}^{2n})$ is a
proper cone in ${\C}$ and when $n=1$ we assume that so is the
case. Then $P_0$ is a closed $\op{}:L^2\to L^2$ with domain ${\cal
D}(P_0)=\langle (x,D)\rangle ^{-2}(L^2)$ and the assumption
(\ref{ev1.1}) is fulfilled. $P_0$ has discrete spectrum and the
\ev{}s are computed in \cite{Sjo74} as recalled in Section
\ref{SectQuad}. They are contained in $P_0({\R}^{2n})$.

\par Let $\lambda _0\in{\C}$ be such an \ev{} and let $E_{\lambda
_0}\subset {\cal D}(P_0)$ be the corresponding space of
generalized eigenvectors. Let $e_1,...,e_{N_0}$ be a basis for
$E_{\lambda _0}$ and let $f_1,...,f_{N_0}\in{\cal S}({\R}^n)$ have
the property that \ekv{ev2.2} {\det ((e_j\vert f_k))\ne 0.} A
possibly natural choice would be to let $f_1,...,f_{N_0}$ be the
dual basis in the space $E^*_{\overline{\lambda }_0}$ of
generalized eigenvectors of $P^*$, associated to the \ev{}
$\overline{\lambda }_0$.

\par Put
$$
R_-u_-=\sum u_-(j)e_j,\ R_+u=((u\vert f_j))\in{\C}^{N_0},
$$
for $u_-=(u_-(j))\in{\C}^{N_0}$. For $\lambda \in{\rm
neigh\,}(\lambda _0)$, the \pb{} \ekv{ev2.3} { (P_0-\lambda
)u+R_-u_-=v,\ R_+u=v_+, } has a unique solution $(u,u_-)\in {\cal
D}(P_0)$, for every $(v,v_+)\in L^2\times {\C}^{N_0}$. In fact,
let $\Pi :L^2\to E_{\lambda _0}$ be the spectral projection and
decompose $u=u'+u''$, $v=v'+v''$, with $u''=\Pi u$, $u'=(1-\Pi )u$
and similarly for $v$. Then the \e{} for $u'$ is $(P_0-\lambda
)u'=v'$ and determines $u'\in {\cal D}(P)$ uniquely. $u''$ is
completely determined by the condition $R_+u''=v_+-R_+u'$, thanks
to the assumption (\ref{ev2.2}). Finally $u_-$ is determined by
$R_-u_-=v'-(P_0-\lambda )u'$.

\par If we introduce the solution \op{}
$${\cal E}=\begin{pmatrix}E &E_+\cr E_- &E_{-+}\end{pmatrix},
\hbox{ by
}
\begin{pmatrix}u\cr u_-\end{pmatrix}={\cal
E}\begin{pmatrix}v\cr v_+\end{pmatrix},$$ then we also know that
\ekv{ev2.4} {
\begin{pmatrix}\Lambda ^{2-k}&\cr 0 &1\end{pmatrix}{\cal E}
\begin{pmatrix}\Lambda ^k&0\cr
0&1\end{pmatrix}\hbox{ is \bdd{}}} for every $k\in{\R}$, when
$\Lambda =\langle (x,D)\rangle $. We also notice that if
$M(\lambda )$ denotes the matrix of ${(\lambda -P_0)_\vert
}_{E_{\lambda _0}}$ \wrt{} the basis $e_1,...,e_{N_0}$, then
\ekv{ev2.5} { E_{-+}(\lambda )=M(\lambda )((e_k\vert f_j))^{-1}. }

\par We now choose $P_0$ as in Proposition \ref{quad}, acting on
 $H_{\Phi
_\epsilon ^0}$, where $\Phi _{\epsilon }^0$ is a quadratic form
with $\Lambda _{\Phi _\epsilon ^0}=\kappa _T(\Lambda _{\epsilon
G^0})$ and $G^0$ is a real quadratic form chosen as in
\ref{gzero}. Here $\epsilon >0$ is small and fixed and the earlier
assumptions are fulfilled with the real phase space replaced by
$\Lambda _{\Phi _\epsilon ^0}$. As in Section \ref{SectQuad} we
now work with the $h$-quantization. Then, if $\lambda _0$ is an
\ev{} of the ($h=1$) quantization, we get the well-posed Grushin
\pb{} \ekv{ev2.6} { (P_0-hz)u+R_-u_-=v,\ R_+u=v_+, } for $z$ in
some fixed \neigh{} of $\lambda _0$. Here, we take \ekv{ev2.7} {
R_+u(j)=(u\vert f_{j,h})_{L^2_{\Phi _\epsilon ^0}},\ R_-u_-\sum
u_-(j)e_{j,h}, } with $f_{j,h}(x)=h^{-{n\over 2}}f_j({x\over
\sqrt{h}})$, and similarly for $e_{j,h}$, so that
$$
R_+={\cal O}(1):L^2_{\Phi _0^\epsilon }\to {\C}^{N_0},\ R_-={\cal
O}(1):{\C}^{N_0}\to H_{\Phi _0^\epsilon },
$$
\ufly{}, when $h\to 0$. More precisely, we have (cf. Proposition
\ref{quad}):
\begin{prop}\label{Propev2.1} \it For every $(v,v_+)\in H_{\Phi
_0^\epsilon }\times {\C}^{N_0}$, the \pb{} (\ref{ev2.6}) has a
unique solution in the same space and the solution satisfies:
$d^2u\in L_{\Phi _0^\epsilon }^2$. Moreover, for every fixed $k\in
{\R}$, we have the apriori estimate \ekv{ev2.8} { h\Vert
(1+{d^2\over h})^{1-k}u\Vert +\vert u_-\vert \le C(\Vert
(1+{d^2\over h})^{-k}v\Vert +h\vert v_+\vert ). }\end{prop}

\preuve When $h=1$, we simply translate the earlier result for
(\ref{ev2.3}) into a result for (\ref{ev2.6}) and get the estimate
\ekv{ev2.9} { \Vert (1+d^2)^{1-k}u\Vert +\vert u_-\vert \le
C(\Vert (1+d^2)^{-k}v\Vert +\vert v_+\vert ). } Now consider
(\ref{ev2.6}) for other values of $h$, and indicate the
$h$-dependence by means of super/sub-scripts.  Let
$Uf(x)=h^{n\over 2}f(\sqrt{h}x)$, so that $U$ is unitary: $H_{\Phi
_0^\epsilon ,h}\to H_{\Phi _0^\epsilon ,1}$. We further have
$$
UP_0^h=hP_0^1U,\ UR_-^h=R_-^1,\ R_+^h=R_+^1U,
$$  and the \pb{} (\ref{ev2.6}) (with general $h$) can be \tf{}ed into
\ekv{ev2.10} { h(P_0^1-z)Uu+R_-^1u_-=Uv,\ R_+^1Uu=v_+, } that we
write as \ekv{ev2.11} { (P_0^1-z)hUu+R_-^1u_-=Uv,\ R_+^1hUu=hv_+.
} Applying (\ref{ev2.9}) to this system, we get
$$
h\Vert (1+d^2)^{1-k}Uu\Vert +\vert u_-\vert \le C(\Vert
(1+d^2)^{-k}Uv\Vert +h\vert v_+\vert ).$$ Here
$$
d(x)Uu(x)=U(d({x\over \sqrt{h}})u(x))=U({d(x)\over \sqrt{h}}u(x)),
$$
and using the unitarity of $U$, we get (\ref{ev2.8}).\fin

\par We can rewrite (\ref{ev2.8}) equivalently as
\ekv{ev2.12} { \Vert (h+d^2)^{1-k}u\Vert +h^{-k}\vert u_-\vert \le
C(\Vert (h+d^2)^{-k}v\Vert +h^{1-k}\vert v_+\vert ). } In the
following, it will be convenient to replace the $f_j$ in the
definition of $R_+^1$ by $\chi _kf_j$, where $\chi _R(x)=\chi
({x\over R})$ for some \sufly{} large $R>0$. Correspondingly,
$f_j^h$ is replaced by $\chi _R({x\over \sqrt{h}})f_j^h$. This
will be only a small modification of $R_+^h$ and does neither
affect the well-posedness of (\ref{ev2.6}) nor the estimates
(\ref{ev2.8}), (\ref{ev2.12}).

\par Mimicking Proposition \ref{Prop7.1}, we have
\medskip
\begin{prop}\label{Propev2.2} Let $\chi _0\in C_0^\infty
({\C}^n)$ be fixed and $=1$ near $x=0$, and fix $k\in{\R}$. Then
for $z$ in a \neigh{} of $\lambda _0$, \indep{} of $k$, we have
the following estimate for the problem  (\ref{ev2.6}) in $H_{\Phi
_0^\epsilon }$ (for $\epsilon >0$ small and fixed and for $h$
\sufly{} small): \ekv{ev2.13} { \Vert (h+d^2)^{1-k}\chi _0u\Vert
+h^{-k}\vert u_-\vert \le C(\Vert (h+d^2)^{-k}\chi _0v\Vert
+h^{1-k}\vert v_+\vert +h^{1\over 2}\Vert 1_Ku\Vert ), } where $K$
is any fixed \neigh{} of ${\it supp\,}\chi _0$.\end{prop}

\preuve \rm Let $\Pi $ denote the orthogonal projection onto the
\hol{} \fu{}s as in Section \ref{SectQuad}. Applying $\Pi \chi _0$
to the first \e{} in (\ref{ev2.6}), we get after some simple
calculations, (using also that $\Pi R_-=R_-$, $u=\Pi u$):
\begin{equation}
\label{ev2.14}
\begin{split}
(P_0-hz)\Pi \chi _0u+R_-u_- & =\Pi \chi _0v+[P_0,\Pi \chi _0]u+\Pi
(1-\chi _0)R_-u_-, \\
R_+\Pi \chi _0u & = v_+-R_+(1-\chi _0)u-R_+(1-\Pi )\chi _0u.
\end{split}
\end{equation}

\par Here (\ref{7.42})  tells us that
$$
\Vert (h+d^2)^{-k}[P_0,\Pi \chi _0]u\Vert \le {\cal O}(h)\Vert
1_Ku\Vert .
$$
Since the $e_j$ decay exponentially and 
$(h+d^2)^{-k}\Pi (h+d^2)^{k}$ is \ufly{} \bdd{} in our weighted $L^2$ space,  
it is also clear that
$$
\Vert (h+d^2)^{-k}\Pi (1-\chi _0)R_-u_-\Vert \le {\cal O}(h^\infty
)\vert u_-\vert ,
$$
and since $\chi f_j$ has compact support, we have $R_+(1-\chi
_0)=0$, when $h>0$ is small enough. We also have $\vert R_+(1-\Pi
)\chi _0u\vert \le {\cal O}(h^\infty )\Vert 1_Ku\Vert $. Applying
this and (\ref{ev2.12}) to the \pb{} (\ref{ev2.14}), we get
\begin{eqnarray}\label{hd2}\Vert (h+d^2)^{1-k}\Pi \chi _0u\Vert
+h^{-k}\vert u_-\vert \le& \Vert (h+d^2)^{-k}\Pi \chi _0v\Vert
+{\cal O}(h)\Vert 1_Ku\Vert +{\cal O}(h^\infty )\vert u_-\vert &
\\ &+h^{1-k}\vert v_+\vert +{\cal O}(h^\infty )\Vert 1_Ku\Vert. &
\nonumber
\end{eqnarray}
According to (\ref{7.36}), we have
$$
\Vert (h+d^2)^{1-k}(1-\Pi )\chi _0u\Vert \le {\cal O}(h^{1\over
2})\Vert 1_Ku\Vert ,
$$
and using this and
$$\Vert (h+d^2)^{-k}\Pi \chi _0v\Vert
\le C\Vert (h+d^2)^{-k}\chi _0v\Vert
 $$
 in (\ref{hd2}), we get
(\ref{ev2.13}). \fin

\remark \label{Remev2.3} Return to the case $h=1$ and choose $P_0$
as after the equation (\ref{ev2.5}). Since $P_0$ is elliptic on
$\Lambda _{\Phi _\epsilon ^0}$ for $0<\epsilon \le \epsilon _0$
with $\epsilon _0$ small enough, an easy deformation argument
shows that the spectrum of $P_0$ on $H_{\Phi _\epsilon ^0}$ is
\indep{} of $\epsilon $, and similarly for the generalized
eigenvectors. The aim of this remark is to show that there exists
a $\delta >0$ such that the generalized eigenvectors $e_j$ satisfy
\ekv{ev2.15} { e_j\in H_{\Phi _0^0-\delta \vert x\vert ^2.} }
Rather than using deformation arguments as elsewhere in this
paper, we shall employ the alternative method of \fop{}s with
complex phase, and more precisely we shall study the evolution
equation associated to $P_0$.

\par Let us first recall some elementary facts from complex symplectic
geometry (as in \cite{Sjo95} and further references given there):
On $\C_x^n\times \C_\xi ^n$, we have the complex symplectic
$(2,0)$-form $\sigma =\sum_1^n d\xi _j\wedge dx_j$ and the real
symplectic forms $\Re \sigma $, $-\Im \sigma $. If $t$ is a \vf{}
on $\C^{2n}$ of type (1,0), we let $\widehat{t}=t+\overline{t}$ be
the associated real \vf{}. Then if $f$ is a \hol{} \fu{}, we let
$H_f$ denote the Hamilton field (of type (1,0)) \wrt{} $\sigma $
and if $g$ is a real-valued $C^1$-\fu{}, we let $H^{\Re \sigma
}_g$, $H^{-\Im \sigma }_g$ denote the Hamilton field of $g$ \wrt{}
$\Re \sigma $ and $-\Im \sigma $ respectively. Then we have the
relations,
$$\widehat{H}_f=H_{\Re f}^{\Re \sigma }=-H_{\Im f}^{-\Im \sigma },
\quad \widehat{H}_{if}=-H_{\Im f}^{\Re \sigma }=-H_{\Re f}^{-\Im
\sigma }.$$

\par Using \fop{}s with quadratic phase in the complex domain,
we see that
if $0\le t\le t_0$, and $u_0\in H_{\Phi _0}$, with $\Phi_0=\Phi
_0^0$, then we can solve the heat equation
$$
{\partial \over \partial t}u(t,x)+P_0u(t,x)=0,\ u(0,x)=u_0(x),
$$
and the solution \op{} $e^{-tP_0}$ is \bdd{} $H_{\Phi _0}\to
H_{\Phi _t}$, where \ekv{ev2.151} { \Lambda _{\Phi _t}=\exp
(t\widehat{H}_{{1\over i}P_0})(\Lambda _{\Phi _0})=\exp (tH_{\Re
P_0}^{-\Im \sigma })(\Lambda _{\Phi _0}). } We further have the
eikonal equation for $\Phi (t,x)=\Phi _t(x)$: \ekv{ev2.152} {
{\partial \Phi \over \partial t}+\Re P_0(x,{2\over i}{\partial
\Phi \over
\partial x})=0,
} corresponding to the \mfld{}
$$
\tau ={\partial \Phi \over \partial t},\ \xi ={2\over i}{\partial
\Phi \over \partial x},
$$
in $\R_{t,\tau }^2\times \C_{x,\xi }^{2n}$, which is Lagrangian
for the symplectic form $d\tau \wedge dt-\Im \sigma $. (To get
(\ref{ev2.152}) at least formally, differentiate
$\norm{u(t,.)}^2_{H_{\Phi(t,.)}}$ with respect to $t$). Since $\Re
P_0$ is constant along the flow of $H_{\Re P_0}^{-\Im \sigma }$,
we know that ${{\Re P_0}_\vert}_{\Lambda _{\Phi _t}}\ge 0$, so
(\ref{ev2.152}) shows that ${\partial \Phi \over \partial t}\le 0$
and hence that \ekv{ev2.153} { t\mapsto \Phi _t(x)\hbox{ is
decreasing}. } Let $L_0\subset \Lambda _{\Phi _0}$ be the
subspace, defined by $\Re P_0=0$ and notice that $\Re P_0\sim {\rm
dist\,}(\cdot ,L_0)^2$ on $\Lambda _{\Phi _0}$. In general, if $f$
is a smooth \fu{} on an IR-\mfld{} $\Lambda $, and $\widetilde{f}$
denotes an almost \hol{} extension of $f$ to a \neigh{} of
$\Lambda $, then at the points where $df$ is real, the Hamilton
field $H_f^{\sigma _\Lambda }$ of $f$ \wrt{} the real symplectic
form $\sigma _\Lambda ={\sigma _\vert}_{\Lambda }$ is equal to
$\widehat{H}_{\widetilde{f}}$. Applying this to $f={1\over i}P_0$,
we get at the points of $L_0$:
$$\widehat{H}_{{1\over i}P_0}
=H_{\Im P_0}^{\sigma _{\Lambda _{\Phi _0}}}.$$ Let $L_t=\exp
(t\widehat{H}_{{1\over i}P_0})(L_0)$ (cf (\ref{ev2.151})). Since
$H_{\Im P_0}^{\sigma _{\Lambda _{\Phi _0}}}$ is transversal to
$L_0$ away from $0$, we deduce that $t\mapsto \Pi _x(L_t)$ moves
transversally to $\Pi _x(L_0)$ away from 0 and since by
(\ref{ev2.152}),
$$
{\partial \Phi \over \partial t}\sim -{\rm dist\,}{\rm \,}(x,\Pi
_x(L_t))^2,
$$
we conclude that for small $t$ \ekv{ev2.154} { \Phi _t(x)\le \Phi
_0(x)-{t^3\over C}\vert x\vert ^2. } If $e_j$ is an eigenvector of
$P_0$: $P_0e_j=\lambda _je_j$, $\lambda _j\in\C$, we first recall
that $e_j\in H_{\Phi _0+\epsilon \vert x\vert ^2}$ for every
$\epsilon >0$, and conclude that $e^{-tP_0}e_j\in H_{\Phi
_t+\epsilon \vert x\vert ^2}$ for every $\epsilon >0$. On the
other hand, $e^{-tP_0}e_j=e^{-t\lambda _j}e_j$, so
$e_j=e^{t\lambda _j}e^{-tP}e_j\in H_{\Phi _t+\epsilon \vert x\vert
^2}$ for every $\epsilon >0$. Taking $t>0$ small but fixed, we
then obtain (\ref{ev2.15}) from (\ref{ev2.154}). If $\lambda _j$
is a multiple \ev{}, we also have to take into account the
possible Jordan blocks in the action of $P_0$ on the corresponding
generalized eigenspace, but this only requires minor modifications
in the argument and we get (\ref{ev2.15}) in general.

\subsection{Estimate for the semi-global problem}
\label{Subev3}

\par We now consider the situation in Section \ref{localsmall}. $P$
is now an $h$-\pop{} acting in $H_{\Phi }=H_{\Phi _\eps}$, and we
define $R_+=R_+^h$, $R_-=R_-^h$ as in the preceding subsection. As
in Section \ref{localsmall}, $P_0$ is now the quadratic
approximation of $P$ at $(0,0)$ and we shall use the fact that
$\Phi _\epsilon =\Phi _\epsilon ^0$ for $\vert x\vert \le
\sqrt{Ah}$ for $A\gg 1$. Recall the estimate (\ref{ev2.13}) for
solutions to (\ref{ev2.6}) (for $P_0$ and with norms in $H_{\Phi
_\epsilon ^0})$).

\par Again, we want to replace the fixed cutoff $\chi _0$ in
(\ref{ev2.13}) by $\chi _0({x\over \sqrt{Ah}})$ and consider the
change of variables $x=\sqrt{Ah}\widetilde{x}$,
$hD_x=\sqrt{Ah}\widetilde{h}D_{\widetilde{x}}$,
$\widetilde{h}=1/A$.
$$
P_0(x,hD_x;h)={h\over
\widetilde{h}}P_0(\widetilde{x},\widetilde{h}D_{\widetilde{x}};h)
=:{h\over
\widetilde{h}}\widetilde{P}_0,
$$
and with $d=d(x)$, $\widetilde{d}=d(\widetilde{x})$:
$$
h+d^2={h\over \widetilde{h}}(\widetilde{h}+\widetilde{d}^2),\
e^{-2\Phi ^0(x)/h}=e^{-2\Phi ^0(\widetilde{x})/\widetilde{h}},
$$
and we relate our unknown \fu{}s by the unitary relation
\ekv{ev3.1} { u(x)=(Ah)^{-{n\over
2}}\widetilde{u}(\widetilde{x}),\ h^{n\over
2}u(x)=\widetilde{h}^{n\over 2}\widetilde{u}(\widetilde{x}). }

With these substitutions, the \pb{} (\ref{ev2.6}) becomes
\ekv{ev3.2} { {h\over
\widetilde{h}}(\widetilde{P}_0-\widetilde{h}z)\widetilde{u}
+\widetilde{R}_-u_-=\widetilde{v},\
\widetilde{R}_+\widetilde{u}=v_+, } and we can apply
(\ref{ev2.13}) to this new \pb{}. A straightforward calculation
gives
\begin{equation} \label{ev3.3}
\begin{split}
 \Vert (h+d^2)^{1-k}\chi
_0({x\over \sqrt{Ah}})u\Vert +h^{-k}\vert u_-\vert  \leq  \ \  &
C\bigg(\Vert
(h+d^2)^{-k}\chi _0({x\over \sqrt{Ah}})v\Vert \\
& +h^{1-k}\vert v_+\vert +{1\over \sqrt{A}}\Vert
(h+d^2)^{1-k}1_K({x\over \sqrt{Ah}})u\Vert \bigg).
\end{split} \end{equation}
 This estimate will
be applied with $k=1/2$.

\par We now return to the full \op{} $P$ (on the FBI-side) and the norms
and scalar products will now be \wrt{} $ e^{-2\Phi /h}L(dx)$,
$\Phi =\Phi _\epsilon $, with $\epsilon >0$ small and fixed.
Recall however that $\Phi =\Phi _0^\epsilon $ in $\vert x\vert \le
\sqrt{Ah}$. We consider the semi-global Grushin \pb{} \ekv{ev3.4}
{ (P-hz)u+R_-u_-=v,\ R_+u=v_+, } in some fixed \bdd{} open set
containing the (projections of the) critical points. For
simplicity, we assume that the critical set is reduced to a single
point, corresponding to $x=0$. Let $\chi \in C_0^\infty $ be equal
to 1 near $0$.

\par Notice that by Remark \ref{Remev2.3}, $e^{-\Phi /h}R_-u_-$ is
exponentially small away from any fixed \neigh{} of $x=0$. Apply
(\ref{7.44}) with $(P-hz)u=v-R_-u_-$: \ekv{ev3.5} { \Vert \Lambda
u\Vert ^2\le C'\Re (\chi v\vert u)+C\Vert \Lambda ^{-1}R_-u_-\Vert
\Vert \Lambda u\Vert +C(\chi _0^2({x\over \sqrt{Ah}})\Lambda
u\vert \Lambda u)+C\Vert (1-\chi )\Lambda u\Vert \Vert \Lambda
u\Vert . } Here $\Lambda$ was defined in (\ref{defLambdaa}). Using
Remark \ref{Remev2.3} it is easy to check that \ekv{ev3.6} { \Vert
\Lambda ^{-1}R_-u_-\Vert \le {C\over \sqrt{h}}\vert u_-\vert , }
and (\ref{ev3.5}) becomes \ekv{ev3.7} { \Vert \Lambda u\Vert ^2\le
C\bigg(\Vert \Lambda ^{-1}v\Vert \Vert \Lambda u\Vert +{1\over
\sqrt{h}}\vert u_-\vert \Vert \Lambda u\Vert +\Vert \Lambda \chi
_0({x\over \sqrt{Ah}})u\Vert ^2+\Vert (1-\chi )\Lambda u\Vert
\Vert \Lambda u\Vert \bigg). } Apply "$2ab\le \alpha a^2+\alpha
^{-1}b^2$" with suitable $\alpha $'s to the 1st, 2nd and the 4th
terms of the \rhs{} and bootstrap away the $\Vert \Lambda u\Vert
^2$ terms. After removing the squares, we get \ekv{ev3.8} { \Vert
\Lambda u\Vert \le C\bigg(\Vert \Lambda ^{-1}v\Vert +{1\over
\sqrt{h}}\vert u_-\vert +\Vert \Lambda \chi _0({x\over
\sqrt{Ah}})u\Vert +\Vert (1-\chi )\Lambda u\Vert \bigg). } Apply
(\ref{ev3.3}) (for the Grushin problem for $P_0$) with $k=1/2$:
\begin{equation*}\begin{split}
& \Vert \Lambda \chi _0({x\over \sqrt{Ah}})u\Vert +h^{-1/2}\vert
u_-\vert \\
&\ \ \  \le  C\bigg(\Vert \Lambda ^{-1}\chi _0({x\over
\sqrt{Ah}})v\Vert +\Vert \Lambda ^{-1}\chi _0({x\over
\sqrt{Ah}})(P-P_0)u\Vert +h^{1\over 2}\vert v_+\vert +{1\over
\sqrt{A}}\Vert \Lambda 1_K({x\over
\sqrt{Ah}})u\Vert \bigg)  \\
& \ \ \ \leq  C\bigg(\Vert \Lambda ^{-1}\chi _0({x\over
\sqrt{Ah}})v\Vert +C(A)h^{1\over 2}\Vert \Lambda u\Vert +h^{1\over
2}\vert v_+\vert +{1\over \sqrt{A}}\Vert \Lambda 1_K({x\over
\sqrt{Ah}})u\Vert \bigg),\end{split}\end{equation*} where we used
(\ref{7.46}), to get the last estimate. Use this estimate in
(\ref{ev3.8}) after adding $h^{-1/2}\vert u_-\vert $ to both
sides:
\begin{equation*}
\begin{split}
\Vert \Lambda u\Vert +h^{-{1\over 2}}\vert u_-\vert \leq C\bigg(
\Vert \Lambda ^{-1}v\Vert +C(A)h^{1\over 2}\Vert \Lambda u\Vert
+h^{1\over 2}\vert v_+\vert +{1\over \sqrt{A}}\Vert \Lambda
1_K({x\over \sqrt{Ah}})u\Vert +\Vert (1-\chi )\Lambda u\Vert
\bigg),
\end{split}\end{equation*}
and choosing first $A$ large enough and then $h>0$ small enough,
we get the basic apriori estimate for the \pb{} (\ref{ev3.4}):
\ekv{ev3.9} { \Vert \Lambda u\Vert +h^{-{1\over 2}}\vert u_-\vert
\le C\big(\Vert \Lambda ^{-1}v\Vert +h^{1\over 2}\vert v_+\vert
+\Vert (1-\chi )\Lambda u\Vert \big). }
\medskip

\subsection{The global Grushin \pb{}.}\label{Subev4}

Now let $P$ be as in Theorem \ref{main}. Applying the inverse
FBI-\tf{} we have the obvious analogue of the Grushin \pb{} and
for that \pb{}, we still have (\ref{ev3.9}) provided that we
define $\Lambda $ to be a suitable $h$-\pop{} whose symbol is
equivalent to $(h+\min (d^2,(Ahd)^{2/3}))^{1/2}$, and interpret
$\chi $ as a pseudodifferential cutoff. From (\ref{ev3.9}), we get
the weaker estimate \ekv{ev4.1}{ h\Vert \psi ^wu\Vert +\vert
u_-\vert \le C(\Vert v\Vert +h\vert v_+\vert +h^{5/6}\Vert (1-\chi
)\psi ^wu\Vert ), } analogous to (\ref{chismall4}). This leads to
\ekv{ev4.2} { h\Vert \psi ^wu\Vert +\vert u_-\vert \le C(\Vert
v\Vert +h\vert v_+\vert +h\Vert q^wu\Vert +{\cal O}(h^2)\Vert
u\Vert ), } which is analogous to (\ref{chismall5}). On the other
hand, we have (\ref{philarge3}), and as in Section
\ref{SectProofMain}, we finally get the global apriori estimate
(analogous to (\ref{finalsmall2})): \ekv{ev4.3} { h\Vert u\Vert
+\vert u_-\vert \le C(\Vert v\Vert +h\vert v_+\vert ). }

 We are therefore exactly in the situation
 of the beginning of this section and from Proposition
 \ref{Propev1.2} we get  that the Grushin problem is well-posed.

\subsection{Asymptotics for $E_{-+}$ and for the \ev{}s. }
\label{Subev5}

\par For simplicity, we continue to assume that ${\cal C}$ is
reduced to a single point, $(0,0)$. We may assume that the global
Grushin \pb{} for the original \op{} $P$, considered in the
preceding subsection, is of the form \ekv{ev5.1} {
(P-hz)u+R_-u_-=v,\ R_+u=v_+, } where $z$ varies in a fixed
\neigh{} of an \ev{} $\lambda _0\in{\C}$, of the quadratic
approximation $P_0$ (with $h=1$) of $P$ at $(0,0)$, and where
\ekv{ev5.2} { R_-u_-=\sum_{j=1}^{N_0}u_-(j)e_j^h(x),
R_+u(j)=(u\vert f_j^h(x)). } Here \ekv{ev5.3} {
e_j^h(x)=h^{-{n\over 4}}e_j({x\over \sqrt{h}}),
f_j^h(x)=h^{-{n\over 4}}f_j({x\over \sqrt{h}}), } and
$e_1,...,e_{N_0}$ form a basis for the generalized eigenspace
$E_{\lambda _0}$ of $P_0$, associated to $\lambda _0$. It is well
known that we may take $e_j$ of the form \ekv{ev5.4} {
e_j(x)=p_j(x)e^{i\Phi _0(x)}, } where $p_j$ is a polynomial and
$\Phi _0(x)$ is a complex quadratic form such that $\Lambda _{\Phi
_0}=\{ (x,\Phi '_0(x))\}$ is the stable outgoing \mfld{} $\Lambda
^0$ for the ${1\over i}H_{P_0}$-flow and (by Remark
\ref{Remev2.3}) we know that \ekv{ev5.5} { \Im \Phi _0\hbox{ is
positive definite.} } We may assume that the $f_j$ have an
analogous form: \ekv{ev5.6} {f_j(x)=q_j(x)e^{i\Psi _0(x)},} with
$q_j$ polynomial and $\Psi _0$ a quadratic form with $\Im \Psi _0$
positive definite.

\par Let $\Lambda _\pm$ be the stable outgoing ($+$) and incoming ($-$)
\mfld{}s through (0,0) for the ${1 \over i}H_p$-flow, where $p$ is
the principal symbol of $P$. Then $\Lambda _\pm$ are complex
Lagrangian \mfld{}s defined to infinite order at $(0,0)$ and
$\Lambda _+^0=T_{(0,0)}\Lambda _+$. Let $\kappa $ be a complex
\ctf{}: ${\rm neigh\,}((0,0);{\C}^{2n})\to {\rm
neigh\,}((0,0);{\C}^{2n})$, mapping $\{ \xi =0\}$ to $\Lambda _+$
and $\{ x=0\}$ to $\Lambda _-$. Let $U$ be a formal elliptic
\fop{} of order 0 quantizing $\kappa $, and consider
$$
U^{-1}PU:=\widetilde{P},
$$
whose symbol is well-defined mod ${\cal O}((x,\xi )^\infty
+h^\infty )$. The principal symbol $\widetilde{p}$ of
$\widetilde{P}$ then vanishes on $\{ x=0\}$ and on $\{ \xi =0\}$
and therefore takes the form \ekv{ev5.7} {
\widetilde{p}=\sum_{\vert \alpha \vert =\vert \beta \vert
=1}a_{\alpha ,\beta }(x,\xi )x^\alpha \xi ^\beta . }

\par Using for simplicity the classical quantization of symbols,
 we get
\ekv{ev5.8} { \widetilde{P}=\sum_{\vert \alpha \vert =\vert \beta
\vert =1}a_{\alpha ,\beta }(x,hD)x^\alpha (hD)^\beta +ha(x,hD;h),
} where $a$ is a classical symbol of order $0$. (We are now
working with formal Taylor series at $(x,\xi )=(0,0)$.)

\par Put
\ekv{ev5.9} {{\cal P}_{\rm hom}^m=\{\sum_{\vert \alpha \vert
=m}b_\alpha \big({x\over \sqrt{h}}\big)^\alpha \}.} Here $b_\alpha
$ will in general be \fu{}s of $h$. When they are not, we say that
$\sum_{\vert \alpha \vert =m}b_\alpha \big({x\over
\sqrt{h}}\big)^\alpha $ is homogeneous of order $0$ in $h$ (or
even \indep{} of $h$, with $x/\sqrt{h}$ viewed as \indep{}
variables). Then in the obvious way,
$$
({x\over \sqrt{h}})^\gamma (\sqrt{h}D)^\delta :{\cal P}_{\rm
hom}^m\to {\cal P}_{\rm hom}^{m+\gamma -\delta }
$$
is homogeneous of degree 0 in $h$.

\par Write
\ekv{ev5.10} { {1\over h}\widetilde{P}= \sum_{\vert \alpha \vert
=\vert \beta \vert =1}a_{\alpha ,\beta }(x,hD)({x\over
\sqrt{h}})^\alpha (\sqrt{h}D)^\beta +a(x,hD;h). } Write $a\sim
\sum_0^\infty h^ja_j$ and Taylor expand $a(x,hD;h)$ at $(0,0)$:
\ekv{ev5.11} { a(x,hD;h)=\sum_{j=0}^\infty \sum_{\gamma ,\delta
}h^{j+{\vert \delta \vert \over 2}+{\vert \gamma \vert \over
2}}{a_{j(\delta )}^{(\gamma )}(0,0)\over \gamma !\delta !}({x\over
\sqrt{h}})^\delta (\sqrt{h}D)^\gamma . } If $\vert \delta \vert
-\vert \gamma \vert =k\in{\bf Z}$, then $\vert \gamma \vert +\vert
\delta \vert =\vert k\vert +2\min (\vert \delta \vert ,\vert
\gamma \vert )$, so the general term in the last sum can be
written
$$
h^{j+{\vert k\vert \over 2}+\min (\vert \gamma \vert ,\vert \delta
\vert )} {a_{j(\delta )}^{(\gamma )}(0,0)\over \gamma !\delta
!}({x\over \sqrt{h}})^\delta (\sqrt{h}D)^\gamma .
$$
In conclusion the block matrix of
$$
a(x,hD;h):\bigoplus_0^\infty {\cal P}_{\rm hom}^m\to
\bigoplus_0^\infty {\cal P}_{\rm hom}^m ,
$$
is $(h^{{\vert j-k\vert \over 2}}A_{j,k})$, where
$A_{j,k}=\sum_{\nu =0}^\infty A_{j,k}^\nu h^\nu $, and
$A_{j,k}^\nu :{\cal P}_{\rm hom}^k \to {\cal P}_{\rm hom}^j$ is
homogeneous of degree 0. (We then say that $A_{j,k}$ is a
classical symbol of order $0$.)

\par The same discussion applies to $a_{\alpha ,\beta }(x,hD)$ and hence
also to $h^{-1}\widetilde{P}$, whose matrix is \ekv{ev5.12} {
(h^{{\vert j-k\vert \over 2}}P_{j,k}),\ P_{j,k}=\sum_{\nu
=0}^\infty P_{j,k}^\nu h^\nu ,} where $P_{j,k}^\nu :{\cal P}_{\rm
hom}^k \to {\cal P}_{\rm hom}^j$ is homogeneous of degree 0. The
leading part of $h^{-1}\widetilde{P}$ is given by
$$
{1\over h}\widetilde{P}_0:=\sum_{\vert \alpha \vert =\vert \beta
\vert =1}a_{\alpha ,\beta }(0,0)({x\over \sqrt{h}})^\alpha
(\sqrt{h}D)^\beta +a_0(0,0),
$$
in the following sense: $h^{-1}\widetilde{P}_0$ has a block
diagonal matrix in $\bigoplus_0^\infty {\cal P}_{\rm hom}^m$, and
$P_{j,j}^0$ is equal to the restriction of $h^{-1}\widetilde{P}_0$
to ${\cal P}_{\rm hom}^j$.

\par Now we shall exploit that the exponent $\vert j-k\vert /2$ in
(\ref{ev5.12}) is an integer precisely when $j$ and $k$ have the
same parity. We therefore introduce \ekv{ev5.13} { {\cal
F}_e=\bigoplus_0^\infty {\cal P}_{\rm hom}^{2k},\ {\cal
F}_o=\bigoplus_0^\infty  {\cal P}_{\rm hom}^{2k+1}. } Then
$h^{-1}\widetilde{P}:{\cal F}_e\oplus {\cal F}_o\to {\cal
F}_e\oplus {\cal F}_o$ has the block diagonal matrix \ekv{ev5.14}
{
\begin{pmatrix}P_{e,e}, &P_{e,o}\cr P_{o,e} &P_{o,o}\end{pmatrix},
} where $P_{e,e},P_{o,o},h^{-1/2}P_{e,o}, h^{-1/2}P_{o,e}$ are
classical symbols of order $0$.

\par The Grushin \pb{} for $\widetilde{P}$ that we obtain from
(\ref{ev5.1}) is \ekv{ev5.15} {
(\widetilde{P}-hz)u+\widetilde{R}_-u_-=v,\ \widetilde{R}_+u=v_+, }
with \ekv{ev5.16} { \widetilde{R}_-=U^{-1}R_-,\
\widetilde{R}_+=R_+U. } We want to decompose $\widetilde{R}_\pm$
into even and odd degrees.

\par Return to $P,P_0$ and notice that $[P_0,\iota ]=0$, where
 $\iota$ is
the involution $\iota(u)(x)=u(-x)$. Consequently, $E_{\lambda _0}$
is invariant under $\iota$ and splits into $E_{\lambda _0}^e\oplus
E_{\lambda _0}^o$, with $\iota =1$ on $E_{\lambda _0}^e$ and
$\iota =-1$ on $E_{\lambda _0}^o$. Let the corresponding
dimensions be $N_e$, $N_o$, so that $N_0=N_e+N_o$. We may assume
that $e_j$ is even for $1\le j\le N_e$ and odd for $N_e+1\le j\le
N_0$, and we may choose $f_j$ with the same properties. Then
$p_j(x),q_j(x)$ are even when $1\le j\le N_e$ and odd otherwise.

\par Now, write
\ekv{ev5.17} { e_j^h(x)=h^{-{n\over 4}-{m_j\over
2}}a_j(x;h)e^{i{\Phi _0(x)\over h}}, } where \ekv{ev5.18} {
a_j\sim\sum_{\nu =0}^\infty a_j^\nu (x)h^\nu ,\hbox{ and
}a_j(x)={\cal O}(\vert x\vert ^{(m_j-2\nu )_+}),} and actually,
$a_j(x;h)=h^{m_j/2}p_j(x/\sqrt{h})$, with $m_j=d^o p_j$. $m_j$ is
even when $1\le j\le N_e$, and odd otherwise. Assume to fix the
ideas that $j\le N_e$. Then \ekv{ev5.19} {
U^{-1}(e_j^h)=h^{-{n\over 4}-{m_j\over
2}}\widetilde{a}_j(x;h)e^{i{F(x)\over h}}, } where
$\widetilde{a}_j$ satisfies (\ref{ev5.18}). Moreover, \ekv{ev5.20}
{ F(x)={\cal O}(x^3), } since $\Lambda _{\Phi _0}$ is tangent to
$\Lambda _\Phi $ so that $\Lambda _F$ is tangent to $\{ \xi =0\}$.

\par Taylor expanding $\widetilde{a}_j$ and $e^{iF/h}=\sum_0^\infty
(iF(x))^k/(k!h^k)$, we see that \ekv{ev5.21} { h^{n\over
4}U^{-1}(e_j^h)\in\bigoplus_0^\infty  {\cal P}_{\rm hom}^m, } and
when $m$ is even, the component in ${\cal P}_{\rm hom}^m$ is a
classical symbol of order 0 (and the order tends to $-\infty $
like $-m/2$, when $m\to\infty $), while the component in ${\cal
P}_{\rm hom}^m$ is of order $h^{1/2}$, when $m$ is odd). The case
$j\ge N_e+1$ is treated similarly, and we conclude that
\ekv{ev5.22} { \widetilde{R}_-=\begin{pmatrix}\widetilde{R}^{ee}_-
&\widetilde{R}^{eo}_-\cr
\widetilde{R}^{oe}_-&\widetilde{R}^{oo}_-\end{pmatrix}
:{\C}^{N_e}\oplus{\C}^{N_o}\to{\cal F}_e\oplus{\cal F}_o, } where
$h^{n/4}\widetilde{R}^{ee}_-$, $h^{n/4}\widetilde{R}^{oo}_-$,
$h^{n/4-1/2}\widetilde{R}^{eo}_-$,
$h^{n/4-1/2}\widetilde{R}^{oe}_-$ are classical symbols of order
0.

\par Next, we do the same work with $\widetilde{R}_+$ and start from
\ekv{ev5.23} { \widetilde{R}_+u(j)=(u\vert \trans{U}(f_{j,h})). }
Possibly after a slight perturbation of $\Psi $, we may assume
that \ekv{ev5.24} { \trans{U}(f_{j,h})=h^{-{n\over
4}}h^{-{\widetilde{m}_j\over 2}}\widetilde{b}_j(x;h)e^{{i\over
h}G(x)}, } where $\widetilde{m}_j$, $\widetilde{b}_j$ have the
same properties as $m_j,\widetilde{a}_j$ above, and
$\widetilde{m}_j$ is even for $1\le j\le N_e$ and odd otherwise.
Moreover $\det G''(0)\ne 0$, and the scalar product in
(\ref{ev5.23}) should be computed as a formal stationary phase
integral. In doing so, we apply the complex Morse lemma (to
$\infty $ order at $x=0$) to reduce $G$ to a quadratic form. If
$\alpha $ is a formal \diffeo{} with $\alpha (0)=0$, and
$Au=\alpha ^*u=u\circ \alpha $, then $A:\bigoplus_0^\infty {\cal
P}_{\rm hom}^m\to \bigoplus_0^\infty {\cal P}_{\rm hom}^m$ has the
same block matrix structure as $h^{-1}\widetilde{P}$ in
(\ref{ev5.12}). From these facts, we get \ekv{ev5.25} {
\widetilde{R}_+=\begin{pmatrix}\widetilde{R}^{ee}_+
&\widetilde{R}^{eo}_+\cr
\widetilde{R}^{oe}_+&\widetilde{R}^{oo}_+\end{pmatrix}: {\cal
F}_e\oplus{\cal F}_o\to {\C}^{N_e}\oplus{\C}^{N_o}, } where
$h^{-n/4}\widetilde{R}^{ee}_+$, $h^{-n/4}\widetilde{R}^{oo}_+$,
$h^{-n/4-1/2}\widetilde{R}^{eo}_+$,
$h^{-n/4-1/2}\widetilde{R}^{oe}_+$ are classical symbols of order
0.

\par Consider the rescaled \pb{} which is equivalent to (\ref{ev5.15}):
\ekv{ev5.26} { ({1\over h}\widetilde{P}-z)u+h^{n\over
4}\widetilde{R}_-u_-=v,\ h^{-{n\over 4}}\widetilde{R}_+u=v_+,} or
in matrix form
$$
\widetilde{{\cal P}}(z)\begin{pmatrix}u\cr u_-\end{pmatrix}=
\begin{pmatrix}v\cr v_+\end{pmatrix}.
$$
Let
$$
{\cal
E}=\begin{pmatrix}\widetilde{E}&\widetilde{E}_+\cr\widetilde{E}_-
&\widetilde{E}_{-+}\end{pmatrix}: \bigoplus_0^\infty {\cal P}_{\rm
hom}^m\oplus{\C}^{N_0}\to \bigoplus_0^\infty {\cal P}_{\rm
hom}^m\oplus{\C}^{N_0},
$$
be the inverse. Decomposing
$$
\widetilde{{\cal P}}(z)=\begin{pmatrix}\widetilde{{\cal P}}^{ee}
&\widetilde{{\cal P}}^{eo}\cr \widetilde{{\cal P}}^{oe} &
\widetilde{{\cal P}}^{oo}\end{pmatrix}: ({\cal
F}_e\oplus{\C}^{N_e})\oplus ({\cal F}_o\oplus{\C}^{N_o})\to ({\cal
F}_e\oplus{\C}^{N_e})\oplus ({\cal F}_o\oplus{\C}^{N_o}),
$$ where $\widetilde{{\cal P}}^{ee}$, $\widetilde{{\cal P}}^{oo}$,
$h^{-{1\over 2}}\widetilde{{\cal P}}^{}$, $h^{-{1\over
2}}\widetilde{{\cal P}}^{}$ are classical symbols of order 0, we
get the same decomposition for ${\cal E}(z)$. In particular,
\ekv{ev5.27} {
\widetilde{E}_{-+}(z)=\begin{pmatrix}\widetilde{E}^{ee}_{-+}
&\widetilde{E}^{eo}_{-+}\cr
\widetilde{E}^{oe}_{-+}&\widetilde{E}^{oo}_{-+}\end{pmatrix}:
{\C}^{N_e}\oplus{\C}^{N_o}\to{\C}^{N_e}\oplus{\C}^{N_o} } has the
same structure. The determinant of this matrix is a classical
symbol of order 0. In fact,
$$
\begin{pmatrix}h^{1\over 2}&0\cr 0&1\end{pmatrix}\widetilde{E}_{-+}(z)
\begin{pmatrix}h^{-{1\over 2}}&0\cr
0&1\end{pmatrix}=\begin{pmatrix}\widetilde{E}_{-+}^{ee}&h^{1\over
2} \widetilde{E}_{-+}^{eo}\cr h^{-{1\over
2}}\widetilde{E}_{-+}^{oe}&\widetilde{E}_{-+}^{oo}\end{pmatrix}
$$
has the same determinant and is a classical symbol of order 0 of the
form
\ekv{ev5.22.5}
{
\det \widetilde{E}_{-+}(\lambda )\sim \det E_{-+}^0(\lambda )+hf_1(\lambda
)+h^2f_2(\lambda )+...,
} 
where $E_{-+}^0(\lambda )$ is the matrix given in (\ref{ev2.5}) (there
denoted without the superscript 0). In particular $E_{-+}^0(\lambda
)=(\lambda -\lambda _0)^{N_0}f(\lambda )$ with $f(\lambda )\ne 0$, in the
space of \hol{} \fu{}s in a \neigh{} of $\lambda _0$. This could also
be deduced from the well known formula
$$
N_0={\rm tr\,}{1\over 2\pi i}\int _{\gamma _0}(\lambda
-P_0)^{-1}d\lambda ={1\over 2\pi i}\int {{d\over d\lambda }\det
E_{-+}^0\over \det E_{-+}^0(\lambda )}d\lambda ,
$$
where $\gamma _0$ is a closed contour around $\lambda _0$. (We have of
course the similar formula for $P, E_{-+}$, permitting to identify the
zeros of $E_{-+}$ and the \ev{}s of $P$, counted with their multiplicities.)

\par 
Using Puiseux series for the partial sums in (\ref{ev5.22.5})  we
conclude that the \ev{}s of $h^{-1}\widetilde{P}$ close to
$\lambda _0$ have complete \asy{} expansions in powers of
$h^{1/N_0}$:
$$
\lambda (h)=\lambda _0+c_1h^{1/N_0}+c_2h^{2/N_0}+... .
$$

\par Finally, it is clear from the construction, that if
$$
{\cal E}=\begin{pmatrix}E&E_+\cr E_-&E_{-+}\end{pmatrix}
$$
is the inverse of the global problem for $P$, then modulo ${\cal
O}(h^\infty )$: \ekv{ev5.28} {
E_{-+}(z;h)=h\widetilde{E}_{-+}(z;h). } and hence the true \ev{}s
of $P$ have the same \asy{} expansions as above. This completes
the proof of Theorem \ref{asympt}.

\Section{The evolution \pb{}.}\label{Section4}
\setcounter{equation}{0}

\par Let $P$ be a closed densely defined un\bdd{} \op{} acting on
 a complex Hilbert
space ${\cal H}$. Assume that the spectrum of $P$ is contained in
\ekv{4.1} { \Re z\ge {1\over C}\langle \Im z\rangle ^\delta -C, }
for some constants $C,\delta >0$. Assume also that \ekv{4.2} {
\Vert (z-P)^{-1}\Vert \le {C\over \langle z\rangle ^\delta
},\hbox{ for }\Re z\le {1\over 2C}\langle \Im z\rangle ^\delta
-2C. } For $t>0$ we put \ekv{4.3} { E(t)={1\over 2\pi
i}\int_\gamma  e^{-tz}(z-P)^{-1} dz, } where $\gamma $ is a
contour to the left of the spectrum which outside a compact set
coincides with the curve \ekv{4.4} { \Re z={1\over 3C}\langle \Im
z\rangle ^\delta , } and oriented in the direction of decreasing
$\Im z$. Clearly the integral converges and defines a \bdd{} \op{}
which depends smoothly on $t$. We have \ekv{4.5} { (\partial
_t+P)E(t)=0,\quad PE(t)=E(t)P. }

\par When $u\in{\cal D}(P)$ (the domain of $P$) we also have
\ekv{4.6} { \lim_{t\to 0}E(t)u=u. } In fact, let $z_0$ be to the
left of $\gamma $ and write \ekv{4.7} { u=(z_0-P)^{-1}v,\ v\in
{\cal H}. } The resolvent identity gives
$$
(z-P)^{-1}(z_0-P)^{-1}={1\over (z-z_0)}((z_0-P)^{-1}-(z-P^{-1}),
$$
so for $t>0$, we have \ekv{4.8} { E(t)u={1\over 2\pi i}\int_\gamma
e^{-tz}{1\over (z-z_0)}(z_0-P)^{-1}vdz-{1\over 2\pi i} \int_\gamma
e^{-tz}{1\over (z-z_0)}(z-P)^{-1}v dz. } Here the first integral
vanishes since we can push the contour to the right and exploit
the decay of the exponential. The second integral allows a limit
when $t\to 0$, so we get \ekv{4.9} { \lim_{t\to 0}E(t)u=-{1\over
2\pi i}\int_\gamma  {1\over (z-z_0)}(z-P)^{-1}vdz. } Here the
integrand is of norm ${\cal O}(\langle z\rangle ^{-1-\delta })$ in
view of (\ref{4.2}) and we can push the contour to the left
(around $z_0$) and apply the residue theorem to get
$$
\lim_{t\to 0}E(t)u=(z_0-P)^{-1}v=u,
$$
and (\ref{4.6}) follows.

\par In the following we assume that $P$ satisfies the assumptions
\bf (H1)--(H5) \rm   so that  Theorem \ref{main} gives  a
localization of the spectrum to a union of a conic \neigh{} of the
open positive axis and a infinite cusp away from the origin. We
introduce 2 contours $\gamma $ and $\widetilde{\gamma }$. Both
contours are given by \ekv{4.10} { \Re z={1\over C_0}h^{2\over
3}\vert \Im z\vert ^{1\over 3} } in the region $\Re z > bh$. Here
$C_0$ and $b$ are positive constants such that $b$ is different
from the real parts of the \ev{}s of the quadratic approximations
of $P$ with $h=1$. In the region $\Re z \le bh$, $\gamma $ is
given by $\Re z= bh$ while $\widetilde{\gamma }$ joins $bh+iC_0^3
b^3 h$ to $bh-iC_0^3 b^3 h$ further to the left so that
$\widetilde{\gamma }$ is entirely to the left of the spectrum of
$P$ while $\gamma $ will have a fixed finite number of \ev{}s,
$\lambda _0,...,\lambda _{N-1}$ to its left. Let $\gamma _{{\rm
int}}$ denote the vertical part of $\gamma $ in the region $\Re
z=bh$ and let $\gamma _{\rm ext}$ denote the part of $\gamma $ in
the region $\Re z \ge bh$.

\par On the exterior piece we have
\ekv{4.11} { \Vert (z-P)^{-1}\Vert \le {{\cal O}(1)\over
h^{2/3}\vert \Im z\vert ^{1/3}} ,} and on the interior piece we
have \ekv{4.12} { \Vert (z-P)^{-1}\Vert \le {{\cal O}(1)\over
h}. } This holds since we have chosen $b$ so
that the distance from $\gamma _{{\rm int}}$ to the spectrum of
$P$ is $\ge h/C$. Further to the left in the region $\Re
z \le bh$, we also have $\Vert (z-P)^{-1}\Vert ={\cal O}(h)$
when ${\rm dist\,}(z,\{ \lambda _0,..,\lambda _{N-1}\})\ge h/C$.

\par Assume for simplicity that  the eigenvalues of the different
quadratic approximations are simple and distinct. Then
\ekv{4.13} { e^{-tP/h}={1\over 2\pi i}\int_{\widetilde{\gamma
}}e^{-tz/h}(z-P)^{-1}dz=\sum _0^{N-1}e^{-t\lambda _j/h}\Pi
_{\lambda _j}+ {1\over 2\pi i}\int_{\gamma }e^{-tz/h}(z-P)^{-1}dz.
} Here $\Pi _{\lambda _j}$ is the (rank one) spectral projection
associated to $\lambda _j$, since the distance from $\lambda _j$ to
the other \ev{}s is $\ge h/C$, $\Pi _{\lambda _j}$ is \ufly{}
\bdd{} in norm when $h\to 0$.

\remark If we drop the assumption on the eigenvalues of the
quadratic approximations, then for instance two eigenvalues
$\lambda _1$ and $\lambda _2$ can be very close together but
separated from the others by $h/C$, and (since we are dealing with
a non-\sa{} \op{}) we can not state that $\Pi _1$ and $\Pi _2$ are
\ufly{} \bdd{} when $h\to 0$. But the sum $\Pi _1+\Pi _2$ will
have this property and so will the term $e^{-\lambda _1t}\Pi
_1+e^{-\lambda _2t}\Pi _2$ in (\ref{4.13}). This kind of situation
will appear when there is a symmetry, and to have a more complete
understanding in that case would include problems about the tunnel
effect.

\medskip

\par We estimate the last integral in (\ref{4.13}) using the
decomposition
$\gamma =\gamma _{{\rm int}}\cup\gamma _{{\rm ext}}$:
$${1\over 2\pi i}\int_{\gamma _{{\rm int}}}e^{-tz/h}(z-P)^{-1}dz=
{\cal O}(h)e^{-{t\over h}bh}{1\over h}={\cal
O}(1)e^{-bt},$$
\begin{align*}
{1\over 2\pi i}\int_{\gamma _{{\rm ext}}}e^{-tz/h}(z-P)^{-1}dz={}&
{\cal O}(1)\int_{C_0^3b^3h}^\infty e^{-{t\over C_0h}h^{2\over 3}
y^{1\over
3}}{1\over h^{2\over 3}y^{1\over 3}}dy\\ ={}&{{\cal
O}(1)\over t^2}\int_{tb}^\infty e^{-x}xdx\\
={}&{\cal O}(1)({1\over t}+{1\over t^2})e^{-tb}.
\end{align*}
Here and below, we let the prefactors ${\cal O}(1)$ depend on
$b,\, C_0$ Combining this with (\ref{4.13}), we get \ekv{4.14} {
e^{-tP/h}=\sum_0^{N-1}e^{-t\lambda _j/h}\Pi _{\lambda _j}+{\cal
O}(1)(1+{1\over t}+{1\over t^2})e^{-tb}. } It is quite possible
that the last estimate improves for small $t$ when we let
$e^{-tP/h}$ act on elements in the domain of $P$. Since $P$ is
accretive with $\Re P \geq -Ch$ we can get rid of the terms
  $1/t$ and $1/t^2$. The proof of Theorem \ref{evolution}
 is complete.

\Section{Application to the Kramers-Fokker-Planck operator}
\label{secFP}

In this section we prove Theorem \ref{mainFP} and compute the
eigenvalues of the Kramers-Fokker-Planck operator with quadratic
potential $V$~:
\begin{equation} \label{defFP5}
 P =v\cdot h\partial _x-V'(x)\cdot h\partial _v+ {\gamma \over
2}(-(h\partial _v)^2+v^2-hn) .
\end{equation}
We will recall the classical procedure to obtain this \op{} by
conjugation from
\begin{equation} \label{eq:P_FP}
 P_{\rm FP}=v\cdot h\partial _x-V'(x)\cdot
h\partial _v -{1\over 2}\gamma h\partial _v\cdot (h\partial
_v+2v).
\end{equation}
In this article $V$ is a $\cc^\infty$ potential with bounded
derivatives of second and higher order, and $x$, $v \in \R^n$.
We suppose that $V$ has a finite number of critical points.

 We observe (see for example
\cite{HerNi03}) that the first two terms form the
Hamilton field $X_0$ of the Hamiltonian $q$ where
$$
q(x,v)={1\over 2}v^2+V(x), \ \ \ \ X_0 = v\cdot h\partial
_x-V'(x)\cdot h\partial _v,
$$
when $v$ is considered as the dual variables of $x$. The Maxwellian is
defined by
$$
  M=e^{-{2\over h}({v^2\over 2}+V(x))} ,
$$
and we get the formally conjugated \op{} $P =
M^{1/2}P_{\rm FP}M^{-1/2}$ in (\ref{defFP5}).

\subsection{Metrics and hypotheses}

In this section we check that the Kramers-Fokker-Planck operator
 satisfies
the hypotheses of the main theorem under the simple assumptions
that $V$ is a Morse function with bounded derivatives of order 2
and higher, such that $|V'(x)| \geq 1/C$ when $|x| \geq C$. Denote
by $(\xi, \eta)$ the variable dual to $(x,v)$. Then
$$
P =p^w - \frac{\gamma h n}{2}, \ \ \ \text{with} \ \  p = {\gamma
\over 2} (v^2 + \eta^2) + iv \cdot \xi -iV'(x) \cdot \eta.
$$
(Note that for this operator $p^w = p(x,hD_x)$). We introduce the
natural weight associated to $P$
$$
\lambda^2(x,\xi,v,\eta) = 1+ ( V'(x))^2 + \xi^2 + v^2 + \eta^2,
$$
and the metric
$$
\Gamma_0 = dx^2 + dv^2 + \frac{ d\xi^2 + d\eta^2}{\lambda^2}.
$$
We note that $\lambda$ is $\cc^\infty$ and that
$$
\lambda \in S(\lambda, \Gamma_0), \ \ \ \lambda' \in S(1,
\Gamma_0) ,
$$
since $V$ is with second derivative bounded. Let us now check that
 $p$ satisfies the symbolic estimates (\ref{symbcal}). Denoting
 $p=p_1 + ip_2$ and $\rho = (x,v,\xi,\eta)$ we get
\begin{equation}
\begin{split}
p_1(\rho) & = {\gamma \over 2} (v^2 + \eta^2) , \\
p_2(\rho) & = v \cdot \xi -V'(x)\cdot \eta , \\
\D p_1(\rho) & = \gamma (0,v,0,\eta) , \\
\D p_2(\rho) & = \gamma ( -V''(x) \cdot \eta,\xi, v, -V'(x)) , \\
\D^2 p_1(\rho) & = \gamma \left(
    \begin{array}{cccc}
    0 &  0 & 0 &0 \\
    0 &  Id & 0 & 0 \\
    0&0&0&0 \\
    0&0&0&  Id
    \end{array} \right), \\
 H_{p_2}(\rho) & =  v \cdot \partial_x - V'(x) \cdot \partial_v
    + \eta \cdot V'' \cdot \partial_\xi - \xi \cdot \partial_\eta , \\
H_{p_2} p_1(\rho) & =  -\gamma \xi \cdot \eta - \gamma V'(x)\cdot v ,\\
\D H_{p_2} p_1(\rho) & = ( - \gamma V''(x)\cdot v, - \gamma V'(x),
-\gamma \eta, -\gamma \xi),  \\
  H_{p_2}^2 p_1(\rho) & =
    -\gamma  v \cdot  V''(x) \cdot v + \gamma (V'(x))^2
    - \gamma \eta \cdot V''(x) \cdot \eta + \gamma \xi^2 .
\end{split}
\end{equation}
We get directly, using that the derivatives of $V$ of order 2
 and higher
are bounded, that
\begin{equation} \label{symbcal3}
\begin{split}
  p_1 & \geq 0, \ \ \
  p  \in S(\lambda^{2},\Gamma_0), \ \ \
  \D p  \in S(\lambda^{},\Gamma_0), \ \ \
  \D^2 p_1  \in S(1,\Gamma_0), \ \ \
 \D H_{p_2} p_1   \in S(\lambda,\Gamma_0).
\end{split}
\end{equation}
Besides, let us denote by $\set{\rho_j}$ the critical points of
$p$. We notice that they are of the form $(x_j,0,0,0)$ where
$\set{x_j}$ are the critical points of $V$. By $\delta^2$ we
denote a $\cc^\infty$ function equivalent to the distance to the
set $\set{\rho_j}$. Then for $\eps_0$ sufficiently small we have
\begin{equation} \label{symbcal4}
\begin{split}
p_1 + \eps_0 H_{p_2}^2 p_1
    & = \gamma \sep{^{} \eps_0(V'(x))^2
    + \eps_0  \xi^2 +  v \cdot (Id/2 - \eps_0 V''(x)) \cdot v
    + \eta  \cdot (Id/2 - \eps_0 V''(x)) \cdot \eta} \\
    & \sim \left\{ \begin{array}{ll}
      \delta^2    & \text{ in a fixed compact set including
      the $\rho_j$s}, \\
      \lambda^2   & \text{ away from a \neigh{} of the $\rho_j$s.}
       \end{array}\right.
\end{split}
\end{equation}
The last thing to check is that the metric $\Gamma_0$ is
(classically in the sense of (\ref{admissible})) admissible.
A simple adaptation of Proposition 5.11 in \cite{HelNi}
shows that $\Gamma_0$ is $cl$-admissible. Note that it is therefore
semiclassically  admissible (in the sense of (\ref{admissible2}))
since in that case weaker assumptions are needed.

As a consequence we can apply to the Kramers-Fokker-Planck
operator $P=p^w
- \gamma h n/2$ the main Theorem \ref{main}. In order to be
complete we compute now the eigenvalues of the quadratic
approximation $P_0$ of $P$ near the critical points.

\subsection{Eigenvalue computation}

Here we will compute explicitly the $\lambda_j$ that
occur in the formula for the spectrum given in
Proposition~\ref{quad}. In addition we compute the constant term that
also contributes to the eigenvalues. Thus we obtain the spectrum up to
$o(h)$. We assume that $p$ has a single critical point at $x=0$ and
that $V$ is quadratic. After a simultaneous orthogonal change of
coordinates in $x$ and in $v$, we may assume that
\begin{equation} \label{eq:simplify_V}
  V(x) = \tfrac{1}{2} \sum_{j=1}^n d_j x_j^2 .
\end{equation}
(The assumption that $V$ is a Morse \fu{} implies that all the $d_j$
are different from 0.)

With this choice of $V$, the operator $P_{\rm FP}$ equals
\begin{equation} \label{eq:P_symmetrization}
  P_{\rm FP} = - \tfrac{1}{2} \gamma h n + (x,v,D_x,D_v)
  W (x,v,D_x,D_v)^T .
\end{equation}
where the matrix $W$ is given by
\[
  W = \left(\begin{matrix}
    0 & 0 & 0              & -\frac{i}{2} h V''_{xx} \\
    0 & 0 & \frac{i}{2}h I & -\frac{i}{2} \gamma h I   \\
    0                       & \frac{i}{2}h I        & 0 & 0 \\
    -\frac{i}{2} h V''_{xx} & -\frac{i}{2} \gamma h  I & 0 &
                                         \frac{1}{2} h^2 \gamma I
  \end{matrix} \right) .
\]
As we explained, the operators $P$ and $P_0$ of Proposition~\ref{quad}
are obtained by conjugation, which corresponds to a complex symplectic
coordinate transformation of the symbol. Since the eigenvalues of the
linearization of the Hamilton flow are invariant under such a
transformation, we can use the unconjugated operator $P_{\rm FP}$ to
compute them.

The matrix $W$ is of the form
$W = \half \left( \begin{matrix} 0 & i h A^T \\ i h A & h^2 B
\end{matrix} \right)$, where the $2n \times 2n$ matrix
$A = \left(\begin{matrix} 0 & I \\ -V''_{xx} & -\gamma I\end{matrix}
\right)$
is the linearization of the vector field component of $P_{\rm FP}$
for $h=1$.
The matrix corresponding to the linearization of the Hamilton field is
given by
\begin{equation}
  W' = \left( \begin{matrix} i h A & h^2 B\\
            0 & - i h A^T \end{matrix} \right) .
\end{equation}
Because of (\ref{eq:simplify_V}), the eigenvalues of $A$ are obtained
simply by diagonalizing the $2 \times 2$ matrices
$\left(\begin{matrix} 0 & 1 \\ -d_j & - \gamma \end{matrix}\right)$.
We find that the eigenvalues of $A$ are given by
\[
  \nu_{j,1}
    = - \frac{\gamma}{2} - \frac{1}{2} \sqrt{\gamma^2 - 4 d_j} ,
     \qquad
  \nu_{j,2}
    = - \frac{\gamma}{2} + \frac{1}{2} \sqrt{\gamma^2 - 4 d_j} ,
\]
with $j=1,\ldots, n$.
Let $s_{j,k}$ denote the sign of the real part of $\nu_{j,k}$
\[
  s_{j,1} = \sgn(\Re(\nu_{j,1})) = -1, \ \ \
  s_{j,2} = \sgn(\Re(\nu_{j,2})) = -\sgn(d_j) .
\]
It follows that the eigenvalues of $W'$ are given by
\[
  i h \nu_{j,1}, \ \ \ i h \nu_{j,2} , \ \ \
  - i h \nu_{j,1}, \ \ \ -i h \nu_{j,2} ,
\]
and that the ones with positive imaginary part are given by
\[
\begin{split}
  i h s_{j,1} \nu_{j,1}
  = & \frac{i}{2} \gamma h + \frac{i}{2} h \sqrt{\gamma^2 - 4 d_j} ,
  \\
  i h s_{j,2} \nu_{j,2}
  = & - \sgn(d_j) \left( - \frac{i}{2} \gamma h
    + \frac{i}{2} h \sqrt{\gamma^2 - 4 d_j} \right) .
\end{split}
\]
The constant term in (\ref{eq:P_symmetrization}) satisfies
$- \tfrac{1}{2} \gamma h n = \tfrac{1}{2} \tr (A)
= \frac{1}{2} \sum_{j=1}^n (\nu_{j,1} + \nu_{j,2})$.
Thus the spectrum of the quadratic operator is given by
\[
  \left\{ h \sum_{j=1}^n
     \left( (\tfrac{1}{2} + \tfrac{1}{2} s_{j,1} + k_{j,1} )
     s_{j,1} \nu_{j,1}
   + (\tfrac{1}{2} + \tfrac{1}{2} s_{j,2} + k_{j,2} )
     s_{j,2} \nu_{j,2} \right) \, ; \;\,
  k_{j,1},k_{j,2} \in \N \right\} .
\]

\remark In the case of quadratic Kramers-Fokker-Planck, the lowest
eigenvalue of the spectrum is $0$ if and only if all the
$s_{j,1},s_{j,2}$ are equal to $-1$, i.e.\ if $x=0$ is a minimum
of $V$.

\medskip


\end{document}